\documentclass[english,toc=bib]{scrartcl}

\usepackage[T1]{fontenc}
\usepackage[utf8]{inputenc}
\usepackage[english]{babel}
\usepackage{xspace}
\usepackage{lmodern}
\usepackage{exscale}
\usepackage[babel]{microtype}
\usepackage{amsmath, amssymb, mathtools, mathrsfs}
\usepackage{graphicx}
\usepackage{xcolor}
\usepackage[xspace]{ellipsis}
\usepackage{array,booktabs}
\usepackage{enumitem}

\usepackage{relsize}
\newcommand{\abbr}[1]{\textsmaller{#1}}

\usepackage{tikz}
\usetikzlibrary{positioning,cd}

\usepackage[style=alphabetic, sorting=nyt, maxnames=10, maxalphanames=4]{biblatex}

\usepackage[numbered]{bookmark}
\usepackage{hyperref}
\hypersetup{
  colorlinks,
  urlcolor = {blue!80!black},
  linkcolor = {red!70!black},
  citecolor = {green!50!black}
}

\usepackage{amsthm,thmtools}
\newtheorem{thrm}[equation]{Theorem}
\newtheorem{prop}[equation]{Proposition}
\newtheorem{lem}[equation]{Lemma}
\newtheorem{coro}[equation]{Corollary}
\theoremstyle{definition}
\newtheorem{defin}[equation]{Definition}
\theoremstyle{remark}

\newtheorem{rmk}[equation]{Remark}
\newtheorem{conj}[equation]{Conjecture}
\makeatletter
\let\c@figure\c@equation
\makeatother

\DeclareMathOperator{\im}{im}
\DeclareMathOperator{\id}{id}
\DeclareMathOperator{\Hom}{Hom}
\DeclareMathOperator{\ev}{ev}

\newcommand{\Q}{\mathbb{Q}}
\newcommand{\R}{\mathbb{R}}
\newcommand{\Z}{\mathbb{Z}}
\newcommand{\N}{\mathbb{N}}
\newcommand{\PP}{\mathtt{P}}
\newcommand{\MM}{\mathtt{M}}
\newcommand{\CC}{\mathtt{C}}
\newcommand{\DD}{\mathtt{D}}
\newcommand{\NN}{\mathtt{N}}
\DeclareMathOperator{\Tw}{Tw}
\newcommand{\vol}{\mathrm{vol}}
\newcommand{\Conf}{\mathrm{Conf}}
\newcommand{\FM}{\mathtt{FM}}
\newcommand{\OmPA}{\Omega_{\mathrm{PA}}}
\newcommand{\APL}{A_{\mathrm{PL}}}
\newcommand{\Gra}{\mathtt{Gra}}
\newcommand{\Graphs}{\mathtt{Graphs}}
\newcommand{\enV}[1][n]{{\mathtt{e}_{#1}^{\vee}}}
\newcommand{\Lie}{\mathtt{Lie}}
\newcommand{\hoLie}{\mathtt{hoLie}}
\newcommand{\Com}{\mathtt{Com}}
\newcommand{\fGC}{\mathrm{fGC}}
\newcommand{\fLoop}{\mathrm{fLoop}}
\newcommand{\EE}{\mathsf{E}}
\newcommand{\GG}[1]{\mathtt{G}_{#1}}
\newcommand{\Embfr}{\operatorname{Emb}^{\mathrm{fr}}}
\newcommand{\Disk}{\mathtt{Disk}}
\newcommand{\hadprod}{\boxtimes}
\newcommand{\CCE}{\mathrm{C}_{*}^{\mathrm{CE}}}
\newcommand{\Zphi}{Z_{\varphi}}
\newcommand{\Ze}{Z_{\varepsilon}}

\tikzset{extv/.style={circle, draw, minimum size=0.5cm, inner sep=0pt}}
\tikzset{intv/.style={circle, fill = black, draw, minimum size = 0.3cm, inner sep = 0pt}}
\tikzset{unkv/.style={circle, fill = gray, draw, minimum size = 0.3cm, inner sep = 0pt}}

\begin{document}

\title{The Lambrechts--Stanley Model of Configuration Spaces}
\author{Najib Idrissi\thanks{Université Paris Diderot, Institut de Mathématiques de Jussieu-Paris Rive Gauche, Sorbonne Paris Cité, CNRS, Sorbonne Université, F-75013 Paris, France. \href{mailto:najib.idrissi-kaitouni@imj-prg.fr}{\nolinkurl{najib.idrissi-kaitouni@imj-prg.fr}}}}
\date{November 23, 2018}

\maketitle

\begin{abstract}
  We prove the validity over $\R$ of a commutative differential graded algebra model of configuration spaces for simply connected closed smooth manifolds, answering a conjecture of Lambrechts--Stanley.
  We get as a result that the real homotopy type of such configuration spaces only depends on the real homotopy type of the manifold.
  We moreover prove, if the dimension of the manifold is at least $4$, that our model is compatible with the action of the Fulton--MacPherson operad (weakly equivalent to the little disks operad) when the manifold is framed.
  We use this more precise result to get a complex computing factorization homology of framed manifolds.
  Our proofs use the same ideas as Kontsevich's proof of the formality of the little disks operads.
\end{abstract}

\tableofcontents

\section*{Introduction}
\label{sec.introduction}

Let $M$ be a closed smooth $n$-manifold and consider the ordered configuration space of $k$ points in $M$:
\[ \Conf_{k}(M) \coloneqq \{ (x_{1}, \dots, x_{k}) \in M^{k} \mid x_{i} \neq x_{j} \; \forall i \neq j \}. \]

Despite their apparent simplicity, configuration spaces remain intriguing.
One of the most basic questions that can be asked about them is the following: if a manifold $M'$ is obtained from $M$ by continuous deformations, then can $\Conf_{k}(M')$ be obtained from $\Conf_{k}(M)$ by continuous deformations?
That is, does the homotopy type of $M$ determine the homotopy type of $\Conf_{k}(M)$?

Without any restriction, this is false: the point $\{0\}$ is homotopy equivalent to the line $\R$, but $\Conf_{2}(\{0\}) = \varnothing$ is not homotopy equivalent to $\Conf_{2}(\R) \neq \varnothing$.
One might wonder if the conjecture becomes true if restricted to closed manifolds.
In 2005, Longoni and Salvatore~\cite{LongoniSalvatore2005} found a counterexample: two closed $3$-manifolds, given by lens spaces, which are homotopy equivalent but whose configuration spaces are not.
This counterexample is not simply connected however.
The question of the homotopy invariance of $\Conf_{k}(-)$ for simply connected closed manifolds remains open to this day.

Here, we do not work with the full homotopy type.
Rather, we restrict ourselves to the \emph{rational} homotopy type.
This amounts, in a sense, to forgetting all the torsion.
Rational homotopy theory can be studied from an algebraic point of view~\cite{Sullivan1977}.
The rational homotopy type of a simply connected space $X$ is fully encoded in a ``model'' of $X$, i.e.\ a commutative differential graded algebra (\abbr{CDGA}) $A$ which is quasi-isomorphic to the \abbr{CDGA} of piecewise polynomial forms $\APL^{*}(X)$.
Due to technical issues, we will in fact work over $\R$.
If $M$ is a smooth manifold, then a real model is a \abbr{CDGA} which is quasi-isomorphic to the \abbr{CDGA} of de Rham forms $\Omega^{*}_{\mathrm{dR}}(M)$.
While this is slightly coarser than the rational homotopy type of $M$, in terms of computations it is often enough.

Thus, our goal is the following: given a model of $M$, deduce an explicit, small model of $\Conf_{k}(M)$.
This explicit model only depends on the model of $M$.
This shows the (real) homotopy invariance of $\Conf_{k}(-)$ on the class of manifolds we consider.
Moreover, this explicit model can be used to perform computations, e.g.\ the real cohomology ring of $\Conf_{k}(M)$, etc.

We focus on simply connected (thus orientable) closed manifolds.
They satisfy Poincaré duality.
Lambrechts and Stanley~\cite{LambrechtsStanley2008} showed that any such manifold admits a model $A$ which satisfies itself Poincaré duality, i.e.\ there is an ``orientation'' $A^{n} \xrightarrow{\varepsilon} \R$ which induces non-degenerate pairings $A^{k} \otimes A^{n-k} \to \R$ for all $k$.
Lambrechts and Stanley~\cite{LambrechtsStanley2008a} built a \abbr{CDGA} $\GG{A}(k)$ out of such a Poincaré duality model (they denote it $F(A,k)$).
If we view $H^{*}(\Conf_{k}(\R^{n}))$ as spanned by graphs modulo Arnold relations, then $\GG{A}(k)$ consists of similar graphs with connected components labeled by $A$, and the differential splits edges.
Lambrechts and Stanley proved that $\GG{A}(k)$ is quasi-isomorphic to $\APL^{*}(\Conf_{k}(M))$ as a dg-module.
They conjectured that this quasi-isomorphism can be enhanced to give a quasi-isomorphism of \abbr{CDGA}s so that $\GG{A}(k)$ defines a rational model of $\Conf_{k}(M)$.
We answer this conjecture by the affirmative in the real setting in the following theorem.

\begin{thrm}[{Corollary~\ref{cor:main-cor}}]
  \label{thm:intermediary}%
  Let $M$ be a simply connected, closed, smooth manifold.
  Let $A$ be any Poincaré duality model of $M$.
  Then for all $k \ge 0$, $\GG{A}(k)$ is a model for the real homotopy type of $\Conf_{k}(M)$.
\end{thrm}

\begin{coro}[{Corollary~\ref{cor.only-depends}}]
  For simply connected closed smooth manifolds, the real homotopy type of $M$ determines the real homotopy type of $\Conf_{k}(M)$.
\end{coro}

Over the past decades, attempts were made to solve the Lambrechts--Stanley conjecture, and results were obtained for special kinds of manifolds, or for low values of $k$.
When $M$ is a smooth complex projective variety, Kriz~\cite{Kriz1994} had previously shown that $\GG{H^{*}(M)}(k)$ is actually a rational \abbr{CDGA} model for $\Conf_{k}(M)$.
The \abbr{CDGA} $\GG{H^{*}(M)}(k)$ is the $\mathsf{E}^{2}$ page of a spectral sequence of Cohen--Taylor~\cite{CohenTaylor1978} that converges to $H^{*}(\Conf_{k}(M))$.
Totaro~\cite{Totaro1996} has shown that for a smooth complex compact projective variety, the spectral sequence only has one nonzero differential.
When $k = 2$, then $\GG{A}(2)$ was known to be a model of $\Conf_{2}(M)$ either when $M$ is $2$-connected~\cite{LambrechtsStanley2004} or when $\dim M$ is even~\cite{CordovaBulens2015}.

Our approach is different than the ones used in these previous works.
We use ideas coming from the theory of operads.
In particular, we consider the operad of little $n$-disks, defined by Boardman--Vogt~\cite{BoardmanVogt1973}, which consists of configuration spaces of small $n$-disks (instead of points) embedded inside the unit $n$-disk.
These spaces of little $n$-disks are equipped with composition products, which are basically defined by inserting a configuration of $l$ little $n$-disks into the $i$th little disk of a configuration of $k$ little $n$-disks, resulting in a configuration of $k+l-1$ little $n$-disks.
The idea is that a configuration of little $n$-disks represents an operation acting on $n$-fold loop spaces, and the operadic composition products of little $n$-disks reflect the composition of such operations.
The configuration spaces of little $n$-disks are homotopy equivalent to the configurations spaces of points in the Euclidean $n$-space $\R^{n}$, but the operadic composition structure does not go through this homotopy equivalence.

In our work, we actually use another model of the little $n$-disk operads, defined using the Fulton--MacPherson compactifications $\FM_{n}(k)$ of the configurations spaces $\Conf_k(\R^n)$~\cite{FultonMacPherson1994,AxelrodSinger1994,Sinha2004}.
This compactification allows us to retrieve, on this collection of spaces $\FM_{n} = \{ \FM_{n}(k) \}$, the operadic composition products which were lost in the configurations spaces $\Conf_k(\R^n)$.
We also use the Fulton--MacPherson compactifications $\FM_{M}(k)$ of the configuration spaces $\Conf_{k}(M)$ associated to a closed manifold $M$.
When $M$ is framed, these compactifications assemble into an operadic right module $\FM_{M}$ over the Fulton--MacPherson operad $\FM_{n}$, which roughly means that we can insert a configuration in $\FM_{n}$ into a configuration in $\FM_{M}$.
We show that the Lambrechts--Stanley model is compatible with this action of the little disks operad, as we explain now.

The little $n$-disks operads are formal~\cite{Kontsevich1999,Tamarkin2003,Petersen2014,LambrechtsVolic2014,FresseWillwacher2015}.
Kontsevich's proof~\cite{Kontsevich1999,LambrechtsVolic2014} of this theorem uses the spaces $\FM_{n}$.
If we temporarily forget about operads, this formality theorem means in particular that each space $\FM_{n}(k)$ is ``formal'', i.e.\ the cohomology $\enV(k) \coloneqq H^{*}(\FM_{n}(k))$ (with a trivial differential) is a model for the real homotopy type of $\FM_{n}(k)$.
To prove Theorem~\ref{thm:intermediary}, we generalize Kontsevich's approach to prove that $\GG{A}(k)$ is a model of $\FM_{M}(k)$.

To establish his result, Kontsevich has to consider fiberwise integrations of forms along a particular class of maps, which are not submersions, but represent the projection map of ``semi-algebraic bundles''.
In order to define such fiberwise integration operations, Kontsevich uses \abbr{CDGA}s of piecewise semi-algebraic (\abbr{PA}) forms $\OmPA^{*}(-)$ instead of the classical \abbr{CDGA}s of de Rham forms.
The theory of \abbr{PA} forms was developed in~\cite{KontsevichSoibelman2000,HardtLambrechtsTurchinVolic2011}.
Any closed smooth manifold $M$ is a semi-algebraic manifold~\cite{Nash1952,Tognoli1973}, and the \abbr{CDGA} $\OmPA^{*}(M)$ is a model for the real homotopy type of $M$.
For the formality of $\FM_{n}$, a descent argument~\cite{GuillenNavarroPascualRoig2005} is available to show that formality over $\R$ implies formality over $\Q$.
However, no such descent argument exists for models with a nontrivial differential such as $\GG{A}$.
Therefore, although we conjecture that our results on real homotopy types descend to $\Q$, we have no general argument ensuring that such a property holds.

The cohomology $\enV = H^{*}(\FM_{n})$ inherits a Hopf cooperad structure from $\FM_{n}$, i.e.\ it is a cooperad (the dual notion of operad) in the category of \abbr{CDGA}s.
The \abbr{CDGA}s of forms $\OmPA^{*}(\FM_{n}(k))$ also inherit a Hopf cooperad structure (up to homotopy).
The formality quasi-isomorphisms between the cohomology algebras $\enV(k)$ and the \abbr{CDGA}s of forms on $\FM_{n}(k)$ are compatible in a suitable sense with this structure.
Therefore the Hopf cooperad $\enV$ fully encodes the rational homotopy type of the operad $\FM_{n}$.

In this paper, we also prove that the Lambrechts-Stanley model $\GG{A}$ determines the real homotopy type of $\FM_{M}$ as a right module over the operad $\FM_{n}$ when $M$ is a framed manifold.
To be precise, our result reads as follows.
\begin{thrm}[{Theorem~\ref{thm.Abis}}]
  \label{thm.A}
  Let $M$ be a framed smooth simply connected closed manifold with $\dim M \ge 4$.
  Let $A$ be any Poincaré duality model of $M$.
  Then the collection $\GG{A} = \{ \GG{A}(k) \}_{k \geq 0}$ forms a Hopf right $\enV$-comodule.
  Moreover the Hopf right comodule $(\GG{A}, \enV)$ is weakly equivalent to $(\OmPA^{*}(\FM_{M}), \OmPA^{*}(\FM_{n}))$.
\end{thrm}

For $\dim M \le 3$, the proof fails (see Proposition~\ref{prop.strong-vanishing}).
However, in this case, the only examples of simply connected closed manifolds are spheres, thanks to Perelman's proof of the Poincaré conjecture~\cite{Perelman2002,Perelman2003}.
We can then directly prove that $\GG{A}(k)$ is a model for $\Conf_{k}(M)$ (see Section~\ref{sec:models-conf-2}).

Our proof of Theorem~\ref{thm.A}, which is inspired by Kontsevich's proof of the formality of the little disks operads, is radically different from the proofs of~\cite{LambrechtsStanley2008a}.
It involves an intermediary Hopf right comodule of labeled graphs $\Graphs_{R}$.
This comodule is similar to a comodule recently studied by Campos--Willwacher~\cite{CamposWillwacher2016}, which corresponds to the case $R = S(\tilde{H}^{*}(M))$.
Note however that the approach of Campos--Willwacher differs from ours.
In comparison to their work, our main contribution is the definition of the quasi-isomorphism between the CDGAs $\OmPA^{*}(\FM_{M}(k))$ and the small, explicit Lambrechts-Stanley model $\GG{A}(k)$, which has the advantage of being finite-dimensional and much more computable than $\Graphs_{S(\tilde{H}(M))}(k)$.

\noindent
\emph{Applications.}
Ordered configuration spaces appear in many places in topology and geometry.
Therefore, thanks to Theorems~\ref{thm:intermediary} and~\ref{thm.A}, the explicit model $\GG{A}(k)$ provides an efficient computational tool in many concrete situations.

To illustrate this, we show how to apply our results to compute factorization homology, an invariant of framed $n$-manifolds defined from an $\mathtt{E}_{n}$-algebra~\cite{AyalaFrancis2015}.
Let $M$ be a framed manifold with Poincaré duality model $A$, and $B$ be an $n$-Poisson algebras, i.e.\ an algebra over the operad $H_{*}(\mathtt{E}_{n})$.
Our results shows that we can compute the factorization homology of $M$ with coefficients in $B$ just from $\GG{A}$ and $B$.
As an application, we compute factorization homology with coefficients in a higher enveloping algebra of a Lie algebra (Proposition~\ref{prop.cmp-knudsen}), recovering a theorem of Knudsen~\cite{Knudsen2016}.

The Taylor tower in the Goodwillie--Weiss calculus of embeddings may be computed in a similar manner~\cite{GoodwillieWeiss1999,BoavidaWeiss2013}.
It follows from a result of~\cite[Section 5.1]{Turchin2013} that $\FM_{M}$ may be used for this purpose.
Therefore our theorem shows that $\GG{A}$ may also be used for computing this Taylor tower.

\noindent
\emph{Roadmap.}
In Section~\ref{sec.backgr-recoll}, we lay out our conventions and recall the necessary background.
This includes dg-modules and \abbr{CDGA}s, (co)operads and their (co)modules, semi-algebraic sets and \abbr{PA} forms.
We also recall basic results on the Fulton--MacPherson compactifications of configuration spaces $\FM_{n}(k)$ and $\FM_{M}(k)$, and the main ideas of Kontsevich's proof of the formality of the little disks operads using the \abbr{CDGA}s of \abbr{PA} forms on the spaces $\FM_{n}(k)$.
We use the formalism of operadic twisting, which we recall, to deal with signs more easily.
Finally, we recollect the necessary background on Poincaré duality \abbr{CDGA}s and the Lambrechts--Stanley \abbr{CDGA}s.
In Section~\ref{sec.model}, we build out of the Lambrechts--Stanley \abbr{CDGA}s a Hopf right $\enV$-comodule $\GG{A}$.

In Section~\ref{sec.label-graphs-stat}, we construct the labeled graph complex $\Graphs_{R}$ which will be used to connect $\GG{A}$ to $\OmPA^{*}(\FM_{M})$.
The construction is inspired by Kontsevich's construction of the unlabeled graph complex $\Graphs_{n}$.
It is done in several steps.
The first step is to consider a graded module of labeled graphs, $\Gra_{R}$.
In order to be able to map $\Gra_{R}$ into $\OmPA^{*}(\FM_{M})$, we recall the construction of what is called a ``propagator'' in the mathematical physics literature.
We then ``twist'' $\Gra_{R}$ to obtain a new object $\Tw \Gra_{R}$, which consists of graphs with two kinds of vertices: ``external'' and ``internal''.
Finally we must reduce our graphs to obtain a new object, $\Graphs_{R}$, by removing all the connected components with only internals vertices in the graphs using a ``partition function'' (a function which resembles the Chern--Simons invariants).

In Section~\ref{sec.proof-theorem}, we prove that the zigzag of Hopf right comodule morphisms between $\GG{A}$ and $\OmPA^{*}(\FM_{M})$ is a weak equivalence.
We first connect our graph complex $\Graphs_{R}$ to the Lambrechts--Stanley \abbr{CDGA}s $\GG{A}$.
This requires vanishing results about the partition function.
Then we end the proof of the theorem by showing that all the morphisms are quasi-isomorphisms.
Finally we study the cases $S^{2}$ and $S^{3}$.

In Section~\ref{sec.fact-homol-e_n}, we use our model to compute factorization homology of framed manifolds and we compare the result to a complex obtained by Knudsen.
In Section~\ref{sec.oriented-surfaces} we work out a variant of our theorem for the only simply connected surface using the formality of the framed little $2$-disks operad, and we present a conjecture about higher dimensional oriented manifolds.

For convenience, we provide a glossary of our main notations at the end of this paper.

\section{Background and recollections}
\label{sec.backgr-recoll}

\subsection{DG-modules and \abbr{CDGA}s}
\label{sec:dg-modules-cdgas}

We consider differential graded modules (dg-modules) over the base field $\R$.
Unless otherwise indicated, (co)homology of spaces is considered with real coefficients.
All our dg-modules will have a cohomological grading, $V = \bigoplus_{n \in \Z} V^{n}$.
All the differentials raise degrees by one: $\deg(dx) = \deg(x) + 1$.
We say that a dg-module is of finite type if it is finite dimensional in each degree.
Let $V[k]$ be the desuspension, defined by $(V[k])^{n} = V^{n+k}$.
For dg-modules $V,W$ and homogeneous elements $v \in V, \, w \in W$, we let $(v \otimes w)^{21} \coloneqq (-1)^{(\deg v)(\deg w)} w \otimes v$ and we extend this linearly to the tensor product.
Moreover, given an element $X \in V \otimes W$, we will often use a variant of Sweedler's notation to express $X$ as a sum of elementary tensors, $X \coloneqq \sum_{(X)} X' \otimes X'' \in V \otimes W$.

We call \textbf{CDGAs} the (graded) commutative unital algebras in dg-modules.
In general, for a \abbr{CDGA} $A$, we let $\mu_{A} : A^{\otimes 2} \to A$ be its product.
For a dg-module $V$, we let $S(V)$ be the free unital symmetric algebra on $V$.

We will need a model category structure on the category of \abbr{CDGA}s.
We use the model category structure given by the general result of~\cite{Hinich1997} for categories of algebras over operads.
The weak equivalences are the quasi-isomorphisms, the fibrations are the surjective morphisms, and the cofibrations are characterized by the left lifting property with respect to acyclic fibrations.
A path object for the initial \abbr{CDGA} $\R$ is given by $\APL^{*}(\Delta^{1}) = S(t,dt)$, the \abbr{CDGA} of polynomials forms on the interval.
It is equipped with an inclusion $\R \lhook\joinrel\xrightarrow{\sim} \APL^{*}(\Delta^{1})$, and two projections $\ev_{0}, \ev_{1} : \APL^{*}(\Delta^{1}) \xrightarrow{\sim} \R$ given by setting $t = 0$ or $t = 1$.
Two morphisms $f,g : A \to B$ with cofibrant source are \textbf{homotopic} if there exists a homotopy $h : A \to B \otimes \APL^{*}(\Delta^{1})$ such that the following diagram commutes:
\[
  \begin{tikzcd}[row sep = small, column sep = large]
    {} & A \ar[dl, "f" swap] \ar[d, "h"] \ar[dr, "g"] \\
    B & B \otimes \APL^{*}(\Delta^{1}) \ar[l, "\id \otimes \ev_{0}" {near start}, "\sim" {near start, swap}] \ar[r, "\id \otimes \ev_{1}" {near start, swap}, "\sim" {near start}] & B
  \end{tikzcd}
  .
\]

Many of the \abbr{CDGA}s that appear in this paper are $\Z$-graded.
However, to deserve the name ``model of $X$'', a \abbr{CDGA} should be connected to $\APL^{*}(X)$ only by $\N$-graded \abbr{CDGA}s.
The next proposition shows that considering this larger category does not change our statement.

\begin{prop}
  \label{prop:grading}
  Let $A$, $B$ be two $\N$-graded \abbr{CDGA}s which are homologically connected, i.e.\ $H^{0}(A) = H^{0}(B) = \R$.
  If $A$ and $B$ are quasi-isomorphic as $\Z$-graded \abbr{CDGA}s, then they also are as $\N$-graded \abbr{CDGA}s.
\end{prop}
\begin{proof}
  This follows from the results of~\cite[§II.6.2]{Fresse2017}.
  Let us temporarily denote $\mathsf{cdga}_{\N}$ the category of $\N$-graded \abbr{CDGA}s ($dg^{*}\Com$ in~\cite{Fresse2017}) and $\mathsf{cdga}_{\Z}$ the category of $\Z$-graded \abbr{CDGA}s ($dg\Com$ in~\cite{Fresse2017}).
  Note that in~\cite{Fresse2017}, $\Z$-graded \abbr{CDGA}s are homologically graded, but we can use the usual correspondence $A^{i} = A_{-i}$ to keep our convention that all dg-modules are cohomologically graded.
  There is an obvious inclusion $\iota : \mathsf{cdga}_{\N} \to \mathsf{cdga}_{\Z}$, which clearly defines a full functor that preserves and reflects quasi-isomorphisms.

  Let $\mathbb{B}^{m}$ be the dg-module $\R$ concentrated in degree $m$, let $\mathbb{E}^{m}$ be the dg-module given by two copies of $\R$ in respective degree $m-1$ and $m$ such that $d_{\mathbb{E}^{m}}$ is the identity of $\R$ in these degrees (hence $\mathbb{E}^{m}$ is acyclic), and let $i : \mathbb{B}^{m} \to \mathbb{E}^{m}$ be the obvious inclusion.
  The model category $\mathsf{cdga}_{\N}$ is equipped with a set of generating cofibrations given by the morphisms $S(i) : S(\mathbb{B}^{m}) \to S(\mathbb{E}^{m})$ and of the morphism $\varepsilon : S(\mathbb{B}^{0}) \to \R$.
  Recall that a cellular complex of generating cofibrations is a \abbr{CDGA} obtained by a sequential colimit $R = \operatorname{colim}_{k} R_{\langle k \rangle}$, where $R_{\langle 0 \rangle} = \R$ and $R_{\langle k+1 \rangle}$ is obtained from $R_{\langle k \rangle}$ by a pushout of generating cofibrations along attaching maps $h: S(\mathbb{B}^m)\rightarrow R_{\langle k \rangle}$.
  In~\cite[§II.6.2]{Fresse2017}, the expression ``connected generating cofibrations'' is used for the generating cofibrations of the form $S(i): S(\mathbb{B}^m)\rightarrow S(\mathbb{E}^m)$ with $m>0$.

  In the proof of \cite[Proposition II.6.2.8]{Fresse2017}, it is observed that, if $A$ is homologically connected, then the attaching map $h: S(\mathbb{B}^0)\rightarrow A$ associated to a generating cofibration $\varepsilon: S(\mathbb{B}^0)\rightarrow\R$ necessarily reduces to the augmentation $\varepsilon: S(\mathbb{B}^0)\rightarrow \R$ followed by the inclusion as the unit $\R \subset A$.
  Thus a pushout of the generating cofibration $\varepsilon : S(\mathbb{B}^0)\rightarrow \R$ reduces to a neutral operation in this case.
  In the proof of \cite[Proposition II.6.2.8]{Fresse2017}, it is deduced from this observation that any homologically connected algebra admits a resolution $R_A \xrightarrow{\sim} A$ such that $R_A$ is a cellular complex of connected generating cofibrations.
  Connected generating cofibrations are also cofibrations in $\mathsf{cdga}_{\Z}$ after applying $\iota$.
  Moreover $\iota$ preserves colimits.
  It follows that $\iota R_{A}$ is cofibrant in $\mathsf{cdga}_{\Z}$ too.

  By hypothesis, $\iota A$ and $\iota B$ are weakly equivalent in $\mathsf{cdga}_{\Z}$, hence $\iota R_{A}$ and $\iota B$ are also weakly equivalent (because $\iota$ clearly preserves quasi-isomorphisms), through a zigzag $\iota R_{A} \xleftarrow{\sim} \cdot \xrightarrow{\sim} \iota B$.
  As $\iota R_{A}$ is cofibrant (and all \abbr{CDGA}s are fibrant), we can find a direct quasi-isomorphism $\iota R_{A} \xrightarrow{\sim} \iota B$ and therefore a zigzag $\iota A \xleftarrow{\sim} \iota R_{A} \xrightarrow{\sim} \iota B$ which only involves $\N$-graded \abbr{CDGA}s.
\end{proof}

\subsection{(Co)operads and their right (co)modules}
\label{sec:coop-their-right}

We assume basic proficiency with Hopf (co)operads and (co)modules over (co)operads, see e.g.~\cite{Fresse2009,LodayVallette2012,Fresse2017}.
We index our (co)operads by finite sets instead of integers to ease the writing of some formulas.
If $W \subset U$ is a subset, we write the quotient $U/W = (U \setminus W) \sqcup \{*\}$, where $*$ represents the class of $W$; note that $U/\varnothing = U \sqcup \{ * \}$.
An operad in dg-modules, for instance, is given by a functor from the category of finite sets and bijections (a symmetric collection) $\PP : U \mapsto \PP(U)$ to the category of dg-modules, together with a unit $\Bbbk \to \PP(\{*\})$ and composition maps $\circ_{W} : \PP(U/W) \otimes \PP(W) \to \PP(U)$ for every pair of sets $W \subset U$, satisfying the usual associativity, unity and equivariance conditions.
Dually, a cooperad $\CC$ is given by a symmetric collection, a counit $\CC(\{*\}) \to \Bbbk$, and cocomposition maps $\circ_{W}^{\vee} : \CC(U) \to \CC(U/W) \otimes \CC(W)$ for every pair $W \subset U$.

Let $\underline{k} = \{1,\dots,k\}$.
We recover the usual notion of a cooperad indexed by the integers by considering the collection $\{\CC(\underline{k})\}_{k \ge 0}$, and the cocomposition maps $\circ_{i}^{\vee} : \CC(\underline{k+l-1}) \to \CC(\underline{k}) \otimes \CC(\underline{l})$ corresponds to $\circ_{\{i, \dots, i+l-1\}}^{\vee}$.

Following Fresse~\cite[§II.9.3.1]{Fresse2017}, a ``Hopf cooperad'' is a cooperad in the category of \abbr{CDGA}s.
We do not assume that (co)operads are trivial in arity zero, but they will satisfy $\PP(\varnothing) = \Bbbk$ (resp.\ $\CC(\varnothing) = \Bbbk$).
Therefore we get (co)operad structures equivalent to the structure of $\Lambda$-(co)operads considered by Fresse~\cite[§II.11]{Fresse2017}, which he uses to model rational homotopy types of operads in spaces satisfying $\PP(0) = *$ (but we do not use this formalism in the sequel).

The result of Proposition~\ref{prop:grading} extends to Hopf cooperads (and to Hopf $\Lambda$-cooperads).
To establish this result, we still use a description of generating cofibrations of $N$-graded Hopf cooperads, which are given by morphisms of symmetric algebras of cooperads $S(i): S(\CC)\rightarrow S(\DD)$, where $i: \CC \rightarrow \DD$ is a dg-cooperad morphism that is injective in positive degrees (see \cite[§II.9.3]{Fresse2017} for details).
In the context of homologically connected cooperads, we can check that the pushout of such a Hopf cooperad morphism along an attaching map reduces to a pushout of a morphism of symmetric algebras of cooperads $S(\CC/\ker(i)) \rightarrow S(\DD)$, where we mod out by the kernel of the map $i: \CC \rightarrow \DD$ in degree $0$.
We deduce from this observation that any homologically connected $N$-graded Hopf cooperad admits a resolution by a cellular complex of generating cofibrations of the form $S(i): S(C)\rightarrow S(D)$, where the map $i$ is injective in all degrees (we again call such a generating cofibration connected).
The category of $\Z$-graded Hopf cooperads inherits a model structure, like the category of $\N$-graded Hopf cooperads considered in \cite[§II.9.3]{Fresse2017}.
Cellular complexes of connected generating cofibrations of $\N$-graded Hopf cooperads define cofibrations in the model category of $\Z$-graded Hopf cooperads yet, as in the proof of Proposition~\ref{prop:grading}.

Given an operad $\PP$, a right $\PP$-module is a symmetric collection $\MM$ equipped with composition maps $\circ_{W} : \MM(U/W) \otimes \PP(W) \to \MM(U)$ satisfying the usual associativity, unity and equivariance conditions.
A right comodule over a cooperad is defined dually.
If $\CC$ is a Hopf cooperad, then a right Hopf $\CC$-comodule is a $\CC$-comodule $\NN$ such that all the $\NN(U)$ are \abbr{CDGA}s and all the maps $\circ_{W}^{\vee}$ are morphisms of \abbr{CDGA}s.

\begin{defin}
  \label{def.mph-comod}
  Let $\CC$ (resp.\ $\CC'$) be a Hopf cooperad and $\NN$ (resp.\ $\NN'$) be a Hopf right comodule over $\CC$ (resp.\ $\CC'$).
  A \textbf{morphism of Hopf right comodules} is a pair $(f_{\NN}, f_{\CC})$ consisting of a morphism of Hopf cooperads $f_{\CC} : \CC \to \CC'$, and a map of Hopf right $\CC'$-comodules $f_{\NN} : \NN \to \NN'$, where $\NN$ has the $\CC$-comodule structure induced by $f_{\CC}$.
  It is a \textbf{quasi-isomorphism} if both $f_{\CC}$ and $f_{\NN}$ are quasi-isomorphisms in each arity.
  A Hopf right $\CC$-module $\NN$ is said to be \textbf{weakly equivalent} to a Hopf right $\CC'$-module $\NN'$ if the pair $(\NN, \CC)$ can be connected to the pair $(\NN', \CC')$ through a zigzag of quasi-isomorphisms.
\end{defin}

The next very general lemma can for example be found in~\cite[Section 5.2]{CamposWillwacher2016}.
Let $\CC$ be a cooperad, and see the \abbr{CDGA} $A$ as an operad concentrated in arity $1$.
Recall that $\CC \circ A = \bigoplus_{i \ge 0} \CC(i) \otimes_{\Sigma_{i}} A^{\otimes i}$ denotes the composition product of operads, where we view $A$ as an operad concentrated in arity $1$.
Then the commutativity of $A$ implies the existence of a distributive law $t : \CC \circ A \to A \circ \CC$, given in each arity by the morphism $t : \CC(\underline{n}) \otimes A^{\otimes n} \to A \otimes \CC(\underline{n})$ given by $x \otimes a_{1} \otimes \dots \otimes a_{n} \mapsto a_{1} \dots a_{n} \otimes x$.

\begin{lem}
  \label{lemma.comod-circle}
  Let $\NN$ be a right $\CC$-comodule, and see $A$ as an operad concentrated in arity $1$.
  Then $\NN \circ A$ is a right $\CC$-comodule through the map:
  \[ \NN \circ A \xrightarrow{\Delta_{\NN} \circ 1} \NN \circ \CC \circ A \xrightarrow{1 \circ t} \NN \circ A \circ \CC.
    \qed \]
\end{lem}

\subsection{Semi-algebraic sets and forms}
\label{sec:semi-algebraic-sets}

Kontsevich's proof of the formality of the little disks operads~\cite{Kontsevich1999} uses the theory of semi-algebraic sets, as developed in~\cite{KontsevichSoibelman2000,HardtLambrechtsTurchinVolic2011}.
A semi-algebraic set is a subset of $\R^{N}$ defined by finite unions of finite intersections of zero sets of polynomials and polynomial inequalities.
By the Nash--Tognoli Theorem~\cite{Nash1952,Tognoli1973}, any closed smooth manifold is algebraic hence semi-algebraic.

There is a functor $\OmPA^{*}$ of ``piecewise semi-algebraic (PA) differential forms'', analogous to de Rham forms.
If $X$ is a compact semi-algebraic set, then $\OmPA^{*}(X) \simeq \APL^{*}(X) \otimes_{\Q} \R$, i.e.\ the \abbr{CDGA} $\OmPA^{*}(X)$ models the real homotopy type of $X$~\cite[Theorem~6.1]{HardtLambrechtsTurchinVolic2011}.

A key feature of \abbr{PA} forms is that it is possible to compute integrals of ``minimal forms'' along fibers of ``\abbr{PA} bundles'', i.e.\ maps with local semi-algebraic trivializations~\cite[Section~8]{HardtLambrechtsTurchinVolic2011}.
A minimal form is of the type $f_{0} df_{1} \wedge \dots \wedge df_{k}$ where $f_{i} : M \to \R$ are semi-algebraic maps.
Given such a minimal form $\lambda$ and a \abbr{PA} bundle $p : M \to B$ with fibers of dimension $r$, there is a new form (which is not minimal in general), also called the pushforward of $\lambda$ along $p$:
\[ p_{*} \lambda \coloneqq \int_{p : M \to B} \lambda \in \OmPA^{k-r}(B). \]

In what follows, we use an extension of the fiberwise integration of minimal forms to the sub-\abbr{CDGA} of ``trivial forms'' given in \cite[Appendix C]{CamposWillwacher2016}.
Briefly recall that trivial forms are integrals of minimal forms along fibers of a trivial \abbr{PA} bundle (see \cite[Definition 81]{CamposWillwacher2016}).
In fact, in Section~\ref{sec.propagator}, we consider a certain form, the ``propagator'', which is not minimal but trivial in this sense, and we apply the extension of the fiberwise integration to this form.

The functor $\Omega^*_{PA}$ is monoidal, but not strongly monoidal, and contravariant.
Thus, given an operad $\PP$ in semi-algebraic sets, $\OmPA^{*}(\PP)$ is an ``almost'' Hopf cooperad and satisfies a slightly modified version of the cooperad axioms, as explained in~\cite[Definition 3.1]{LambrechtsVolic2014}.
Cooperadic structure maps are replaced by zigzags $\OmPA^{*}(\PP(U)) \xrightarrow{\circ_{W}^{*}} \OmPA^{*}(\PP(U/W) \times \PP(W)) \xleftarrow{\sim} \OmPA^{*}(\PP(U/W)) \times \OmPA^{*}(\PP(W))$ (where the second map is the Künneth morphism).
If $\CC$ is a Hopf cooperad, an ``almost'' morphism $f: \CC \to \OmPA^{*}(\PP)$ is a collection of \abbr{CDGA} morphisms $f_{U} : \CC(U) \to \OmPA^{*}(\PP(U))$ for all $U$, such that the following diagrams commute:
\[ \begin{tikzcd}[row sep = small]
    \CC(U) \ar[rr, "\circ_{W}^{\vee}"] \ar[d, "f_{U}"] && \CC(U/W) \otimes \CC(W) \ar[d, "f_{U/W} \otimes f_{W}"] \\
     \OmPA^{*}(\PP(U)) \ar[r, "\circ_{W}^{*}"] & \OmPA^{*}(\PP(U/W) \times \PP(W)) & \OmPA^{*}(\PP(U/W)) \otimes \OmPA^{*}(\PP(W)) \ar[l, "\sim" swap]
  \end{tikzcd} \]
Similarly, if $\MM$ is a $\PP$-module, then $\OmPA^{*}(\MM)$ is an ``almost'' Hopf right comodule over $\OmPA^{*}(\PP)$.
If $\NN$ is a Hopf right $\CC$-comodule, where $\CC$ is a cooperad equipped with an ``almost'' morphism $f : \CC \to \OmPA^{*}(\PP)$, then an ``almost'' morphism $g : \NN \to \OmPA^{*}(\MM)$ is a collection of \abbr{CDGA} morphisms $g_{U} : \NN(U) \to \OmPA^{*}(\MM(U))$ that make the following diagrams commute:
\[ \begin{tikzcd}[row sep = small]
    \NN(U) \ar[rr, "\circ_{W}^{\vee}"] \ar[d, "g_{U}"] && \NN(U/W) \otimes \CC(W) \ar[d, "g_{U/W} \otimes f_{W}"] \\
     \OmPA^{*}(\MM(U)) \ar[r, "\circ_{W}^{*}"] & \OmPA^{*}(\MM(U/W) \times \PP(W)) & \OmPA^{*}(\MM(U/W)) \otimes \OmPA^{*}(\PP(W)) \ar[l, "\sim" swap]
  \end{tikzcd} \]
We will generally omit the adjective ``almost'', keeping in mind that some commutative diagrams are a bit more complicated than at first glance.

\begin{rmk}
  There is a construction $\Omega^{*}_{\sharp}$ that turns a simplicial operad $\PP$ into a Hopf cooperad and such that a morphism of Hopf cooperads $\CC \to \Omega^{*}_{\sharp}(\PP)$ is the same thing as an ``almost'' morphism $\CC \to \APL^{*}(\PP)$, where $\APL^{*}$ is the functor of Sullivan forms~\cite[Section II.10.1]{Fresse2017}.
  Moreover there is a canonical collection of maps $(\Omega^{*}_{\sharp}(\PP))(U) \to \APL^{*}(\PP(U))$, which are weak equivalences if $\PP$ is a cofibrant operad.
  This functor is built by considering the right adjoint of the functor on operads induced by the Sullivan realization functor, which is monoidal.
  A similar construction can be extended to $\OmPA^{*}$ and to modules over operads.
  This construction allows us to make sure that the cooperads and comodules we consider truly encode the rational or real homotopy type of the initial operad or module (see \cite[§II.10.2]{Fresse2017}).
\end{rmk}

\subsection{Little disks and related objects}
\label{sec.little-disks}

The little disks operad $\mathtt{E}_{n}$ is a topological operad initially introduced by May and Boardman--Vogt~\cite{May1972,BoardmanVogt1973} to study iterated loop spaces.
Its homology $\mathtt{e}_{n} := H_{*}(\mathtt{E}_{n})$ is described by a theorem of Cohen~\cite{Cohen1976}: it is either the operad governing associative algebras for $n = 1$, or $n$-Poisson algebras for $n \ge 2$.
We also consider the linear dual $\enV \coloneqq H^{*}(\mathtt{E}_{n})$, which is a Hopf cooperad.

In fact, we use the Fulton--MacPherson operad $\FM_{n}$, which is an operad in spaces weakly equivalent to the little disks operad $\mathtt{E}_{n}$.
The components $\FM_{n}(k)$ are compactifications of the configuration spaces $\Conf_{k}(\R^{n})$, defined by using a real analogue due to Axelrod--Singer~\cite{AxelrodSinger1994} of the Fulton--MacPherson compactifications~\cite{FultonMacPherson1994}.
The idea of this compactification is to allow configurations where points become ``infinitesimally close''.
Then one uses insertion of such infinitesimal configurations to define operadic composition products on the spaces $\FM_n(k)$.
We refer to~\cite{Sinha2004} for a detailed treatment and to~\cite[Sections~5.1--5.2]{LambrechtsVolic2014} for a clear summary.
In both references, the name $C[k]$ is used for what we call $\FM_{n}(k)$.

The first two spaces $\FM_{n}(\varnothing) = \FM_{n}(\underline{1}) = *$ are singletons, and $\FM_{n}(\underline{2}) = S^{n-1}$ is a sphere.
We let the volume form of $\FM_{n}(\underline{2})$ be:
\begin{equation}
  \label{eq.def-vol}
  \vol_{n-1} \in \OmPA^{n-1}(S^{n-1}) = \OmPA^{n-1}(\FM_{n}(\underline{2}))
\end{equation}

The space $\FM_{n}(\underline{k})$ is a semi-algebraic stratified manifold, of dimension $nk-n-1$ for $k \geq 2$, and of dimension $0$ otherwise.
For $u \neq v \in U$, we can define the projection maps that forget all but two points in the configuration, $p_{uv} : \FM_{n}(U) \to \FM_{n}(\underline{2})$.
These projections are semi-algebraic bundles.

If $M$ is a manifold, the configuration space $\Conf_{k}(M)$ can similarly be compactified to give a space $\FM_{M}(\underline{k})$.
By forgetting points, we again obtain projection maps, for $u,v \in U$:
\begin{align}
  \label{eq:proj-fm-m}
  p_{u} : \FM_{M}(U) & \twoheadrightarrow \FM_{M}(\underline{1}) = M,
  & p_{uv} : \FM_{M}(U) & \twoheadrightarrow \FM_{M}(\underline{2}).
\end{align}
The two projections $p_{1}$ and $p_{2}$ are equal when restricted $\partial \FM_{M}(\underline{2})$, and they define a sphere bundle of rank $n-1$,
\begin{equation}
  \label{eq:sphere-bundle}
  p : \partial \FM_{M}(\underline{2}) \twoheadrightarrow M.
\end{equation}

When $M$ is framed, the collection of spaces $\FM_{M}$ assemble to form a topological right module over $\FM_{n}$, with composition products defined by insertion of infinitesimal configurations.
Moreover in this case, the sphere bundle $p : \partial \FM_{M}(\underline{2}) \to M$ is trivialized by:
\begin{equation}
  \label{eq:bundle-trivial}
  M \times S^{n-1} \cong \FM_{M}(\underline{1}) \times \FM_{n}(\underline{2}) \xrightarrow{\circ_{1}} \partial \FM_{M}(\underline{2}).
\end{equation}

Recall from Section~\ref{sec:semi-algebraic-sets} that we can endow $M$ with a semi-algebraic structure.
It is immediate that $\FM_{M}(\underline{k})$ is a stratified semi-algebraic manifold of dimension $nk$.
Moreover, the proofs of~\cite[Section~5.9]{LambrechtsVolic2014} can be adapted to show that the projections $p_{U} : \FM_{M}(U \sqcup V) \to \FM_{M}(U)$ are \abbr{PA} bundles.

\subsection{Operadic twisting}
\label{sec:operadic-twisting}

We will make use of the ``operadic twisting'' procedure in what follows~\cite{DolgushevWillwacher2015}.
Let us now recall this procedure, in the case of cooperads.

Let $\Lie_{n}$ be the operad governing shifted Lie algebras.
A $\Lie_{n}$-algebra is a dg-module $\mathfrak{g}$ equipped with a Lie bracket $[-,-] : \mathfrak{g}^{\otimes 2} \to \mathfrak{g}[1-n]$ of degree $1-n$, i.e.\ we have $[\mathfrak{g}^{i}, \mathfrak{g}^{j}] \subset \mathfrak{g}^{i+j+(1-n)}$.

\begin{rmk}
  The degree convention is such that there is an embedding of operads $\Lie_{n} \to H_{*}(\FM_{n})$, i.e.\ Poisson $n$-algebras are $\Lie_{n}$-algebras.
  The usual Lie operad is $\Lie_{1}$.
  This convention is consistent with~\cite{Willwacher2014}.
  However in~\cite{Willwacher2016}, the notation is $\Lie^{(n)} = \Lie_{n+1}$.
  In~\cite{DolgushevWillwacher2015}, only the unshifted operad $\Lie = \Lie_{1}$ is considered.
\end{rmk}

The operad $\Lie_{n}$ is quadratic Koszul (see e.g.~\cite[Section~13.2.6]{LodayVallette2012}), and as such admits a cofibrant resolution $\hoLie_{n} \coloneqq \Omega(K(\Lie_{n}))$, where $\Omega$ is the cobar construction and $K(\Lie_{n})$ is the Koszul dual cooperad of $\Lie$.
Algebras over $\hoLie_{n}$ are (shifted) $L_{\infty}$-algebras, also known as homotopy Lie algebras, i.e.\ dg-modules $\mathfrak{g}$ equipped with higher brackets $[-,\dots,-]_{k} : \mathfrak{g}^{\otimes k} \to \mathfrak{g}[3-k-n]$ (for $k \ge 1$) satisfying the classical $L_{\infty}$ equations.

Let $\CC$ be a cooperad (with finite-type components in each arity) equipped with a map to the dual of $\hoLie_{n}$.
This map can equivalently be seen as a Maurer--Cartan element in the following dg-Lie algebra~\cite[Section~6.4.2]{LodayVallette2012}:
\begin{equation}
  \label{eq:hom-sigma}
  \Hom_{\Sigma}(K(\Lie_{n}), \CC^{\vee})
  \coloneqq \biggl( \prod_{i \ge 0} \bigl( \CC^{\vee}(i) \otimes \R[-n]^{\otimes i} \bigr)^{\Sigma_{i}}[n], \partial, [-,-] \biggr),
\end{equation}
where we used the explicit description of the Koszul dual $K(\Lie_{n})$ as a shifted version of the cooperad encoding cocommutative coalgebras.
Given $f,g \in \Hom_{\Sigma}(K(\Lie_{n}), \CC^{\vee})$, their bracket is $[f,g] = f \star g \mp g \star f$, where $\star$ is given by:
\[ f \star g : K(\Lie_{n}) \xrightarrow{\text{cooperad}} K(\Lie_{n}) \circ K(\Lie_{n}) \xrightarrow{f \circ g} \CC^{\vee} \circ \CC^{\vee} \xrightarrow{\text{operad}} \CC^{\vee}. \]
An element $\mu \in \Hom_{\Sigma}(K(\Lie_{n}), \CC^{\vee})$ is said to satisfy the Maurer--Cartan equation if $\partial \mu + \mu \star \mu = 0$.
Such an element is called a twisting morphism in~\cite[Section~6.4.3]{LodayVallette2012}, and the equivalence with morphisms $\hoLie_{n} \to \CC^{\vee}$ (or dually $\CC \to \hoLie_{n}^{\vee}$) is~\cite[Theorem~6.5.7]{LodayVallette2012}.
In the sequel, we will alternate between the two points of view, morphisms or Maurer--Cartan elements.

There is an action of the symmetric group $\Sigma_{i}$ on $\underline{i} = \{ 1, \dots, i\}$.
As a graded module, the twist of $\CC$ with respect to $\mu$ is given by:
\begin{equation}
  \label{eq:def-tw-c}
  \Tw \CC(U) := \bigoplus_{i \geq 0} \bigl( \CC(U \sqcup \underline{i}) \otimes \R[n]^{\otimes i} \bigr)_{\Sigma_{i}}.
\end{equation}
The symmetric collection $\Tw \CC$ inherits a cooperad structure from $\CC$.
The differential of $\Tw \CC$ is the sum of the internal differential of $\CC$ with a differential coming from the action of $\mu$ that we now explain.
The action of $\mu$ is threefold, and the total differential of $\Tw \CC(U)$ can be expressed as:
\begin{equation}
  \label{eq:total-diff-tw-c}
  d_{\Tw\CC} \coloneqq d_{\CC} + (- \cdot \mu) + (- \cdot \mu_{1}) + (\mu_{1} \cdot -).
\end{equation}
Let us now explain these notations.
Let $i \ge 0$ be some fixed integer and let us describe the action of $\mu$ on $\CC(U \sqcup \underline{i}) \subset \Tw \CC(U)$ (up to degree shifts).
In what follows, for a set $J \subset \underline{i}$, we let $j \coloneqq \# J$, and $\underline{i}/J \cong \underline{i+j-1}$.

Recall that $\mu$ is a formal sum of elements $\CC(\underline{j})^{\vee}$ for $j \ge 0$.
The first action $(- \cdot \mu)$ is the sum over all subsets $J \subset \underline{i}$ of the following cocompositions:
\begin{equation}
  \label{eq:act-mu-1}
  \CC(U \sqcup \underline{i}) \xrightarrow{\circ_{J}^{\vee}} \CC(U \sqcup \underline{i}/J) \otimes \CC(J) \xrightarrow{\id \otimes \mu} \CC(U \sqcup \underline{i}/J) \otimes \R \cong \CC(U \sqcup \underline{i+j-1}).
\end{equation}

For the two other terms, we need the element $\mu_{1} \in \prod_{j \ge 0} \CC(\underline{j} \sqcup \{*\})^{\vee}$.
It is the sum over all possible ways of distinguishing one input of $\mu$ in each arity.
(Distinguishing one input does not respect the invariants in the definition of Equation~\eqref{eq:hom-sigma}, but taking the sum over all possible ways does.)

The second action $(- \cdot \mu_{1})$ is then the sum of the following cocompositions, over all subsets $J \subset \underline{i}$ and over all $* \in U$ (where we use the obvious bijection $U / \{*\} \cong U$):
\begin{equation}
  \label{eq:act-mu-2}
  \CC(U \sqcup \underline{i}) \xrightarrow{\circ_{\{u\} \sqcup J}^{\vee}} \CC\bigl( (U \sqcup \underline{i}) / (\{*\} \sqcup J) \bigr) \otimes \CC(\{*\} \sqcup J) \xrightarrow{\id \otimes \mu_{1}} \CC(U \sqcup \underline{i+j-1}),
\end{equation}

Finally, the third action $(\mu_{1} \cdot -)$ is the sum over all subsets $J \subset \underline{i}$ of the cocompositions (where we use the obvious bijection $(U \sqcup I) / (U \sqcup J) = \{*\} \sqcup I \setminus J$):
\begin{equation}
  \label{eq:act-mu-3}
  \CC(U \sqcup \underline{i}) \xrightarrow{\circ_{U \sqcup J}^{\vee}} \CC(\{*\} \sqcup \underline{i} \setminus J) \otimes \CC(U \sqcup J) \xrightarrow{\mu_{1} \otimes \id} \CC(U \sqcup J),
\end{equation}

\begin{lem}
  \label{lem:tw-hopf}
  If $\CC$ is a Hopf cooperad satisfying $\CC(\varnothing) = \Bbbk$, then $\Tw\CC$ inherits a Hopf cooperad structure.
\end{lem}
\begin{proof}
  To multiply an element of $\CC(U \sqcup I) \subset \Tw\CC(U)$ with an element of $\CC(U \sqcup J) \subset \Tw\CC(U)$, we use the maps $\CC(V) \xrightarrow{\circ_{\varnothing}^{\vee}} \CC(V / \varnothing) \otimes \CC(\varnothing) \cong \CC(V \sqcup \{*\})$ iterated several times, to obtain elements in $\CC(U \sqcup I \sqcup J)$ and the product.
\end{proof}

Moreover, we will need to twist right comodules over cooperads.
This construction is found (for operads) in~\cite[Appendix~C.1]{Willwacher2016}.
Let us fix a cooperad $\CC$ and a twist $\Tw \CC$ with respect to $\mu$ as above.
Given a right $\CC$-comodule $\MM$, we can also twist it with respect to $\mu$, as follows.
As a graded module, the object $\Tw\MM(U)$ is defined by:
\[ \Tw \MM(U) \coloneqq \prod_{i \ge 0} \bigl( \MM(U \sqcup \underline{i}) \otimes (\R[n])^{\otimes i} \bigr)_{\Sigma_{i}}. \]
The comodule structure is inherited from $\MM$.
The total differential is the sum:
\begin{equation}
  \label{eq:total-diff-tw-comod}
  d_{\Tw\MM} \coloneqq d_{\MM} + (- \cdot \mu) + (- \cdot \mu_{1}),
\end{equation}
where  $(- \cdot \mu)$ and $(- \cdot \mu_{1})$ are as in Equations~\eqref{eq:act-mu-1} and~\eqref{eq:act-mu-2} but using the comodule structure.
Note that $\MM$ is only a right module, so there can be no term $(\mu_{1} \cdot -)$ in this differential.
Lemma~\ref{lem:tw-hopf} has an immediate extension:
\begin{lem}
  \label{lem:tw-hopf-comod}
  If $\CC$ is a Hopf cooperad satisfying $\CC(\varnothing) = \Bbbk$ and $\MM$ is a Hopf right $\CC$-comodule, then $\Tw \MM$ inherits a Hopf right $(\Tw\CC)$-comodule structure.
  \qed
\end{lem}

\subsection{Formality of the little disks operad}
\label{sec.formality}

Kontsevich's proof of the formality of the little disks operads~\cite[Section~3]{Kontsevich1999}, can be summarized by the fact that $\OmPA^{*}(\FM_{n})$ is weakly equivalent to $\enV$ as a Hopf cooperad.
For detailed proofs, we refer to~\cite{LambrechtsVolic2014}.

We outline this proof here as we will mimic its pattern for our theorem.
The idea of the proof is to construct a Hopf cooperad $\Graphs_{n}$.
The elements of $\Graphs_{n}$ are formal linear combinations of special kinds of graphs, with two types of vertices, numbered ``external'' vertices and unnumbered ``internal'' vertices.
The differential is defined combinatorially by edge contraction.
It is built in such a way that there exists a zigzag $\enV \xleftarrow{\sim} \Graphs_{n} \xrightarrow{\sim} \OmPA^{*}(\FM_{n})$.
The first map is the quotient by the ideal of graphs containing internal vertices.
The second map is defined using integrals along fibers of the \abbr{PA} bundles $\FM_{n}(U \sqcup I) \to \FM_{n}(U)$ which forget some points in the configuration.
An induction argument shows that the first map is a quasi-isomorphism, and the second map is easily seen to be surjective on cohomology.

In order to deal with signs more easily, we use (co)operadic twisting (Section~\ref{sec:operadic-twisting}).
Thus the Hopf cooperad $\Graphs_{n}$ is not the same as the Hopf cooperad $\mathcal{D}$ from~\cite{LambrechtsVolic2014}, see Remark~\ref{rmk:comp-d-graphs}.

\noindent
\emph{The cohomology of $\mathtt{E}_{n}$.}
The cohomology $\enV(U) = H^{*}(\mathtt{E}_{n}(U))$ has a classical presentation due to Arnold~\cite{Arnold1969} and Cohen~\cite{Cohen1976}.
We have
\begin{equation}
  \label{eq.def-enV}
  \enV(U) = S(\omega_{uv})_{u, v \in U} / I,
\end{equation}
where the generators $\omega_{uv}$ have cohomological degree $n-1$, and the ideal $I$ encoding the relations is generated by the polynomials (called Arnold relations):
\begin{equation}
  \label{eq:5}
  \omega_{uu} = 0; \; \omega_{vu} = (-1)^{n} \omega_{uv}; \; \omega_{uv}^{2} = 0; \; \omega_{uv} \omega_{vw} + \omega_{vw} \omega_{wu} + \omega_{wu} \omega_{uv} = 0.
\end{equation}

The cooperadic structure maps are given by (where $[u], [v] \in U/W$ are the classes of $u$ and $v$ in the quotient):
\begin{align}
  \circ^{\vee}_{W} : \enV(U)
  & \to \enV(U/W) \otimes \enV(W),
  & \omega_{uv}
  & \mapsto
  \begin{cases}
    1 \otimes \omega_{uv}, & \text{if } u,v \in W; \\
    \omega_{[u][v]} \otimes 1, & \text{otherwise}.
  \end{cases}
  \label{eq.coop-gra-n}
\end{align}

\noindent
\emph{Graphs with only external vertices.}
The intermediary cooperad of graphs, $\Graphs_{n}$, is built in several steps.
In the first step, define a cooperad of graphs with only external vertices, with generators $e_{uv}$ of degree $n-1$:
\begin{equation}
  \label{eq.def-gra-n}
  \Gra_{n}(U) = \bigl( S(e_{uv})_{u, v \in U} / (e_{uv}^{2} = e_{uu} = 0, e_{vu} = (-1)^{n} e_{uv}), d = 0 \bigr).
\end{equation}
The definition of $\Gra_{n}(U)$ is almost identical to the definition of $\enV(U)$, except that we do not kill the Arnold relations.

The \abbr{CDGA} $\Gra_{n}(U)$ is spanned by words of the type $e_{u_{1} v_{1}} \dots e_{u_{r} v_{r}}$.
Such a word can be viewed as a graph with $U$ as the set of vertices, and an edge between $u_{i}$ and $v_{i}$ for each factor $e_{u_{i} v_{i}}$.
For example, $e_{uv}$ is a graph with a single edge from $u$ to $v$ (see Equation~\eqref{fig.exa-graph} for another example).
Edges are oriented, but for even $n$ an edge is identified with its mirror (so we can forget orientations), while for odd $n$ it is identified with the opposite of its mirror.
In pictures, we do not draw orientations, keeping in mind that for odd $n$, they are necessary to get precise signs.
Graphs with double edges or edges between a vertex and itself are set to zero.
Given such a graph, its set of edges $E_{\Gamma} \subset \binom{U}{2}$ is well-defined.
The vertices of these graphs are called ``external'', in contrast with the internal vertices that are going to appear in the next part.

\begin{equation}
  \label{fig.exa-graph}
  e_{12} e_{13} e_{56} =
  \begin{gathered} \begin{tikzpicture}
      \node [extv] (1) {1};
      \node [extv, right = .5cm of 1] (2) {2};
      \node [extv, above = .5cm of 1] (3) {3};
      \node [extv, right = .5cm of 3] (4) {4};
      \node [extv, right = .5cm of 2] (5) {5};
      \node [extv, above = .5cm of 5] (6) {6};
      \draw (1) -- (2); \draw (1) -- (3); \draw (5) -- (6);
    \end{tikzpicture} \end{gathered}
  \in \Gra_{n}(\underline{6})
\end{equation}

The multiplication of the \abbr{CDGA} $\Gra_{n}(U)$, from this point of view, consists of gluing two graphs along their vertices.
The cooperadic structure map $\circ^{\vee}_{W} : \Gra_{n}(U) \to \Gra_{n}(U/W) \otimes \Gra_{n}(W)$ maps a graph $\Gamma$ to $\pm \Gamma_{U/W} \otimes \Gamma_{W}$ such that $\Gamma_{W}$ is the full subgraph of $\Gamma$ with vertices $W$ and $\Gamma_{U/W}$ collapses this full subgraph to a single vertex.
On generators, $\circ_{W}^{\vee}$ is defined by a formula which is in fact identical to Equation~\eqref{eq.coop-gra-n}, replacing $\omega_{??}$ by $e_{??}$.
This implies that the cooperad $\Gra_{n}$ maps to $\enV$ by sending $e_{uv}$ to $\omega_{uv}$.

There is a morphism $\omega' : \Gra_{n} \to \OmPA^{*}(\FM_{n})$ given on generators by:
\begin{align}
  \label{eq:2}
  \omega' : \Gra_{n}(U)
  & \to \OmPA^{*}(\FM_{n}(U)),
  &
  \Gamma
  & \mapsto \bigwedge_{(u,v) \in E_{\Gamma}} p_{uv}^{*}(\vol_{n-1}),
\end{align}
where $p_{uv} : \FM_{n}(U) \to \FM_{n}(\underline{2})$ is the projection map defined in Section~\ref{sec.little-disks}, and $\vol_{n-1}$ is the volume form of $\FM_{n}(2) \cong S^{n-1}$ from Equation~\eqref{eq.def-vol}.

\noindent
\emph{Twisting.}
The second step of the construction is cooperadic twisting, using the procedure outlined in Section~\ref{sec:operadic-twisting}.
The Hopf cooperad $\Gra_{n}$ maps into $\Lie_{n}^{\vee}$ as follows.
The cooperad $\Lie_{n}^{\vee}$ is cogenerated by $\Lie_{n}^{\vee}(\underline{2})$, and on cogenerators the cooperad map is given by sending $e_{12} \in \Gra_{n}(\underline{2})$ to the cobracket in $\Lie_{n}^{\vee}(\underline{2})$ and all the other graphs to zero.
This map to $\Lie_{n}^{\vee}$ yields a map to $\hoLie_{n}^{\vee}$ by composition with the canonical map $\Lie_{n}^{\vee} \xrightarrow\sim \hoLie_{n}^{\vee}$.
In the dual basis, the corresponding Maurer--Cartan element $\mu$ is given by:
\begin{equation}
  \label{fig.maurer-cartan-element}
  \mu \coloneqq e_{12}^{\vee} =
  \begin{tikzpicture}[baseline=(1.base)]
    \node[extv] (1) {$1$};
    \node[extv] (2) [right = 1cm of 1] {$2$};
    \draw (1) -- (2);
  \end{tikzpicture}
  \in \Gra_{n}^{\vee}(\underline{2})
\end{equation}
The cooperad $\Gra_{n}$ satisfies $\Gra_{n}(\varnothing) = \R$.
Thus by Lemma~\ref{lem:tw-hopf}, $\Tw \Gra_{n}$ inherits a Hopf cooperad structure, which we now explicitly describe.

\begin{figure}[htbp]
  \centering
  \[
  \begin{gathered} \begin{tikzpicture}
      \node[intv] (int) {};
      \node[extv] (1) [above = 0.3cm of int] {$1$};
      \node[extv] (2) [below left = 0.4cm of int] {$2$};
      \node[extv] (3) [below right = 0.4cm of int] {$3$};
      \draw (1) -- (int); \draw (2) -- (int); \draw (3) -- (int);
    \end{tikzpicture} \end{gathered}
  \xmapsto{\; d \;}
  \pm
  \begin{gathered} \begin{tikzpicture}
      \node (int) {};
      \node[extv] (1) [above = 0.3cm of int] {$1$};
      \node[extv] (2) [below left = 0.4cm of int] {$2$};
      \node[extv] (3) [below right = 0.4cm of int] {$3$};
      \draw (1) -- (2); \draw (2) -- (3);
    \end{tikzpicture} \end{gathered}
  \pm
  \begin{gathered} \begin{tikzpicture}
      \node (int) {};
      \node[extv] (1) [above = 0.3cm of int] {$1$};
      \node[extv] (2) [below left = 0.4cm of int] {$2$};
      \node[extv] (3) [below right = 0.4cm of int] {$3$};
      \draw (2) -- (3); \draw (3) -- (1);
    \end{tikzpicture} \end{gathered}
  \pm
  \begin{gathered} \begin{tikzpicture}
      \node (int) {};
      \node[extv] (1) [above = 0.3cm of int] {$1$};
      \node[extv] (2) [below left = 0.4cm of int] {$2$};
      \node[extv] (3) [below right = 0.4cm of int] {$3$};
      \draw (3) -- (1); \draw (1) -- (2);
    \end{tikzpicture} \end{gathered}
  \]
  \caption{The differential of $\Tw\Gra_n$. This particular example shows that the Arnold relation (the RHS) is killed up to homotopy.}
  \label{fig.diff-gra-n}
\end{figure}
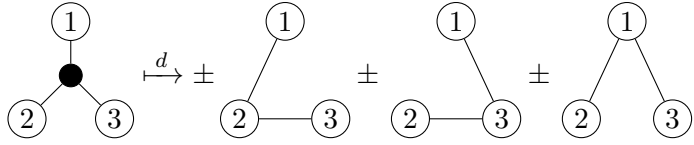

The dg-module $\Tw\Gra_{n}(U)$ is spanned by graphs with two types of vertices: external vertices, which correspond to elements of $U$ and that we will picture as circles with the name of the label in $U$ inside, and indistinguishable internal vertices, corresponding to the elements of $\underline{i}$ in Equation~\eqref{eq:def-tw-c} and that we will draw as black points.
For example, the graph inside the differential in the left hand side of Figure~\ref{fig.diff-gra-n} represents an element of $\Tw \Gra_{n}(U)$ with $U = \{1,2,3\}$ and $i = 1$.
The degree of an edge is still $n-1$, the degree of an external vertex is still $0$, and the degree of an internal vertex is $-n$.

The product of $\Tw \Gra_{n}(U)$ glues graphs along their external vertices only.
Compared to Lemma~\ref{lem:tw-hopf}, this coincides with adding isolated internal vertices (by iterating the cooperad structure map $\circ_{\varnothing}^{\vee}$) and gluing along all vertices.

Let us now describe the differential adapted from~\cite[Section~6.4]{LambrechtsVolic2014} (see Remark~\ref{rmk:d-hat-diff} for the differences).
We first give the final result, then we explain how it is obtained from the description in Section~\ref{sec:operadic-twisting}.
An edge is said to be \textbf{contractible} if it connects any vertex to an internal vertex, except if it connects a univalent internal vertices to a vertex which is not a univalent internal vertex.
The differential of a graph $\Gamma$ is the sum:
\[ d\Gamma = \sum_{\substack{e \in E_{\Gamma} \\ \text{contractible}}} \pm \Gamma / e, \]
where $\Gamma/e$ is $\Gamma$ with $e$ collapsed, and $e$ ranges over all contractible edges.

Let us now explain how to compare this with the description in Section~\ref{sec:operadic-twisting}, see also~\cite[Appendix~I.3]{Willwacher2014} for a detailed description.
Recall that the Maurer--Cartan element $\mu$ (Equation~\eqref{fig.maurer-cartan-element}) is equal to $1$ on the graph with exactly two vertices and one edge, and vanishes on all other graphs.
Recall from Equation~\eqref{eq:total-diff-tw-c} that the differential of $\Tw \Gra_{n}$ has three terms: $(- \cdot \mu) + (- \cdot \mu_{1}) + (\mu_{1} \cdot -)$, plus the differential of $\Gra_{n}$ which vanishes.
Let $\Gamma$ be some graph. Then $d\Gamma = \Gamma \cdot \mu + \Gamma \cdot \mu_{1} + \mu_{1} \cdot \Gamma$ where:
\begin{itemize}
\item The element $\Gamma \cdot \mu$ is the sum over all ways of collapsing a subgraph $\Gamma' \subset \Gamma$ with only internal vertices, the result being $\mu(\Gamma') \Gamma / \Gamma'$.
  This is nonzero only if $\Gamma'$ has exactly two vertices and one edge.
  Thus this summand corresponds to contracting all edges between two (possibly univalent) internal vertices in $\Gamma$.
\item The element $\Gamma \cdot \mu_{1}$ is the sum over all ways of collapsing a subgraph $\Gamma' \subset \Gamma$ with exactly one external vertex (and any number of internal vertices), with result $\mu(\Gamma') \Gamma / \Gamma'$.
  This summands corresponds to contracting all edges between one external vertex and one internal (possibly univalent) vertex.
\item The element $\mu_{1} \cdot \Gamma$ is the sum over all ways of collapsing a subgraph $\Gamma' \subset \Gamma$ containing all the external vertices, with result $\mu_{1}(\Gamma/\Gamma') \Gamma'$.
  The coefficient $\mu_{1}(\Gamma/\Gamma')$ can only be nonzero if $\Gamma$ is obtained from $\Gamma'$ by adding a univalent internal vertex.
  A careful analysis of the signs~\cite[Appendix~I.3]{Willwacher2014} shows that this cancels out with the contraction of edges connected to univalent internal vertices from the other two summands, unless \emph{both} endpoints of the edge are univalent and internal (and hence disconnected from the rest of the graph), in which cases the same term appears three times, and only two cancel out (see \cite[Fig.~3]{Willwacher2014} for the dual picture).
\end{itemize}

\begin{defin}
  \label{def:internally}
  A graph is \textbf{internally connected} if it remains connected when the external vertices are deleted.
  It is easily checked that as a commutative algebra, $\Tw \Gra_{n}(U)$ is freely generated by such graphs.
\end{defin}

The morphisms $\enV \gets \Gra_{n} \xrightarrow{\omega'} \OmPA^{*}(\FM_{n})$ extend along the inclusion $\Gra_{n} \subset \Tw\Gra_{n}$ as follows.
The extended morphism $\Tw\Gra_{n} \to \enV$ simply sends a graph with internal vertices to zero.
We need to check that this commutes with the differential.
We thus need to determine when a graph with internal vertices (sent to zero) can have a differential with no internal vertices (possibly sent to a nonzero element in $\enV$).
The differential decreases the number of internal vertices by exactly one.
So by looking at generators (internally connected graphs) we can look at the case of graphs with a single internal vertex connected to some external vertices.
Either the internal vertex is univalent, but then the edge is not contractible and the differential vanishes.
Or the internal vertex is connected to more than one external vertices.
In this case, one check that the differential of the graph is zero modulo the Arnold relations, (see~\cite[Introduction]{LambrechtsVolic2014} and Figure~\ref{fig.diff-gra-n} for an example).

The extended morphism $\omega : \Tw\Gra_{n} \to \OmPA^{*}(\FM_{n})$ (see~\cite[Definition~14]{Kontsevich1999} and~\cite[Chapter~9]{LambrechtsVolic2014} where the analogous integral is denoted $\widehat{I}$) sends a graph $\Gamma \in \Gra_{n}(U \sqcup I) \subset \Tw\Gra_{n}(U)$ to:
\begin{equation}
  \label{eq:4}
  \omega(\Gamma) \coloneqq \int_{\FM_{n}(U \sqcup I) \xrightarrow{p_{U}} \FM_{n}(U)} \omega'(\Gamma) = (p_{U})_{*}(\omega'(\Gamma)),
\end{equation}
where $p_{U}$ is the projection that forgets the points of the configuration corresponding to $I$, and the integral is an integral along the fiber of this \abbr{PA} bundle (see Section~\ref{sec:semi-algebraic-sets}).
Note that the volume form on the sphere is minimal, hence $\omega'(\Gamma)$ is minimal and therefore we can compute this integral.

\begin{rmk}
  \label{rmk:d-hat-diff}
  This Hopf cooperad is different from the module of diagrams $\widehat{\mathcal{D}}$ introduced in~\cite[Section~6.2]{LambrechtsVolic2014}: $\Tw \Gra_{n}$ is the quotient of $\widehat{\mathcal{D}}$ by graphs with multiple edges and loops.
  The analogous integral $\widehat{I} : \widehat{\mathcal{D}} \to \OmPA^{*}(\FM_{n})$ is from \cite[Chapter~9]{LambrechtsVolic2014}.
  It vanishes on graphs with multiple edges and loops by \cite[Lemmas 9.3.5, 9.3.6]{LambrechtsVolic2014}, so $\omega$ is well-defined.
  Moreover the differential is slightly different.
  In~\cite{LambrechtsVolic2014} some kind of edges, called ``dead ends''~\cite[Definition~6.1.1]{LambrechtsVolic2014}, are not contractible.
  When restricted to graphs without multiple edges or loops, these are edges connected to univalent internal vertices.
  But in $\Tw\Gra_{n}$, edges connecting two internal vertices that are both univalent are contractible (see \cite[Fig.~3]{Willwacher2014} for the dual picture).
  This does not change $\widehat{I}$, which vanishes on graphs with univalent internal vertices anyway \cite[Lemma 9.3.8]{LambrechtsVolic2014}.
  Note that $\widehat{\mathcal{D}}$ is \emph{not} a Hopf cooperad \cite[Example~7.3.2]{LambrechtsVolic2014} due to multiple edges.
\end{rmk}

\noindent
\emph{Reduction.}
The cooperad $\Tw\Gra_{n}$ does not have the homotopy type of the cooperad $\OmPA^{*}(\FM_{n})$.
It is reduced by quotienting out all the graphs with connected components consisting exclusively of internal vertices.
This is a bi-ideal generated by $\Tw \Gra_{n}(\varnothing)$, thus the resulting quotient is a Hopf cooperad:
\[ \Graphs_{n} \coloneqq \Tw \Gra_{n} / \bigl( \Tw \Gra_{n}(\varnothing) \bigr). \]

\begin{rmk}
  \label{rmk:comp-d-graphs}
  This Hopf cooperad is not isomorphic to the Hopf cooperad $\mathcal{D}$ from~\cite[Section~6.5]{LambrechtsVolic2014}.
  We allow internal vertices of any valence, whereas in $\mathcal{D}$ internal vertices must be at least trivalent.
  There is a quotient map $\Graphs_{n} \to \mathcal{D}$, which is a quasi-isomorphism by~\cite[Proposition~3.8]{Willwacher2014}.
  The statement of~\cite[Proposition~3.8]{Willwacher2014} is actually about the dual operads, but as we work over a field and the spaces we consider have finite-type cohomology, this is equivalent.
  The notation is also different: the couple $(\Graphs_{n}, \mathtt{fGraphs}_{n,c})$ in~\cite{Willwacher2014} denotes $(\mathcal{D}^{\vee}, \Graphs_{n}^{\vee})$ in~\cite{LambrechtsVolic2014}.
\end{rmk}

One checks that the two morphisms $\enV \gets \Tw\Gra_{n} \to \OmPA^{*}(\FM_{n})$ factor through the quotient (the first one because graphs with internal vertices are sent to zero, the second one because $\omega$ vanishes on graphs with only internal vertices by~\cite[Lemma~9.3.7]{LambrechtsVolic2014}).
The resulting zigzag $\enV \gets \Graphs_{n} \to \OmPA^{*}(\FM_{n})$ is then a zigzag of weak equivalence of Hopf cooperads thanks to the proof of~\cite[Theorem~2]{Kontsevich1999} (or~\cite[Theorem~8.1]{LambrechtsVolic2014} and the discussion at the beginning of~\cite[Chapter~10]{LambrechtsVolic2014}), combined with the comparison between $\mathcal{D}$ and $\Graphs_{n}$ from~\cite[Proposition~3.8]{Willwacher2014} (see Remark~\ref{rmk:comp-d-graphs}).

\subsection{Poincaré duality \abbr{CDGA} models}
\label{sec.poinc-dual-models}

The model for $\OmPA^{*}(\FM_{M})$ relies on a Poincaré duality model of $M$.
We mostly borrow the terminology and notation from~\cite{LambrechtsStanley2008}.

Fix an integer $n$ and let $A$ be a connected \abbr{CDGA} (i.e.\ $A = \R \oplus A^{\ge 1}$).
An \textbf{orientation} on $A$ is a linear map $A^{n} \to \R$ satisfying $\varepsilon \circ d = 0$ (which we often view as a chain map $A \to \R[-n]$) such that the induced pairing
\begin{align}
  \label{eq:7}
  \langle-,-\rangle : A^{k} \otimes A^{n-k}
  & \to \R,
  & a \otimes b
  & \mapsto \varepsilon(ab)
\end{align}
is non-degenerate for all $k$.
This implies that $A = A^{\leq n}$, and that $\varepsilon : A^{n} \to \R$ is an isomorphism.
The pair $(A,\varepsilon)$ is called a \textbf{Poincaré duality \abbr{CDGA}}.
If $A$ is such a Poincaré duality \abbr{CDGA}, then so is its cohomology.
The following ``converse'' has been shown by Lambrechts--Stanley.
\begin{thrm}[Direct corollary of {Lambrechts--Stanley~\cite[Theorem~1.1]{LambrechtsStanley2008}}]
  \label{thm.lambrechts-stanley}
  Let $M$ be a simply connected semi-algebraic closed oriented manifold.
  Then there exists a zigzag of quasi-isomorphisms of \abbr{CDGA}s
  \[ A \xleftarrow{\;\rho\;} R \xrightarrow{\;\sigma\;} \OmPA^{*}(M), \]
  such that $A$ is a Poincaré duality \abbr{CDGA} of dimension $n$, $R$ is a quasi-free \abbr{CDGA} generated in degrees $\geq 2$, $\sigma$ factors through the sub-\abbr{CDGA} of trivial forms.
\end{thrm}

\begin{proof}
  We refer to Section~\ref{sec:semi-algebraic-sets} for a reminder on trivial forms.
  We pick a minimal model $R$ of the manifold $M$ (over $\mathbb{R}$) and a quasi-isomorphism from $R$ to the sub-CDGA of trivial forms in $\Omega_{PA}^*(M)$, which exists because the sub-CDGA of trivial forms is quasi-isomorphic to $\Omega_{PA}^*(M)$ (see Section 1.3), and hence, is itself a real model for $M$.
  We compose this new quasi-isomorphism this the inclusion to eventually get a quasi-isomorphism $\sigma: R\rightarrow\Omega_{PA}^*(M)$ which factors through the sub-CDGA of trivial forms, and we set $\varepsilon = \int_M\sigma(-): R\rightarrow\mathbb{R}[-n]$.
  The CDGA $R$ is of finite type because $M$ is a closed manifold.
  Hence, we can apply the Lambrechts-Stanley Theorem~\cite[Theorem~1.1]{LambrechtsStanley2008} to the pair $(R,\epsilon)$ to get the Poincaré duality algebra $A$ of our statement.
\end{proof}

Let $A$ be a Poincaré duality \abbr{CDGA} of finite type and let $\{a_{i}\}$ be a homogeneous basis of $A$.
Consider the dual basis $\{ a_{i}^{*} \}$ with respect to the duality pairing, i.e.\ $\varepsilon(a_{i} a_{j}^{\*}) = \delta_{ij}$ is given by the Kronecker symbol.
Then the \textbf{diagonal cocycle} is defined by the following formula and is independent of the chosen homogeneous basis (see e.g.~\cite[Definition~8.16]{FelixOpreaTanre2008}:
\begin{equation}
  \label{eq.def-delta}
  \Delta_{A} \coloneqq \sum_{i} (-1)^{|a_{i}|} a_{i} \otimes a_{i}^{*} \in A \otimes A.
\end{equation}
The element $\Delta_{A}$ is a cocycle of degree $n$ (this follows from $\varepsilon \circ d = 0$).
It satisfies $\Delta_{A}^{21} = (-1)^{n} \Delta_{A}$ (where $(-)^{21}$ is defined in Section~\ref{sec:dg-modules-cdgas}).
For all $a \in A$ it satisfies the equation $(a \otimes 1) \Delta_{A} = (1 \otimes a) \Delta_{A}$.
There is a volume form,
\[ \vol_{A} \coloneqq \varepsilon^{-1}(1_{\R}) \in A^{n}. \]
The product $\mu_{A} : A \otimes A \to A$ sends $\Delta_{A}$ to $\chi(A) \cdot \vol_{A}$, where $\chi(A)$ is the Euler characteristic of $A$.
We will need the following technical result later.

\begin{prop}
  \label{lemma.zigzag-a-r-omega}
  One can choose the zigzag of Theorem~\ref{thm.lambrechts-stanley} such there exists a symmetric cocycle $\Delta_{R} \in R \otimes R$ of degree $n$ satisfying $(\rho \otimes \rho)(\Delta_{R}) = \Delta_{A}$.
  If $\chi(M) = 0$ we can moreover choose it so that $\mu_{R}(\Delta_{R}) = 0$.
\end{prop}

\begin{proof}
  We follow closely the proof of~\cite{LambrechtsStanley2008} to obtain the result. 
  Recall that the proof of~\cite{LambrechtsStanley2008} has two different cases: $n \le 6$, where the manifold is automatically formal and hence $A = H^{*}(M)$, and $n \ge 7$, where the \abbr{CDGA} is built out of an inductive argument.
  We split our proof along these two cases.

  Let us first deal with the case $n \ge 7$.
  When $n \geq 7$, the proof of Lambrechts and Stanley builds a zigzag of weak equivalences $A \xleftarrow{\rho} R \gets R' \to \OmPA^{*}(M)$, where $R'$ is the minimal model of $M$, the \abbr{CDGA} $R$ is obtained from $R'$ by successively adjoining generators of degree $\geq n/2+1$, and the Poincaré duality \abbr{CDGA} $A$ is a quotient of $R$ by an ideal of ``orphans''.
  We let $\varepsilon : R' \to \R[-n]$ be the composite $R' \to \OmPA^{*}(M) \xrightarrow{\int_{M}} \R[-n]$.

  The minimal model $R'$ is quasi-free, and since $M$ is simply connected it is generated in degrees $\geq 2$.
  The \abbr{CDGA} $R$ is obtained from $R'$ by a cofibrant cellular extension, adjoining cells of degree greater than $2$.
  It follows that $R$ is cofibrant and quasi-freely generated in degrees $\geq 2$.
  Composing with $R' \to \OmPA^{*}(M)$ yields a morphism $\sigma : R \to \OmPA^{*}(M)$ and we therefore get a zigzag $A \gets R \to \OmPA^{*}(M)$.

  The morphism $\rho$ is a quasi-isomorphism, so there exists some cocycle $\tilde\Delta \in R \otimes R$ such that $\rho(\tilde\Delta) = \Delta_{A} + d\alpha$ for some $\alpha$.
  By surjectivity of $\rho$ (it is a quotient map) there is some $\beta$ such that $\rho(\beta) = \alpha$; we let $\Delta' = \tilde\Delta - d\beta$, and now $\rho(\Delta') = \Delta_{A}$.

  Let us assume for the moment that $\chi(M) = 0$.
  Then the cocycle $\mu_{R}(\Delta') \in R$ satisfies $\rho(\mu_{R}(\Delta')) = \mu_{A}(\Delta_{A}) = 0$, i.e.\ it is in the kernel of $\rho$.
  It follows that the cocycle $\Delta'' = \Delta' - \mu_{R}(\Delta') \otimes 1$ is still mapped to $\Delta_{A}$ by $\rho$, and satisfies $\mu_{R}(\Delta'') = 0$.
  If $\chi(M) \neq 0$ we just let $\Delta'' = \Delta'$.
  Finally we symmetrize $\Delta''$ to get the $\Delta_{R}$ of the lemma, which satisfies all the requirements.

  Let us now deal with the case $n \le 6$.
  The \abbr{CDGA} $\OmPA^{*}(M)$ is formal~\cite[Proposition 4.6]{NeisendorferMiller1978}.
  We choose $A = (H^{*}(M), d_{A} = 0)$, and $R$ to be the minimal model of $M$, which maps into both $A$ and $\OmPA^{*}(M)$ by quasi-isomorphisms.
  The rest of the proof is now identical to the previous case.
\end{proof}

\subsection{The Lambrechts--Stanley \abbr{CDGA}s}
\label{sec.kriz-lambr-stanl}

We now give the definition of the \abbr{CDGA} $\GG{A}(k)$ from~\cite[Definition~3.4]{LambrechtsStanley2008a}, where it is called $F(A,k)$.

Let $A$ be a Poincaré duality \abbr{CDGA} of dimension $n$ and let $k$ be an integer.
For $1 \leq i \neq j \leq k$, let $\iota_{i} : A \to A^{\otimes k}$ be defined by $\iota_{i}(a) = 1^{\otimes i-1} \otimes a \otimes 1^{\otimes k-i-1}$, and let $\iota_{ij} : A \otimes A \to A^{\otimes k}$ be given by $\iota_{ij}(a \otimes b) = \iota_{i}(a) \cdot \iota_{j}(b)$.
Recalling the description of $\enV$ in Equation~\eqref{eq.def-enV}, the \abbr{CDGA} $\GG{A}(k)$ is defined by:
\begin{equation}
  \label{eq:9}
  \GG{A}(k) \coloneqq \bigl( A^{\otimes k} \otimes \enV(k) / ( \iota_{i}(a) \cdot \omega_{ij} = \iota_{j}(a) \cdot \omega_{ij} ), d \omega_{ij} = \iota_{ij}(\Delta_{A}) \bigr).
\end{equation}

The fact that this is well-defined is proved in~\cite[Lemma~3.2]{LambrechtsStanley2008a}.
We will call these \abbr{CDGA}s the Lambrechts--Stanley \abbr{CDGA}s, or \textbf{LS CDGAs} for short.
For example $\GG{A}(0) = \R$, $\GG{A}(1) = A$, and $\GG{A}(2)$ is isomorphic to:
\[ \GG{A}(2) \cong \bigl( (A \otimes A) \oplus (A \otimes \omega_{12}), d(a \otimes \omega_{12}) = (a \otimes 1) \cdot \Delta_{A} = (1 \otimes a) \cdot \Delta_{A}). \]

Recall that there always exists a Poincar\'e duality model of $M$ (Section~\ref{sec.poinc-dual-models}).
When $M$ is a simply connected closed manifold, a theorem of Lambrechts--Stanley~\cite[Theorem~10.1]{LambrechtsStanley2008a} implies that for any such $A$,
\begin{equation}
  \label{eq.ea-fmm-same-cohom}
  H^{*}(\GG{A}(k); \mathbb{Q}) \cong H^{*}(\FM_{M}(k); \mathbb{Q}) \; \text{as graded modules}.
\end{equation}

\section{The Hopf right comodule model \texorpdfstring{$\GG{A}$}{G\_A}}
\label{sec.model}

In this section we describe the Hopf right $\enV$-comodule derived from the \abbr{LS CDGA}s of Section~\ref{sec.kriz-lambr-stanl}.
From now on we fix a simply connected smooth closed manifold $M$.
Following Section~\ref{sec.little-disks}, we endow $M$ with a fixed semi-algebraic structure.
Note that for now, we do not impose any further conditions on $M$, but a key argument (Proposition~\ref{prop.strong-vanishing}) will require $\dim M \ge 4$.
We also fix a arbitrary Poincaré duality \abbr{CDGA} model $A$ of $M$.
We then define the right comodule structure of $\GG{A}$ as follows, using the cooperad structure of $\enV$ given by Equation~\eqref{eq.coop-gra-n}:

\begin{prop}
  \label{prop.comodule-eA}
  If $\chi(M) = 0$, then the following maps are well-defined on $\GG{A} = \{ \GG{A}(k) \}_{k \geq 0}$ and endow it with a Hopf right $\enV$-comodule structure:
  \begin{align}
    \circ_{W}^{\vee} : A^{\otimes U} \otimes \enV(U)
    & \to \bigl( A^{\otimes (U/W)} \otimes \enV(U/W) \bigr) \otimes \enV(W),
    \notag \\
    (a_{u})_{u \in U} \otimes \omega
    & \mapsto
      \underbrace{ ((a_{u})_{u \in U \setminus W} \otimes \textstyle{\prod_{w \in W} a_{w}}) }_{\in A^{\otimes (U/W)}}
      \otimes
      \underbrace{\circ_{W}^{\vee}(\omega)}_{\mathclap{\in \enV(U/W) \otimes \enV(W)}}.
      \label{eq.comodule}
  \end{align}
\end{prop}

In informal terms, $\circ_{W}^{\vee}$ multiplies together all the elements of $A$ indexed by $W$ on the $A^{\otimes U}$ factor and indexes the result by $* \in U/W$, while it applies the cooperadic structure map of $\enV$ on the other factor.
Note that if $W = \varnothing$, then $\circ_{W}^{\vee}$ adds a factor of $1_{A}$ (the empty product) indexed by $* \in U/\varnothing = U \sqcup \{*\}$.

\begin{proof}
  We split the proof in three parts: factorization of the maps through the quotient, compatibility with the differential, and compatibility of the maps with the cooperadic structure of $\enV$.

  Let us first prove that the comodule structure maps we wrote factor through the quotient.
  Since $A$ is commutative and $\enV$ is a Hopf cooperad, the maps of the proposition commute with multiplication.
  The ideals defining $\GG{A}(U)$ are multiplicative ideals.
  Hence it suffices to show that the maps~\eqref{eq.comodule} take the generators $(\iota_{u}(a) - \iota_{v}(a)) \cdot \omega_{uv}$ of the ideal to elements of the ideal in the target.
  We simply check each case, using Equations~\eqref{eq.coop-gra-n} and~\eqref{eq.comodule}:
  \begin{itemize}
  \item If $u,v \in W$, then $\circ_{W}^{\vee}(\iota_{u}(a) \omega_{uv}) = \iota_{*}(a) \otimes \omega_{uv}$, which is also equal to $\circ_{W}^{\vee}(\iota_{v}(a) \omega_{uv})$.
  \item Otherwise, we have $\circ_{W}^{\vee}(\iota_{u}(a) \omega_{uv}) = \iota_{[u]}(a) \omega_{[u][v]} \otimes 1$, which is equal to $\iota_{[v]}(a) \omega_{[u][v]} \otimes 1 = \circ_{W}^{\vee}(\iota_{v}(a) \omega_{uv})$ modulo the relations.
  \end{itemize}

  Let us now prove that they are compatible with the differential.
  It is again sufficient to prove this on generators.
  The equality $\circ_{W}^{\vee}(d(\iota_{u}(a))) = d(\circ_{W}^{\vee}(\iota_{u}(a)))$ is immediate.
  For $\omega_{uv}$ we again check the three cases.
  Recall that since our manifold has vanishing Euler characteristic, $\mu_{A}(\Delta_{A}) = 0$.
  \begin{itemize}
  \item If $u,v \in W$, then $\circ_{W}^{\vee}(d\omega_{uv}) = \iota_{*}(\mu_{A}(\Delta_{A})) = 0$, while by definition $d(\circ_{W}^{\vee}(\omega_{uv})) = d(1 \otimes \omega_{uv}) = 0$.
  \item Otherwise, $\circ_{W}^{\vee}(d \omega_{uv}) = \iota_{[u][v]}(\Delta_{A}) \otimes 1$, which is equal to $d (\circ_{W}^{\vee}(\omega_{uv})) = d(\omega_{[u][v]} \otimes 1)$.
  \end{itemize}

  We finally prove that the structure maps are compatible with the cooperad structure of $\enV$.
  Let $\Com^{\vee}$ be the cooperad governing cocommutative coalgebras.
  It follows from Lemma~\ref{lemma.comod-circle} that $\Com^{\vee} \circ A = \{ A^{\otimes k}\}_{k \ge 0}$ inherits a $\Com^{\vee}$-comodule structure.
  Therefore the arity-wise tensor product (see~\cite[Section 5.1.12]{LodayVallette2012}, where this operation is called the Hadamard product) $(\Com^{\vee} \circ A) \hadprod \enV \coloneqq \{ A^{\otimes k} \otimes \enV(k) \}_{k \geq 0}$ is a $(\Com^{\vee} \hadprod \enV)$-comodule.
  The cooperad $\Com^{\vee}$ is the unit of $\hadprod$.
  Hence the $(\Com^{\vee} \circ A) \hadprod \enV$ is an $\enV$-comodule.
  It remains to make the easy check that the resulting comodule maps are given by Equation~\eqref{eq.comodule}.
\end{proof}

\section{Labeled graph complexes}
\label{sec.label-graphs-stat}

In this section we construct the intermediary comodule, $\Graphs_{R}$, used to prove our theorem, where $R$ is a suitable cofibrant \abbr{CDGA} quasi-isomorphic to $A$ and $\OmPA^{*}(M)$ (Theorem~\ref{thm.lambrechts-stanley}).
We will construct a zigzag of \abbr{CDGA}s of the form:
\[ \GG{A} \gets \Graphs_{R} \to \OmPA^{*}(\FM_{M}). \]

The construction of $\Graphs_{R}$ follows the same pattern as the construction of $\Graphs_{n}$ in Section~\ref{sec.formality}, but with the vertices of the graph labeled by elements of $R$.
The differential moreover mimics the definition of the differential of $\GG{A}$, together with a differential that mimics the one of $\Graphs_{n}$.

If $\chi(M) = 0$, then the collections $\GG{A}$ and $\Graphs_{R}$ are Hopf right comodules respectively over $\enV$ and over $\Graphs_{n}$, and the left arrow is a morphism of comodules between $(\GG{A}, \enV)$ and $(\Graphs_{R}, \Graphs_{n})$.
When $M$ is moreover framed, $\OmPA^{*}(\FM_{M})$ is a Hopf right comodule over $\OmPA^{*}(\FM_{n})$, and the right arrow is then a morphism from $(\Graphs_{R}, \Graphs_{n})$ to $(\OmPA^{*}(\FM_{M}), \OmPA^{*}(\FM_{n}))$.

In order to deal with signs more easily and make sure that the differential squares to zero, we want to use the formalism of operadic twisting, as in the definition of $\Graphs_{n}$.
But when $\chi(M) \neq 0$ there is no comodule structure, so we make a detour through graphs with loops (Section~\ref{sec.graphs-with-loops} below), see Remark~\ref{rmk:detour}.

\subsection{Graphs with loops and multiple edges}
\label{sec.graphs-with-loops}

We first define a variant $\Graphs^{\circlearrowleft}_{n}$ of $\Graphs_{n}$, where graphs are allowed to have ``loops'' (also sometimes known as ``tadpoles'') and multiple edges, see \cite[Section~3]{Willwacher2014}.
For a finite set $U$, the \abbr{CDGA} $\Gra^{\circlearrowleft}_{n}(U)$ is presented by (where the generators have degree $n-1$):
\[ \Gra^{\circlearrowleft}_{n}(U) \coloneqq \bigl( S(e_{uv})_{u, v \in U} / (e_{vu} = (-1)^{n} e_{uv}), d = 0 \bigr). \]

The difference with Equation~\eqref{eq.def-gra-n} is that we no longer set $e_{uu} = e_{uv}^{2} = 0$.
Note that $\Gra^{\circlearrowleft}_{n}(U)$ is actually free as a \abbr{CDGA}: given an arbitrary linear order on $U$, $\Gra_{n}^{\circlearrowleft}(U)$ is freely generated by the generators $\{ e_{uv}\}_{u \le v \in U}$.

\begin{rmk}
  When $n$ is even, $e_{uv}^{2} = 0$ since $\deg e_{uv} = n-1$ is odd; and when $n$ is odd, the relation $e_{uu} = (-1)^{n} e_{uu}$ implies $e_{uu} = 0$.
  In other words, for even $n$, there are no multiple edges, and for odd $n$, there are no loops~\cite[Remark~3.1]{Willwacher2014}.
\end{rmk}

The dg-modules $\Gra_{n}^{\circlearrowleft}(U)$ form a Hopf cooperad, like $\Gra_{n}$, with cocomposition given by a formula similar to the definition of Equation~\eqref{eq.coop-gra-n}:
\begin{align}
  \circ^{\vee}_{W} : \Gra_{n}^{\circlearrowleft}(U)
  & \to \Gra_{n}^{\circlearrowleft}(U/W) \otimes \Gra_{n}^{\circlearrowleft}(W),
    \notag \\
  e_{uv}
  & \mapsto
  \begin{cases}
    e_{**} \otimes 1 + 1 \otimes e_{uv}, & \text{if } u,v \in W; \\
    e_{[u][v]} \otimes 1, & \text{otherwise}.
  \end{cases}
\end{align}

This new cooperad has a graphical description similar to $\Gra_{n}$.
The cooperad $\Gra_{n}$ is the quotient of $\Gra^{\circlearrowleft}_{n}$ by the ideal generated by the loops and the multiple edges.
The difference in the cooperad structure is that when we collapse a subgraph, we sum over all ways of choosing whether edges are in the subgraph or not; if they are not, then they yield a loop.
For example:

\begin{equation}
  \label{eq:exa-coop-gra-loop}
  \begin{tikzpicture}[baseline=.3cm]
    \node [extv] (1) {1};
    \node [extv, right = 1cm of 1] (2) {2};
    \node [extv, above right = 0.5cm of 1] (3) {3};
    \draw[red] (1) -- (2);
    \draw (2) -- (3);
  \end{tikzpicture}
  \xmapsto{\; \circ_{\{1,2\}}^{\vee} \;}
  \left(
  \begin{tikzpicture}[baseline=.33cm]
    \node [extv] (1) {$*$};
    \node [extv, above = 0.4cm of 1] (3) {3};
    \draw (1) to (3);
    \draw[red] (1) to[loop right] (1);
  \end{tikzpicture}
  \otimes
  \begin{tikzpicture}[baseline=(1.base)]
    \node [extv] (1) {1};
    \node [extv, right = 0.5cm of 1] (2) {2};
  \end{tikzpicture}
  \right)
  +
  \left(
  \begin{tikzpicture}[baseline=.33cm]
    \node [extv] (1) {$*$};
    \node [extv, above = 0.4cm of 1] (3) {3};
    \draw (1) to (3);
  \end{tikzpicture}
  \otimes
  \begin{tikzpicture}[baseline=(1.base)]
    \node [extv] (1) {1};
    \node [extv, right = 0.5cm of 1] (2) {2};
    \draw[red] (1) -- (2);
  \end{tikzpicture}
  \right)
\end{equation}

The element $\mu \coloneqq e_{12}^{\vee} \in (\Gra^{\circlearrowleft}_{n})^{\vee}(\underline{2})$ still defines a morphism $\Gra_{n}^{\circlearrowleft} \to \hoLie_{n}^{\vee}$, which allows us to define the twisted Hopf cooperad $\Tw \Gra^{\circlearrowleft}_{n}$.
It has a graphical description similar to $\Tw\Gra_{n}$ with internal and external vertices.
Finally we can quotient by graphs containing connected component consisting exclusively of internal vertices to get a Hopf cooperad:
\[ \Graphs^{\circlearrowleft}_{n} \coloneqq \Tw \Gra_{n}^{\circlearrowleft} / (\text{connected components with only internal vertices}). \]

\begin{rmk}
  The Hopf cooperad $\Tw\Gra_{n}^{\circlearrowleft}$ is slightly different from $\widehat{\mathcal{D}}$ from \cite[Section 6]{LambrechtsVolic2014}.
  First the cocomposition is different, and the first term of the RHS in Equation~\eqref{eq:exa-coop-gra-loop} would not appear in $\widehat{\mathcal{D}}$.
  The differential is also slightly different: an edge connected to two univalent internal vertices -- hence disconnected from the rest of the graph -- is contractible here (see \cite[Section 3]{Willwacher2014} and Remark~\ref{rmk:d-hat-diff}).
  This fixes the failure of $\widehat{\mathcal{D}}$ to be a cooperad \cite[Example 7.3.2]{LambrechtsVolic2014}.
\end{rmk}

\subsection{Labeled graphs with only external vertices: \texorpdfstring{$\Gra_{R}$}{Gra\_R}}
\label{sec.external-vertices}

We construct a collection of \abbr{CDGA}s $\Gra_{R}$, corresponding to the first step in the construction of $\Graphs_{n}$ of Section~\ref{sec.formality}.
We first apply the formalism of Section~\ref{sec.poinc-dual-models} to $\OmPA^{*}(M)$ in order to obtain a Poincaré duality \abbr{CDGA} out of $M$, thanks to Theorem~\ref{thm.lambrechts-stanley}.
We thus fix a zigzag of quasi-isomorphisms $A \xleftarrow{\rho} R \xrightarrow{\sigma} \OmPA^{*}(M)$, where $A$ is a Poincaré duality \abbr{CDGA}, $R$ is a cofibrant \abbr{CDGA}, and $\sigma$ factors through the sub-\abbr{CDGA} of trivial forms (see Section~\ref{sec:semi-algebraic-sets}).

Recall the definition of the diagonal cocycle $\Delta_{A} \in (A \otimes A)^{n}$ from Equation~\eqref{eq.def-delta}.
Recall also Proposition~\ref{lemma.zigzag-a-r-omega}, where we fixed a symmetric cocycle $\Delta_{R} \in (R \otimes R)^{n}$ such that $(\rho \otimes \rho)(\Delta_{R}) = \Delta_{A}$.
Moreover recall that if $\chi(M) = 0$, then $\mu_{A}(\Delta_{A}) = 0$, and we choose $\Delta_{R}$ such that $\mu_{R}(\Delta_{R}) = 0$ too.

\begin{defin}
  \label{def.gra-r}
  Let \abbr{CDGA} of \textbf{labeled graphs with loops} on the set $U$ be:
  \[ \Gra^{\circlearrowleft}_{R}(U) \coloneqq \bigl( R^{\otimes U} \otimes \Gra^{\circlearrowleft}_{n}(U), d e_{uv} = \iota_{uv}(\Delta_{R}) \bigr). \]
\end{defin}
This \abbr{CDGA} is well-defined because $\Gra_{n}^{\circlearrowleft}(U)$ is free as a \abbr{CDGA}, hence $\Gra^{\circlearrowleft}_{R}(U)$ is a relative Sullivan algebra in the terminology of~\cite[Section~14]{FelixHalperinThomas2001}.

\begin{rmk}
  This definition is valid for any \abbr{CDGA} $R$ and any symmetric cocycle $\Delta_{R}$.
  \label{rmk:r-hypotheses}
  We need $R$ as in Proposition~\ref{lemma.zigzag-a-r-omega} to connect $\Gra_{R}^{\circlearrowleft}$ with $\GG{A}$ and $\OmPA^{*}(\FM_{M})$.
\end{rmk}

\begin{rmk}
  \label{rmk.loop-multi}
  It follows that the differential of a loop is $d e_{uu} = \iota_{uu}(\Delta_{R}) = \iota_{u}(\mu_{R}(\Delta_{R}))$, which is zero when $\chi(M) = 0$.
\end{rmk}

\begin{prop}
  \label{prop.gra-tad-comod}
  The collection $\Gra^{\circlearrowleft}_{R}(U)$ forms a Hopf right $\Gra^{\circlearrowleft}_{n}$-comodule.
\end{prop}

This is true even if $\chi(M) \neq 0$ thanks to the introduction of the loops..

\begin{proof}
  The proof of this proposition is almost identical to the proof of Proposition~\ref{prop.comodule-eA}.
  If we forget the extra differential (keeping only the internal differential of $R$), then $\Gra^{\circlearrowleft}_{R}$ is the arity-wise tensor product $(\Com^{\vee} \circ R) \hadprod \Gra^{\circlearrowleft}_{n}$, which is automatically a Hopf $\Gra^{\circlearrowleft}_{n}$-right comodule.
  Checking the compatibility with the differential involves almost exactly the same equations as Proposition~\ref{prop.comodule-eA}, except that when $u,v \in W$ we have:
  \[ \circ_{W}^{\vee}(d(e_{uv})) = \iota_{*}(\mu_{R}(\Delta_{R})) \otimes 1 = d(e_{**} \otimes 1 + 1 \otimes e_{uv}) = d(\circ_{W}^{\vee}(e_{uv})), \]
  where $de_{**} = \iota_{*}(\mu_{R}(\Delta_{R}))$ by Remark~\ref{rmk.loop-multi}, and $d(1 \otimes e_{uv}) = 0$ by definition.
\end{proof}

We now give a graphical interpretation of Definition~\ref{def.gra-r}, in the spirit of Section~\ref{sec.graphs-with-loops}.
We view $\Gra^{\circlearrowleft}_{R}(U)$ as spanned by graphs with $U$ as set of vertices, and each vertex has a label which is an element of $R$.
The $\Gra^{\circlearrowleft}_{n}$-comodule structure collapses subgraphs as before, and the label of the collapsed vertex is the product of all the labels in the subgraph.
An example of graph in $\Gra_{R}^{\circlearrowleft}(\underline{3})$ is given by (where $x,y,z \in R$):
\begin{equation}
  \label{eq:exa-gra-r-loop}
  \begin{tikzpicture}[baseline=(1.base), every label/.style={font=\small}]
    \node [extv, label = {$x_{1}$}] (1) {1};
    \node [extv, label = {$x_{2}$}, right = 1cm of 1] (2) {2};
    \node [extv, label = {$x_{3}$}, right = 1cm of 2] (3) {3};
    \draw (1) -- (2) (1) to[bend right] (3) (1) to[loop left] (1);
  \end{tikzpicture}
\end{equation}

The product glues two graphs along their vertices, multiplying the labels in the process.
The differential of $\Gamma$, as defined in Definition~\ref{def.gra-r}, is the sum of $d_{R}$, the internal differential of $R$ acting on each label (one at a time), together with the sum over the edges $e \in E_{\Gamma}$ of the graph $\Gamma \setminus e$ with that edge removed and the labels of the endpoints multiplied by the factors of $\Delta_{R} = \sum_{(\Delta_{R})} \Delta_{R}' \otimes \Delta_{R}'' \in R \otimes R$, where we use Sweedler's notation (Section~\ref{sec:dg-modules-cdgas}).
We will often write $d_{\mathrm{split}}$ for this differential, to contrast it with the differential that contracts edges which will occur in the complex $\Tw\Gra^{\circlearrowleft}_{R}$ defined later on.
If $e$ is a loop, then in the corresponding term of $d\Gamma$ the vertex incident to $e$ has its label multiplied by $\mu_{R}(\Delta_{R})$, while the loop is removed.
For example, we have (gray vertices can be either internal or external and $x,y \in R$):

\[ \begin{tikzpicture}[baseline]
    \node[unkv, label = {\scriptsize $x$}] (1) {};
    \node[unkv, label = {\scriptsize $y$}] (2) [right = 1cm of 1] {};
    \node (x1) [left = 0.5cm of 1] {};
    \node (y1) [above left = .7cm of 1] {};
    \node (z1) [below left = .7cm of 1] {};
    \node (x2) [right = 0.5cm of 2] {};
    \node (y2) [above right = .7cm of 2] {};
    \node (z2) [below right = .7cm of 2] {};
    \draw (1) -- (2);
    \draw (x1) -- (1); \draw (y1) -- (1); \draw(z1) -- (1);
    \draw (x2) -- (2); \draw (y2) -- (2); \draw(z2) -- (2);
  \end{tikzpicture}
  \xmapsto{\; d \;}
  \sum_{(\Delta_R)}
  \begin{tikzpicture}[baseline]
    \node[unkv, label = {\scriptsize $x \Delta_{R}'$}] (1) {};
    \node[unkv, label = {\scriptsize $y \Delta_{R}''$}] (2) [right = 1cm of 1] {};
    \node (x1) [left = 0.5cm of 1] {};
    \node (y1) [above left = .7cm of 1] {};
    \node (z1) [below left = .7cm of 1] {};
    \node (x2) [right = 0.5cm of 2] {};
    \node (y2) [above right = .7cm of 2] {};
    \node (z2) [below right = .7cm of 2] {};
    \draw (x1) -- (1); \draw (y1) -- (1); \draw(z1) -- (1);
    \draw (x2) -- (2); \draw (y2) -- (2); \draw(z2) -- (2);
  \end{tikzpicture}
  .
\]

If $\chi(M) \neq 0$, we cannot directly map $\Gra^{\circlearrowleft}_{R}$ to $\OmPA^{*}(\FM_{M})$, as the Euler class in $\OmPA^{*}(M)$ would need to be the boundary of the image of the loop $e_{11} \in \Gra^{\circlearrowleft}_{R}(\underline{1})$.
We thus define a sub-\abbr{CDGA} which will map to $\OmPA^{*}(\FM_{M})$ whether $\chi(M)$ vanishes or not.

\begin{defin}
  \label{def.no-loops}
  For a given finite set $U$, let $\Gra_{R}(U)$ be the submodule of $\Gra_{R}^{\circlearrowleft}(U)$ spanned by graphs without loops.
\end{defin}

One has to be careful with the notation.
While $\Gra_{R}^{\circlearrowleft}(U) = R^{\otimes U} \otimes \Gra_{n}^{\circlearrowleft}$, it is not true that $\Gra_{R}(U) = R^{\otimes U} \otimes \Gra_{n}(U)$: in $\Gra_{n}(U)$, multiple edges are forbidden, whereas they are allowed in $\Gra_{R}(U)$.

\begin{prop}
  The space $\Gra_{R}(U)$ is a sub-\abbr{CDGA} of $\Gra_{R}^{\circlearrowleft}(U)$.
  If $\chi(M) = 0$ the collection $\Gra_{R}$ assembles to form a Hopf right $\Gra_{n}$-comodule.
\end{prop}
\begin{proof}
  Clearly, neither the splitting part of the differential nor the internal differential coming from $R$ can create new loops, nor can the product of two graphs without loops contain a loop, thus $\Gra_{R}(U)$ is indeed a sub-\abbr{CDGA} of $\Gra_{R}^{\circlearrowleft}(U)$.
  If $\chi(M) = 0$, the proof that $\Gra_{R}$ is a $\Gra_{n}$-comodule is almost identical to the proof of Proposition~\ref{prop.gra-tad-comod}, except that we need to use $\mu_{R}(\Delta_{R}) = 0$ to check that $d(\circ_{W}^{\vee}(e_{uv})) = \circ^{\vee}_{W}(d(e_{uv}))$ when $u,v \in W$.
\end{proof}

\subsection{The propagator}
\label{sec.propagator}

To define $\omega' : \Gra_{R} \to \OmPA(\FM_{M})$, we need a ``propagator'' $\varphi \in \OmPA^{n-1}(\FM_{M}(2))$, for which a reference is~\cite[Section 4]{CattaneoMneev2010}.

Recall from Equation~\eqref{eq:proj-fm-m} the projections $p_{u} : \FM_{M}(U) \to M$ and $p_{uv} : \FM_{M}(U) \to \FM_{M}(\underline{2})$.
Recall moreover the sphere bundle $p : \partial \FM_{M}(\underline{2}) \to M$ defined in Equation~\eqref{eq:sphere-bundle}, which is trivial when $M$ is framed, with the isomorphism $M \times S^{n-1} \xrightarrow{\circ_{1}} \partial \FM_{M}(\underline{2})$ from Equation~\eqref{eq:bundle-trivial}.
We denote by $(p_{1}, p_{2}) : \FM_{M}(\underline{2}) \to M \times M$ the product of the two canonical projections.

\begin{prop}[{\cite[Propositions~7 and~87]{CamposWillwacher2016}}]
  \label{def.propagator}
  There exists a form $\varphi \in \OmPA^{n-1}(\FM_{M}(\underline{2}))$ such that $\varphi^{21} = (-1)^{n} \varphi$, $d \varphi = (p_{1}, p_{2})^{*}((\sigma \otimes \sigma)(\Delta_{R}))$ and such that the restriction of $\varphi$ to $\partial \FM_{M}(2)$ is a global angular form, i.e.\ it is a volume form of $S^{n-1}$ when restricted to each fiber.
  When $M$ is framed one can moreover choose $\varphi|_{\partial \FM_{M}(2)} = 1 \times \vol_{S^{n-1}} \in \Omega^{n-1}_{PA}(M \times S^{n-1})$.
  This propagator can moreover be chosen to be a trivial form ((see Section~\ref{sec:semi-algebraic-sets}).
\end{prop}

The proofs of~\cite{CamposWillwacher2016} relies on earlier computations given in~\cite{CattaneoMneev2010}, where this propagator is studied in detail.
One can see from the proofs of~\cite[Section 4]{CattaneoMneev2010} that $d\varphi$ can in fact be chosen to be any pullback of a form cohomologous to the diagonal class $\Delta_{M} \in \OmPA^{n}(M \times M)$.
We will make further adjustments to the propagator $\varphi$ in Proposition~\ref{prop.phi-phi}.
Recall $p_{u}$, $p_{uv}$ from Equation~\eqref{eq:proj-fm-m}.

\begin{prop}
  \label{prop.zigzag-gra-r}
  There is a morphism of collections of \abbr{CDGA}s given by:
  \begin{align*}
    \omega' : \Gra_{R}
    & \to \OmPA(\FM_{M}),
    & \begin{cases}
      \bigotimes_{u \in U} x_{u} \in R^{\otimes U}
      & \mapsto \bigwedge_{u \in U} p_{u}^{*}(\sigma(x_{u})), \\
      e_{uv}
      & \mapsto p_{uv}^{*}(\varphi).
    \end{cases}
  \end{align*}
  Moreover, if $M$ is framed, then $\omega'$ defines a morphism of comodules, where $\omega' : \Gra_{n} \to \OmPA^{*}(\FM_{n})$ was defined in Section~\ref{sec.formality}:
  \[ (\Gra_{R}, \Gra_{n}) \xrightarrow{(\omega', \omega')} (\OmPA^{*}(\FM_{M}), \OmPA^{*}(\FM_{n})) \]
\end{prop}

\begin{proof}
  The property $d \varphi = (p_{1}, p_{2})^{*}((\sigma \otimes \sigma)(\Delta_{R}))$ shows that the map $\omega'$ preserves the differential.
  Let us now assume that $M$ is framed to prove that this is a morphism of right comodules.
  Cocomposition commutes with $\omega'$ on the generators coming from $A^{\otimes U}$, since the comodule structure of $\OmPA^{*}(\FM_{M})$ multiplies together forms that are pullbacks of forms on $M$:
  \[ \circ_{W}^{\vee}(p_{u}^{*}(x)) = \begin{cases}
      p_{u}^{*}(x) \otimes 1 & \text{if } u \not\in W; \\
      p_{*}^{*}(x) \otimes 1 & \text{if } u \in W.
    \end{cases} \]

  We now check the compatibility of the cocomposition $\circ_{W}^{\vee}$ with $\omega'$ on the generator $\omega_{uv}$, for some $W \subset U$.
  \begin{itemize}
  \item If one of $u,v$, or both, is not in $W$, then the equality $\circ_{W}^{\vee}(\omega'(e_{uv})) = (\omega' \otimes \omega')(\circ_{W}^{\vee}(e_{uv})).$ is clear.
  \item  Otherwise suppose $\{u,v\} \subset W$.
    We may assume that $U = W = \underline{2}$ (it suffices to pull back the result along $p_{uv}$ to get the general case), so that we are considering the insertion of an infinitesimal configuration $M \times \FM_{n}(\underline{2}) \to \FM_{M}(\underline{2})$.
    This insertion factors through the boundary $\partial \FM_{M}(\underline{2})$.
    We have (see Definition~\ref{def.propagator}):
    \[ \circ_{\underline{2}}^{\vee}(\varphi) = 1 \otimes \vol_{S^{n-1}} \in \OmPA^{*}(M) \otimes \OmPA^{*}(\FM_{n}(\underline{2})) = \OmPA^{*}(M) \otimes \OmPA^{*}(S^{n-1}).
    \]
    Going back to the general case, we find:
    \[ \circ_{W}^{\vee}(\omega'(e_{uv})) = \circ_{W}^{\vee}(p_{uv}^{*}(\varphi)) = 1 \otimes p_{uv}^{*}(\vol_{S^{n-1}}), \]
    which is indeed the image of $\circ_{W}^{\vee}(\omega_{uv}) = 1 \otimes \omega_{uv}$ by $\omega' \otimes \omega'$.
    \qedhere
\end{itemize}
\end{proof}

\subsection{Labeled graphs with internal and external vertices: \texorpdfstring{$\Tw\Gra_{R}$}{Tw Gra\_R}}
\label{sec.twisting-twgra_r}

The general framework of operadic twisting, recalled in Section~\ref{sec:operadic-twisting}, shows that to twist a right (co)module, one only needs to twist the (co)operad.
Since our cooperad is one-dimensional in arity zero, the comodule inherits a Hopf comodule structure too (Lemma~\ref{lem:tw-hopf-comod}).

\begin{defin}
  The \textbf{twisted labeled graph comodule} $\Tw\Gra^{\circlearrowleft}_{R}$ is a Hopf right $(\Tw\Gra^{\circlearrowleft}_{n})$-comodule obtained from $\Gra^{\circlearrowleft}_{R}$ by twisting with respect to the Maurer--Cartan element $\mu \in (\Gra^{\circlearrowleft}_{n})^{\vee}(\underline{2})$ of Section~\ref{sec.formality}.
\end{defin}

We now explicitly describe this comodule in terms of graphs.
The dg-module $\Tw \Gra^{\circlearrowleft}_{R} (U)$ is spanned by graphs with two kinds of vertices, external vertices corresponding to elements of $U$, and indistinguishable internal vertices (usually drawn in black).
The degree of an edge is $n-1$, the degree of an external vertex is $0$, while the degree of an internal vertex is $-n$.
All the vertices are labeled by elements of $R$, and their degree is added to the degree of the graph.

The Hopf structure glues two graphs along their external vertices, multiplying labels in the process.
The differential is a sum of three terms
\[ d = d_{R} + d_{\mathrm{split}} + d_{\mathrm{contr}}. \]
The first part is the internal differential coming from $R$, acting on each label separately.
The second part comes from $\Gra^{\circlearrowleft}_{R}$ and splits edges, multiplying by $\Delta_{R}$ the labels of the endpoints.
The third part is similar to the differential of $\Tw \Gra^{\circlearrowleft}_{n}$: it contracts all \textbf{contractible} edges, i.e.\ edges connecting an internal vertex to another vertex of either kind.
When an edge is contracted, the label of the resulting vertex is the product of the labels of of the endpoints of the former edge (see Figure~\ref{fig.diff-gra-n}).
This result comes from the twisting construction (see the definition in Equation~\eqref{eq:total-diff-tw-comod}).
For example, we have:

\begin{equation}
    \biggl(
    \begin{tikzpicture}[baseline, scale=.5]
      \node [extv] (e) {1};
      \node [intv, right = .5cm of 1, label = {\small $x$}] (i) {};
      \draw (e) -- (i);
    \end{tikzpicture}
    \biggr)
    \xmapsto{\; d \;}
    \biggl(
    \begin{tikzpicture}[baseline]
      \node [extv] (e) {1};
      \node [intv, right = .5cm of 1, label = {\small $d_{R}x$}] (i) {};
    \end{tikzpicture}
    \biggr)
    \pm
    \sum_{(\Delta_{R})}
    \biggl(
    \begin{tikzpicture}[baseline]
      \node [extv, label = {\small $\Delta'_{R}$}] (e) {1};
      \node [intv, right = .5cm of 1, label = {\small $x \Delta''_{R}$}] (i) {};
    \end{tikzpicture}
    \biggr)
    \pm
    \begin{tikzpicture}[baseline]
      \node [extv, label = {\small $x$}] (e) {1};
    \end{tikzpicture}
    \quad (x \in R).
  \label{eq:exa-diff-tw-gra-r}
\end{equation}

\begin{rmk}
  \label{rmk.tw-gra-r-dead-end}
  An edge connected to a univalent internal vertex is contractible in $\Tw\Gra^{\circlearrowleft}_{R}$, though this is not the case in $\Tw\Gra_{n}^{\circlearrowleft}$.
  Indeed, if we go back to the definition of the differential in a twisted comodule (Equation~\eqref{eq:total-diff-tw-comod}), we see that the Maurer--Cartan element $\mu$ (Equation~\eqref{fig.maurer-cartan-element}) only acts on the right of the graph.
  Therefore, there is no term to cancel out the contraction of such edges, as was the case in $\Tw\Gra_{n}$ (see the discussion in Section~\ref{sec.formality} about the differential).
  In Equation~\eqref{eq:exa-diff-tw-gra-r}, the only edge would not be considered as contractible in $\Tw\Gra_{n}$ if we forgot the labels, but it is in $\Tw\Gra_{R}$.
\end{rmk}

Finally, the comodule structure is similar to the cooperad structure of $\Tw\Gra^{\circlearrowleft}_{n}$: for $\Gamma \in \Gra^{\circlearrowleft}_{R}(U \sqcup I) \subset \Tw\Gra^{\circlearrowleft}_{R}(U)$, the cocomposition $\circ_{W}^{\vee}(\Gamma)$ is the sum over tensors of the type $\pm \Gamma_{U/W} \otimes \Gamma_{W}$, where $\Gamma_{U/W} \in \Gra^{\circlearrowleft}_{R}(U/W \sqcup J)$, $\Gamma_{W} \in \Gra_{n}(W \sqcup J')$, $J \sqcup J' = I$, and there exists a way of inserting $\Gamma_{W}$ in the vertex $*$ of $\Gamma_{U/W}$ and reconnecting edges to get $\Gamma$ back.
See the following example of cocomposition $\circ_{\{1\}}^{\vee} : \Tw \Gra_{R}(\underline{1}) \to \Tw\Gra_{R}(\underline{1}) \otimes \Tw\Gra_{R}(\underline{1})$, where $x,y \in R$:

\begin{equation*}
  \begin{gathered} \begin{tikzpicture}
      \node [extv, label={above right}:{\footnotesize $x$}] (1) {1};
      \node [intv, label={\footnotesize $y$}, above = 0.5cm of 1] (i) {};
      \draw (1) -- (i);
    \end{tikzpicture} \end{gathered}
  \;
  \xmapsto{\; \circ_{\{1\}}^{\vee} \;}
  \;
  \left(
  \begin{gathered} \begin{tikzpicture}
      \node [extv, label={above right}:{\footnotesize $x$}] (1) {$*$};
      \node [intv, label={\footnotesize $y$}, above = 0.5cm of 1] (i) {};
      \draw (1) -- (i);
    \end{tikzpicture} \end{gathered}
  \otimes
  \begin{gathered} \begin{tikzpicture}
      \node [extv] (1) {1};
    \end{tikzpicture} \end{gathered}
  \right)
  \pm
  \left(
  \begin{gathered} \begin{tikzpicture}
      \node [extv, label={\footnotesize $xy$}] (1) {$*$};
      \draw (1) to [loop right] (1);
    \end{tikzpicture} \end{gathered}
  \otimes
  \begin{gathered} \begin{tikzpicture}
      \node [extv] (1) {1};
      \node [intv, above = 0.5cm of 1] (i) {};
    \end{tikzpicture} \end{gathered}
  \right)
  \pm
  \left(
  \begin{gathered} \begin{tikzpicture}
      \node [extv, label={\footnotesize $xy$}] (1) {$*$};
    \end{tikzpicture} \end{gathered}
  \otimes
  \begin{gathered} \begin{tikzpicture}
      \node [extv] (1) {1};
      \node [intv, above = 0.5cm of 1] (i) {};
      \draw (1) -- (i);
    \end{tikzpicture} \end{gathered}
  \right)
\end{equation*}

\begin{lem}
  \label{lemma.trick}
  The subspace $\Tw\Gra_{R}(U) \subset \Tw\Gra^{\circlearrowleft}_{R}(U)$ spanned by graphs with no loops is a sub-\abbr{CDGA}.
\end{lem}
\begin{proof}
  It is clear that this defines a subalgebra.
  We need to check that it is preserved by the differential, i.e.\ that the differential cannot create new loops if there are none in a graph.
  This is clear for the internal differential coming from $R$ and for the splitting part of the differential.
  The contracting part of the differential could create a loop from a double edge.
  However for even $n$ multiple edges are zero for degree reasons, and for odd $n$ loops are zero because of the antisymmetry relation (see Remark~\ref{rmk.loop-multi}).
\end{proof}

Note that despite the notation, $\Tw\Gra_{R}$ is a priori not defined as the twisting of the $\Gra_{n}$-comodule $\Gra_{R}$: when $\chi(M) \neq 0$, the collection $\Gra_{R}$ is not even a $\Gra_{n}$-comodule.
However, the following proposition is clear and shows that we can get away with this abuse of notation:

\begin{prop}
  If $\chi(M) = 0$, then $\Tw\Gra_{R}$ assembles to a right Hopf $(\Tw\Gra_{n})$-comodule, isomorphic to the twisting of the right Hopf $\Gra_{n}$-comodule $\Gra_{R}$ of Definition~\ref{def.no-loops}.
  \qed
\end{prop}

\begin{rmk}
  \label{rmk:detour}
  We could have defined the algebra $\Tw\Gra_{R}$ explicitly in terms of graphs, and defined the differential $d$ using an ad-hoc formula. The difficult part would have then been to check that $d^{2} = 0$ (involving difficult signs), which is a consequence of the general operadic twisting framework.
\end{rmk}

\subsection{The map \texorpdfstring{$\omega : \Tw \Gra_{R} \to \OmPA^{*}(\FM_{M})$}{w : Tw Gra\_R -> Omega(M)}}
\label{sec:map-omega}

This section is dedicated to the proof of the following proposition.

\begin{prop}
  \label{prop.mph-tw-om}
  There is a morphism of collections of \abbr{CDGA}s $\omega : \Tw\Gra_{R} \to \OmPA^{*}(\FM_{M})$ extending $\omega'$, given on a graph $\Gamma \in \Gra_{R}(U \sqcup I) \subset \Tw\Gra_{R}(U)$ by:
  \[ \omega(\Gamma) \coloneqq \int_{p_{U} : \FM_{M}(U \sqcup I) \to \FM_{M}(U)} \omega'(\Gamma) = (p_{U})_{*}(\omega'(\Gamma)).
  \]
  Moreover, if $M$ is framed, then this defines a morphism of Hopf right comodules:
  \[ (\omega, \omega) : (\Tw\Gra_{R}, \Tw\Gra_{n}) \to (\OmPA^{*}(\FM_{M}), \OmPA^{*}(\FM_{n})). \]
\end{prop}

Recall that in general, it is not possible to consider integrals along fibers of arbitrary \abbr{PA} forms, see~\cite[Section~9.4]{HardtLambrechtsTurchinVolic2011}.
However, here, the image of $\sigma$ is included in the sub-\abbr{CDGA} of trivial forms in $\OmPA^{*}(M)$, and the propagator is a trivial form (see Proposition~\ref{def.propagator}), therefore the integral $(p_{U})_{*}(\omega'(\Gamma))$ exists.
  
The proof of the compatibility with the Hopf structure and, in the framed case, the comodule structure, is formally similar to the proof of the same facts about $\omega : \Tw\Gra_{n} \to \OmPA^{*}(\FM_{n})$.
We refer to~\cite[Sections 9.2, 9.5]{LambrechtsVolic2014}.
The proof is exactly the same proof, but writing $\FM_{M}$ or $\FM_{n}$ instead of $C[-]$ and $\varphi$ instead of $\vol_{S^{n-1}}$ in every relevant sentence, and recalling that when $M$ is framed, we choose $\varphi$ such that $\circ_{\underline{2}}^{\vee}(\varphi) = 1 \otimes \vol_{S^{n-1}}$.

The proof that $\omega$ is a chain map is different albeit similar.
We recall Stokes' formula for integrals along fibers of semi-algebraic bundles.
If $\pi : E \to B$ is a semi-algebraic bundle, the fiberwise boundary $\pi^{\partial} : E^{\partial} \to B$ is the bundle with
\[ E^{\partial} \coloneqq \bigcup_{b \in B} \partial \pi^{-1}(b). \]

\begin{rmk}
  The space $E^{\partial}$ is neither $\partial E$ nor $\bigcup_{b \in B} \pi^{-1}(b) \cap \partial E$ in general.
  (Consider for example the projection on the first coordinate $[0,1]^{\times 2} \to [0,1]$.)
\end{rmk}

Stokes' formula, in the semi-algebraic context, is~\cite[Proposition 8.12]{HardtLambrechtsTurchinVolic2011}:
\[ d \left( \int_{\pi : E \to B} \alpha \right) = \int_{\pi : E \to B} d\alpha \pm \int_{\pi^{\partial} : E^{\partial} \to B} \alpha_{| E^{\partial}}. \]

If we apply this formula to compute $d \omega(\Gamma)$, we find that the first term is:
\begin{equation}
  \label{eq.int-d-omega}
  \int_{p_{U}} d \omega'(\Gamma) = \int_{p_{U}} \omega'(d_{R} \Gamma + d_{\mathrm{split}} \Gamma) = \omega(d_{R} \Gamma + d_{\mathrm{split}} \Gamma),
\end{equation}
since $\omega'$ was a chain map.
It thus remain to check that the second term satisfies:
\[ \int_{p_{U}^{\partial} : \FM_{M}^{\partial}(U \sqcup I) \to \FM_{M}(U)} \omega'(\Gamma) = \int_{p_{U}} \omega'(d_{\mathrm{contr}} \Gamma) = \omega(d_{\mathrm{contr}} \Gamma). \]

The fiberwise boundary of the projection $p_{U} : \FM_{n}(U \sqcup I) \to \FM_{n}(U)$ is rather complex~\cite[Section 5.7]{LambrechtsVolic2014}, essentially due to the quotient by the affine group in the definition of $\FM_{n}$ which lowers dimensions.
We will not repeat its explicit decomposition into cells as we do not need it here.

The fiberwise boundary of $p_{U} : \FM_{M}(U \sqcup I) \to \FM_{M}(U)$ is simpler.
Our definitions mimick the description of~\cite[Section~5.7]{LambrechtsVolic2014}.
Let $V = U \sqcup I$.
The interior of $\FM_{M}(U)$ is the space $\Conf_{U}(M)$, and thus $\FM^{\partial}_{M}(V)$ is the closure of $(\partial \FM_{M}(V)) \cap \pi^{-1}(\Conf_{U}(M))$.
Let the set of ``boundary faces'' be given by:
\[ \mathcal{BF}_{M}(V,U) = \{ W \subset V \mid \# W \geq 2 \text{ and } \# W \cap U \leq 1 \}. \]
This set indexes the strata of the fiberwise boundary of $p_{U}$.
The idea is that a configuration is in the fiberwise boundary iff it is obtained by an insertion map $\circ_{W}$ with $W \in \mathcal{BF}_{M}(V,U)$.
In the description of $\FM_{n}^{\partial}(V)$, similar boundary faces, denoted $\mathcal{BF}(V,U)$, appear.
But there, there was an additional part which corresponds to $U \subset W$.
Unlike the case of $\FM_{n}$, for $\FM_{M}$ the image of $p_{U} (- \circ_{W} -)$ is always included in the boundary of $\FM_{M}(U)$ when $U \subset W$.
We follow a pattern similar to the one used in the proof of~\cite[Proposition 5.7.1]{LambrechtsVolic2014}.

\begin{lem}
  The subspace $\FM_{M}^{\partial}(V) \subset \FM_{M}(V)$ is equal to:
  \[ \bigcup_{\mathclap{W \in \mathcal{BF}_{M}(V,U)}} \im \bigl( \circ_{W} : \FM_{M}(V/W) \times \FM_{n}(W) \to \FM_{M}(V) \bigr). \]
\end{lem}

\begin{proof}
  Let $\operatorname{cls}$ denote the closure operator.
  Since $\Conf_{U}(M)$ is the interior of $\FM_{M}(U)$ and $p : \FM_{M}(V) \to \FM_{M}(U)$ is a bundle, it follows that the fiberwise boundary $\FM_{M}^{\partial}$ is obtained as the closure of the preimage of the interior (see the corresponding statement in the proof of~\cite[Proposition~5.7.1]{LambrechtsVolic2014}), i.e.:
  \begin{align*}
    \FM_{M}^{\partial}(V)
    & = \operatorname{cls} \left( \FM_{M}^{\partial}(V) \cap p^{-1}(\Conf_{U}(M)) \right) \\
    & = \operatorname{cls} \left( \partial \FM_{M}(V) \cap p^{-1}(\Conf_{U}(M)) \right).
  \end{align*}

  The boundary $\partial \FM_{M}(V)$ is the union of the subsets $\im(\circ_{W})$ for $\# W \geq 2$ (note that the case $W = V$ is included, unlike for $\FM_{n}$).
  If $\# W \cap U \geq 2$, which is equivalent to $W \not\in \mathcal{BF}_{M}(V,U)$, then $\im(p_{U}(- \circ_{W} -)) \subset \partial \FM_{M}(U)$, because if a configuration belongs to this image then at least two points of $U$ are infinitesimally close.
  Therefore:
  \begin{align*}
    \operatorname{cls} \bigl( \partial \FM_{M}(V) \cap p^{-1}(\Conf_{U}(M)) \bigr)
    & = \operatorname{cls} \biggl( \bigcup_{\hspace{-1em} \mathrlap{\# W \geq 2}} \im (\circ_{W}) \cap p^{-1}(\Conf_{U}(M)) \biggr) \\
    & = \operatorname{cls} \biggl( \bigcup_{\hspace{-1em} \mathrlap{\# W \in \mathcal{BF}_{M}(V,U)}} \im (\circ_{W}) \cap p^{-1}(\Conf_{U}(M))\biggr) \\
    & = \bigcup_{\mathclap{W \in \mathcal{BF}_{M}(V,U)}} \operatorname{cls} \bigl( \im(\circ_{W}) \cap p^{-1}(\Conf_{U}(M)) \bigr) \\
    & = \bigcup_{\mathclap{W \in \mathcal{BF}_{M}(V,U)}} \im(\circ_{W}).
      \qedhere
  \end{align*}
\end{proof}

\begin{lem}
  \label{lemma.int-partial-omega}
  For a given graph $\Gamma \in \Tw\Gra_{R}(U)$, the integral over the fiberwise boundary is given by:
  \[ \int_{p_{U}^{\partial}} \omega'(\Gamma)_{| \FM_{M}^{\partial}(V)} = \omega(d_{\mathrm{contr}} \Gamma). \]
\end{lem}

\begin{proof}
  The maps $\circ_{W} : \FM_{M}(V/W) \times \FM_{n}(W) \to \FM_{M}(V)$ are smooth injective maps and their domains are compact, thus they are homeomorphisms onto their images.
  Recall that $\# W \geq 2$ for $W \in \mathcal{BF}_{M}(V,U)$; hence $\dim \FM_{n}(W) = n \#W - n - 1$.
  The dimension of the image of $\circ_{W}$ is then:
  \begin{align*}
    \dim \im(\circ_{W})
    & = \dim \FM_{M}(V/W) + \dim \FM_{n}(W) \\
    & = n \# (V/W) + (n \#W - n - 1) \\
    & = n \#V - 1,
  \end{align*}
  i.e.\ the image is of codimension $1$ in $\FM_{M}(V)$.
  It is also easy to check that if $W \neq W'$, then $\im(\circ_{W}) \cap \im(\circ_{W'})$ is of codimension strictly bigger than $1$.

  We now fix $W \in \mathcal{BF}_{M}(V,U)$.
  Since $\# W \cap U \leq 1$, the composition $U \subset V \to V/W$ is injective and identifies $U$ with a subset of $V/W$.
  There is then a forgetful map $p_{U}' : \FM_{M}(V/W) \to \FM_{M}(U)$.
  We then have a commutative diagram:
  \begin{equation}
    \begin{tikzcd}
      \FM_{M}(V/W) \times \FM_{n}(W) \dar{\circ_{W}} \rar{p_{1}} & \FM_{M}(V/W) \dar{p_{U}'} \\
      \FM_{M}(V) \rar{p_{U}} & \FM_{M}(U)
    \end{tikzcd}
    .
    \label{eq:6}
  \end{equation}

  It follows that $p_{U}(- \circ_{W} -) = p'_{U} \circ p_{1}$ is the composite of two semi-algebraic bundles, hence it is a semi-algebraic bundle itself~\cite[Proposition 8.5]{HardtLambrechtsTurchinVolic2011}.
  Combined with the fact about codimensions above, we can therefore apply the summation formula~\cite[Proposition 8.11]{HardtLambrechtsTurchinVolic2011}:
  \begin{equation}
    \label{eq.sum-chain-map}
    \int_{p_{U}^{\partial}} \omega'(\Gamma) = \sum_{W \in \mathcal{BF}_{M}(V,U)} \int_{p_{U}(- \circ_{W} -)} \omega'(\Gamma)_{| \FM_{M}(V/W) \times \FM_{n}(W)}.
  \end{equation}

  Now we can directly adapt the proof of Lambrechts and Volić.
  For a fixed $W$, by~\cite[Proposition 8.13]{HardtLambrechtsTurchinVolic2011}, the corresponding summand is equal to $\pm \omega(\Gamma_{V/W}) \cdot \int_{\FM_{n}(W)} \omega'(\Gamma_{W})$, where
  \begin{itemize}
  \item $\Gamma_{V/W} \in \Tw\Gra_{R}(U)$ is the graph with $W$ collapsed to a vertex and $U \hookrightarrow V/W$ is identified with its image;
  \item $\Gamma_{W} \in \Tw\Gra_{n}(W)$ is the full subgraph of $\Gamma$ with vertices $W$ and the labels removed.
  \end{itemize}

  The vanishing lemmas in the proof of Lambrechts and Volić then imply that the integral $\int_{\FM_{n}(W)} \omega'(\Gamma_{W})$ is zero unless $\Gamma_{W}$ is the graph with exactly two vertices and one edge, in which case the integral is equal to $1$.
  In this case, $\Gamma_{V/W}$ is the graph $\Gamma$ with one edge connecting an internal vertex to some other vertex collapsed.
  The sum runs over all such edges, and dealing with signs carefully we see that Equation~\eqref{eq.sum-chain-map} is precisely equal to $\omega(d_{\mathrm{contr}} \Gamma)$.
\end{proof}

We can now finish proving Proposition~\ref{prop.mph-tw-om}.
We combine Equation~\eqref{eq.int-d-omega} and Lemma~\ref{lemma.int-partial-omega}, and apply Stokes' formula to $d \omega(\Gamma)$ to show that it is equal to $\omega(d \Gamma) = \omega(d_{R} \Gamma + d_{\mathrm{split}} \Gamma) + \omega(d_{\mathrm{contr}} \Gamma)$.

\subsection{Reduced labeled graphs: \texorpdfstring{$\Graphs_R$}{Graphs\_R}}
\label{sec.reduct-graphs_r}

The last step in the construction of $\Graphs_{R}$ is the reduction of $\Tw\Gra_{R}$ so that it has the right cohomology.
We borrow the terminology of Campos--Willwacher~\cite{CamposWillwacher2016} for the next two definitions.

\begin{defin}
  \label{def:fgc}
  The \textbf{full graph complex} $\fGC_{R}$ is the \abbr{CDGA} $\Tw\Gra_{R}(\varnothing)$.
  It consists of labeled graphs with only internal vertices, and the product is disjoint union of graphs.
\end{defin}

\begin{rmk}
  \label{rmk:full-graph}
  The adjective ``full'' refers to the fact that graphs are possibly disconnected and have vertices of any valence in $\fGC_{R}$.
\end{rmk}

As an algebra, $\fGC_{R}$ is free and generated by connected graphs.
In general we will call \textbf{internal components} the connected components of a graph that only contain internal vertices.
The full graph complex naturally acts on $\Tw\Gra_{R}(U)$ by adding extra internal components.

\begin{defin}
  \label{def:part-func-phi}
  The \textbf{partition function} $\Zphi : \fGC_{R} \to \R$ is the restriction of $\omega : \Tw\Gra_{R} \to \OmPA^{*}(\FM_{M})$ to $\fGC_{R} = \Tw\Gra_{R}(\varnothing) \to \OmPA^{*}(\FM_{M}(\varnothing)) = \OmPA^{*}(\mathrm{pt}) = \R$.
\end{defin}

\begin{rmk}
  \label{rmk:name-Z}
  The expression ``partition function'' comes from the mathematical physics literature, more specifically Chern--Simons invariant theory, where it refers to the partition function of a quantum field theory.
\end{rmk}

By the double-pushforward formula~\cite[Proposition~8.13]{HardtLambrechtsTurchinVolic2011} and Fubini's theorem~\cite[Proposition~8.15]{HardtLambrechtsTurchinVolic2011}, $\Zphi$ is an algebra morphism and
\begin{equation}
  \label{eq.z-r-prod}
  \forall \gamma \in \fGC_{R}, \forall \Gamma \in \Tw\Gra_{R}(U), \; \omega(\gamma \cdot \Gamma) = \Zphi(\gamma) \cdot \omega(\Gamma).
\end{equation}

\begin{defin}
  \label{def:graphs-phi}
  Let $\R_{\varphi}$ be the $\fGC_{R}$-module of dimension $1$ induced by $\Zphi : \fGC_{R} \to \R$.
  The \textbf{reduced graph comodule} $\Graphs^{\varphi}_{R}$ is the tensor product:
  \[ \Graphs^{\varphi}_{R}(U) \coloneqq \R_{\varphi} \otimes_{\fGC_{R}} \Tw\Gra_{R}(U). \]
\end{defin}

In other words, a graph of the type $\Gamma \sqcup \gamma$ containing an internal component $\gamma \in \fGC_{R}$ is identified with $\Zphi(\gamma) \cdot \Gamma$.
It is spanned by representative classes of graphs with no internal connected component; we call such graphs \textbf{reduced}.
The notation is meant to evoke the fact that $\Graphs^{\varphi}_{R}$ depends on the choice of the propagator $\varphi$, unlike the collection $\Graphs^{\varepsilon}_{R}$ that will appear in Section~\ref{sec.constr-graphs_n-to}.

\begin{prop}
  \label{prop.mph-1}
  The map $\omega : \Tw\Gra_{R}(U) \to \OmPA^{*}(\FM_{M}(U))$ defined in Proposition~\ref{prop.mph-tw-om} factors through the quotient defining $\Graphs^{\varphi}_{R}$.

  If $\chi(M) = 0$, then $\Graphs^{\varphi}_{R}$ forms a Hopf right $\Graphs_{n}$-comodule.
  If moreover $M$ is framed, then the map $\omega$ defines a Hopf right comodule morphism.
\end{prop}

\begin{proof}
  Equation~\eqref{eq.z-r-prod} immediately implies that $\omega$ factors through the quotient.

  The vanishing lemmas shows that if $\Gamma \in \Tw\Gra_{n}(U)$ has internal components, then $\omega(\Gamma)$ vanishes by~\cite[Lemma~9.3.7]{LambrechtsVolic2014}, so it is straightforward to check that if $\chi(M) = 0$, then $\Graphs^{\varphi}_{R}$ becomes a Hopf right comodule over $\Graphs_{n}$.
  It is also clear that for $M$ framed, the quotient map $\omega$ remains a Hopf right comodule morphism.
\end{proof}

\begin{prop}[{\cite[Lemma 3]{CattaneoMneev2010}}]
  \label{prop.phi-phi}
  The propagator $\varphi$ can be chosen such that the following additional property (P4) holds:
  \[ \int_{p_{1} : \FM_{M}(\underline{2}) \to  \FM_{M}(\underline{1}) = M} p_{2}^{*}(\sigma(x)) \wedge \varphi = 0, \; \forall x \in R. \tag{P4} \]
\end{prop}

From now on and until the end, we assume that $\varphi$ satisfies (P4).

\begin{rmk}
  \label{rmk.pfive}
  The additional property (P5) of the paper mentioned above would be helpful in order to get a direct morphism $\Graphs^{\varphi}_{R} \to \GG{A}$, because then the partition function would vanish on all connected graphs with at least two vertices. However we run into difficulties when trying to adapt the proof in the setting of \abbr{PA} forms, mainly due to the lack of an operator $d_{M}$ acting on $\OmPA^{*}(M \times N)$ differentiating ``only in the first slot''.
\end{rmk}

\begin{coro}
  \label{cor.weak-vanishing}
  The morphism $\omega$ vanishes on graphs containing univalent internal vertices.
\end{coro}
\begin{proof}
  Let $\Gamma \in \Gra_{R}(U \sqcup I) \subset \Tw\Gra_{R}(U)$ be a graph with a univalent internal vertex $u \in I$, labeled by $x$, and let $v$ be the only vertex connected to $u$.
  Let $\tilde\Gamma$ be the full subgraph of $\Gamma$ on the set of vertices $U \sqcup I \setminus \{u\}$.
  Then using~\cite[Propositions 8.10 and 8.15]{HardtLambrechtsTurchinVolic2011} (in a way similar to the end of the proof of~\cite[Lemma 9.3.8]{LambrechtsVolic2014}), we find:
  \begin{align*}
    \omega(\Gamma)
    & = \int_{\FM_{M}(U \sqcup I) \to \FM_{M}(U)} \omega'(\Gamma) \\
    & = \int_{\FM_{M}(U \sqcup I) \to \FM_{M}(U)} \omega'(\tilde\Gamma) p_{uv}^{*}(\varphi) p_{u}^{*}(\sigma(x)) \\
    & = \int_{\FM_{M}(U \sqcup I \setminus \{u\}) \to \FM_{M}(U)} \omega'(\tilde\Gamma) \wedge p_{v}^{*} \left( \int_{\FM_{M}(\{u,v\}) \to \FM_{M}(\{v\})} p_{uv}^{*}(\varphi) p_{u}^{*}(\sigma(x)) \right),
  \end{align*}
  which vanishes by (P4) in Proposition~\ref{prop.phi-phi}.
\end{proof}

Almost everything we have done so far works in full generality.
We now prove a fact which sets a class of manifolds apart.

\begin{prop}
  \label{prop.strong-vanishing}
  Assume that $M$ is simply connected and that $\dim M \ge 4$.
  Then the partition function $\Zphi$ vanishes on any connected graph with no bivalent vertices labeled by $1_{R}$ and containing at least two vertices.
\end{prop}

\begin{rmk}
  If $\gamma \in \fGC_{R}$ has only one vertex, labeled by $x$, then $\Zphi(\gamma) = \int_{M} \sigma(x)$ which can be nonzero.
\end{rmk}

\begin{proof}
  Let $\gamma \in \fGC_{R}$ be a connected graph with at least two vertices and no bivalent vertices labeled by $1_{R}$.
  By Corollary~\ref{cor.weak-vanishing}, we can assume that all the vertices of $\gamma$ are at least bivalent.
  By hypothesis, if a vertex is bivalent then it is labeled by an element of $R^{> 0} = R^{\geq 2}$.

  Let $k = i + j$ be the number of vertices of $\gamma$, with $i$ vertices that are at least trivalent and $j$ vertices that are bivalent and labeled by $R^{\geq 2}$.
  It follows that $\gamma$ has at least $\frac{1}{2} (3i + 2j)$ edges, all of degree $n-1$.
  Since bivalent vertices are labeled by $R^{\geq 2}$, their labels contribute by at least $2j$ to the degree of $\gamma$.
  The (internal) vertices contribute by $-kn$ to the degree, and the other labels have a nonnegative contribution.
  Thus:
  \begin{align*}
    \deg \gamma
    & \geq \biggl( \frac{3}{2} i + j \biggr) (n-1) + 2j - kn
    = \biggl( \frac{3}{2} k - \frac{3}{2}j + j \biggr) (n-1) + 2j - kn \\
    & = \frac{1}{2} \bigl( k (n-3) - j (n-5) \bigr).
  \end{align*}
  This last number is always positive for $0 \leq j \leq k$: it is an affine function of $j$, and it is positive when $j = 0$ and $j = k$ (recall that $n \geq 4$).
  The degree of $\gamma \in \fGC_{R}$ must be zero for the integral defining $\Zphi(\gamma)$ to be the integral of a top form of $\FM_{M}(\underline{k})$ and hence possibly nonzero.
  But by the above computation, $\deg \gamma > 0 \implies \Zphi(\gamma) = 0$.
 \end{proof}
 
\begin{rmk}
  When $n = 3$, the manifold $M$ is the $3$-sphere $S^{3}$ by Perelman's proof of the Poincaré conjecture~\cite{Perelman2002,Perelman2003}.
  The partition function $\Zphi$ is conjectured to be trivial on $S^{3}$ for a proper choice of framing, thus bypassing the need for the above degree counting argument.
  See also Proposition~\ref{prop.result-s3}.
\end{rmk}

We will also need the following technical property of $\fGC_{R}$.
\begin{lem}
  \label{lem:fgc-cofibrant}
  The \abbr{CDGA} $\fGC_{R}$ is cofibrant.
\end{lem}
\begin{proof}
  We filter $\fGC_{R}$ by the number of edges, defining $F_{s}\fGC_{R}$ to be the submodule of $\fGC_{R}$ spanned by graphs of $\gamma$ such that all the connected components $\gamma$ have at most $s$ edges.
  The differential of $\fGC_{R}$ can only decrease ($d_{\mathrm{split}}$ and $d_{\mathrm{contr}}$) or leave constant ($d_{R}$) the number of edges.
  Moreover $F_{s}\fGC_{R}$ is clearly stable under products (disjoint unions), hence $F_{s}\fGC_{R}$ is a sub-\abbr{CDGA} of $\fGC_{R}$.
  It is also clear that $\fGC_{R} = \operatorname{colim}_{s} F_{s}\fGC_{R}$.
  We will prove that $F_{0}\fGC_{R}$ is cofibrant, and that each $F_{s}\fGC_{R} \subset F_{s+1}\fGC_{R}$ is a cofibration.

  The \abbr{CDGA} $F_{0}\fGC_{R}$ is the free \abbr{CDGA} on graphs with a single vertex labeled by $R$.
  In other words, $F_{0}\fGC_{R} = S(R, d_{R})$ is the free symmetric algebra on the dg-module $R$, and any free symmetric algebra on a dg-module is cofibrant.

  Let us now show that $F_{s}\fGC_{R} \subset F_{s+1}\fGC_{R}$ is a cofibration for any $s \ge 0$.
  We will show that it is in fact a ``relative Sullivan algebra''~\cite[Section~14]{FelixHalperinThomas2001}.
  As a \abbr{CDGA}, we have $F_{s+1}\fGC_{R} = (F_{s}\fGC_{R} \otimes S(V_{s+1}), d)$, where $V_{s+1}$ is the graded module of connected graphs with exactly $s+1$ edges.
  Let us now show the Sullivan condition.

  Recall that $R$ is obtained from the minimal model of $M$ by a relative Sullivan extension, hence it is itself a Sullivan algebra~\cite[Section~12]{FelixHalperinThomas2001}.
  In other words, $R = (S(W), d)$ where $W$ is increasingly and exhaustively filtered by $W(-1) = 0 \subset \dots \subset W(t) \subset \dots \subset W$ such that $d(W(t)) \subset S(W(t-1))$.
  This induces a filtration on $R$ by defining $R(t) \coloneqq \bigoplus_{t_{1} + \dots + t_{r} = t} \bigl( V(t_{1}) \otimes \dots \otimes V(t_{r}) \bigr)_{\Sigma_{r}}$.

  This in turns induces an increasing and exhaustive filtration on $V_{s+1}$ by submodules $V_{s+1}(t)$ as follows.
  A connected graph $\gamma \in V_{s+1}$ is in $V_{s+1}(t)$ if each label $x_{i} \in R$ of a vertex of $\gamma$ belongs to the filtration $R(t_{i})$ such that $\sum t_{i} = t$.
  It is then immediate to check that $d(V_{s+1}(t+1)) \subset V_{s} \otimes S(V_{s+1}(t))$.
  Indeed, if $\gamma \in V_{s+1}(t+1)$, then $d_{\mathrm{split}}\gamma$ and $d_{\mathrm{contr}}\gamma \in V_{s}$, because both strictly decrease the number of edges.
  And $d_{R}\gamma \in V_{s+1}(t)$ because the internal differential of $R$ decreases the filtration of $R$.
\end{proof}

\section{From the model to forms via graphs}
\label{sec.proof-theorem}

In this section we connect $\GG{A}$ to $\OmPA^{*}(\FM_{M})$ and we prove that the connecting morphisms are quasi-isomorphisms.
We assume that $M$ is a simply connected closed smooth manifold with $\dim M \ge 4$ (see Proposition~\ref{prop.strong-vanishing}).

\subsection{Construction of the morphism to \texorpdfstring{$\GG{A}$}{G\_A}}
\label{sec.constr-graphs_n-to}

\begin{prop}
  \label{prop.gra-r-g-a}
  For each finite set $U$, there is a \abbr{CDGA} morphism $\rho'_{*} : \Gra_{R}(U) \to \GG{A}(U)$ given by $\rho$ on the $R^{\otimes U}$ factor and sending the generators $e_{uv}$ to $\omega_{uv}$ on the $\Gra_{n}$ factor.
  When $\chi(M) = 0$, this defines a Hopf right comodule morphism $(\Gra_{R}, \Gra_{n}) \to (\GG{A}, \enV)$.
  \qed
\end{prop}

If we could find a propagator for which property (P5) held (see Remark~\ref{rmk.pfive}), then we could just send all graphs containing internal vertices to zero and obtain an extension $\Graphs^{\varphi}_{R} \to \GG{A}$.
Since we cannot assume that (P5) holds, the definition of the extension is more complex.
However we still have Proposition~\ref{prop.strong-vanishing}, and homotopically speaking, graphs with bivalent vertices are irrelevant.

\begin{defin}
  Let $\fGC^{0}_{R}$ be the quotient of $\fGC_{R}$ defined by identifying a disconnected vertex labeled by $x$ with the number $\varepsilon_{A}(\rho(x))$.
\end{defin}

\begin{lem}
  \label{lem.fgczero}
  The subspace $I \subset \fGC^{0}_{R}$ spanned by graphs with at least one univalent vertex, or at least one bivalent vertex labeled by $1_{R}$, or at least one label in $\ker(\rho : R \to A)$, is a \abbr{CDGA} ideal.
\end{lem}

\begin{proof}
  It is clear that $I$ is an algebra ideal.
  Let us prove that it is a differential ideal.
  If one of the labels of $\Gamma$ is in $\ker \rho$, then so do all the summands of $d \Gamma$, because $\ker \rho$ is a \abbr{CDGA} ideal of $R$.

  If $\Gamma$ contains a bivalent vertex $u$ labeled by $1_{R}$, then so does $d_{R} \Gamma$. In $d_{\mathrm{split}} \Gamma$, splitting one of the two edges connected to $u$ produces a univalent vertex and hence vanishes in $\fGC_{R}^{0}$ because the label is $1_{R}$.
  In $d_{\mathrm{contr}} \Gamma$, the contraction of the two edges connected to $u$ cancel each other.

  Finally let us prove that if $\Gamma$ has a univalent vertex $u$, then $d\Gamma$ lies in $I$.
  It is clear that $d_{R} \Gamma \in I$.
  Contracting or splitting the only edge connected to the univalent vertex could remove the univalent vertex.
  Let us prove that these two summands cancel each other up to $\ker \rho$.
  
  It is helpful to consider the case pictured in Equation~\eqref{eq:exa-diff-tw-gra-r}.
  Let $y$ be the label of the univalent vertex $u$, and let $x$ be the label of the only vertex incident to $u$.
  Contracting the edge yields a new vertex labeled by $xy$.
  Due to the definition of $\fGC^{0}_{R}$, splitting the edge yields a new vertex labeled by $\alpha \coloneqq \sum_{(\Delta_{R})} \varepsilon(\rho(x \Delta''_{R})) y \Delta'_{R}$.
  We thus have $\rho(\alpha) = \rho(x) \cdot \sum_{(\Delta_{A})} \pm \varepsilon_{A}(\rho(y) \Delta_{A}'') \Delta_{A}'$.

  It is a standard property of the diagonal class that $\sum_{(\Delta_{A})} \pm \varepsilon_{A}(a \Delta_{A}'') \Delta_{A}' = a$ for all $a \in A$ (this property is a direct consequence of the definition in Equation~\eqref{eq.def-delta}).
  Applied to $a = \rho(y)$, it follows from the previous equation that $\rho(\alpha) = \pm \rho(xy)$; examining the signs, this summand cancels from the summand that comes from contracting the edge.
\end{proof}

\begin{defin}
  The algebra $\fGC'_{R}$ is the quotient of $\fGC^{0}_{R}$ by the ideal $I$.
\end{defin}

Note that $\fGC'_{R}$ is also free as an algebra, with generators given by connected graphs with no isolated vertices, nor univalent vertices, nor bivalent vertices labeled by $1_{R}$, and where the labels lie in $R / \ker(\rho) = A$.

\begin{defin}
  \label{def:floop}
  A circular graph is a graph in the shape of a circle and where all vertices are labeled by $1_{R}$, i.e.\ graphs of the type $e_{12} e_{23} \dots e_{(k-1)k} e_{k1}$.
  Let $\fLoop_{R} \subset \fGC^{0}_{R}$ be the submodule spanned by graphs whose connected components either have univalent vertices or are equal to a circular graphs.
\end{defin}

\begin{lem}
  The submodule $\fLoop_{R}$ is a sub-\abbr{CDGA} of $\fGC_{R}^{0}$.
\end{lem}
\begin{proof}
  The submodule $\fLoop_{R}$ is stable under products (disjoint union) by definition, so we just need to check that it is stable under the differential.
  Thanks to the proof of Lemma~\ref{lem.fgczero}, in $\fGC_{R}^{0}$, if a graph contains a univalent vertex, then so do all the summands of its differential.
  On a circular graph, the internal differential of $R$ vanish, because all labels are equal to $1_{R}$.
  Contracting an edge in a circular graph yields another circular graph, and splitting an edge yields a graph with univalent vertices, which belongs to $\fLoop_{R}$.
\end{proof}

\begin{prop}
  \label{prop.fgcr-cofib}
  The sequence $\fLoop_{R} \to \fGC^{0}_{R} \to \fGC'_{R}$ is a homotopy cofiber sequence of \abbr{CDGA}s.
\end{prop}

\begin{proof}
  The \abbr{CDGA} $\fGC^{0}_{R}$ is freely generated by connected labeled graphs with at least two vertices.
  It is a quasi-free extension of $\fLoop_{R}$ by the algebra generated by graphs that are not circular and that do not contain any univalent vertices.
  The homotopy cofiber of the inclusion $\fLoop_{R} \to \fGC^{0}_{R}$ is this algebra $\fGC''_{R}$, generated by graphs that are not circular and do not contain any univalent vertices, together with a differential induced by the quotient $\fGC^{0}_{R} / (\fLoop_{R})$.

  Let us note that the quotient map $\fGC^{0}_{R} \to \fGC'_{R} = \fGC^{0}_{R} / I$ vanishes on $\fLoop_{R}$, because $\fLoop_{R}$ is included in $R$.
  Thus we have a diagram:
  \[
    \begin{tikzcd}[row sep = small, column sep = small]
      0 \ar[r] & \fLoop_{R} \ar[d,hook] \ar[r,hook] & \fGC_{R}^{0} \ar[r,two heads] \ar[d, equal] & \fGC''_{R} \coloneqq \fGC^{0}_{R} / \fLoop_{R} \ar[r] \ar[d, two heads] & 0
      \\
      0 \ar[r] & I \ar[r,hook] & \fGC_{R}^{0} \ar[r,two heads] & \fGC'_{R} \coloneqq \fGC^{0}_{R} / I \ar[r] & 0
    \end{tikzcd}
  \]
  Let us prove that the morphism $\fGC''_{R} \to \fGC'_{R}$ is a quasi-isomorphism.
  Define an increasing filtration on both algebras by letting $F_{s}\fGC'_{R}$ (resp.\ $F_{s}\fGC''_{R}$) be the submodule spanned by graphs $\Gamma$ such that $\# \text{edges} - \# \text{vertices} \leq s$.
  The splitting part of the differential strictly decreases the filtration, so only $d_{R}$ and $d_{\mathrm{contr}}$ remain on the first page of the associated spectral sequences.
  
  One can then filter by the number of edges.
  On the first page of the spectral sequence associated to this new filtration, there is only the internal differential $d_{R}$.
  Thus on the second page, the vertices are labeled by $H^{*}(R) = H^{*}(M)$. The contracting part of the differential decreases the new filtration by exactly one, and so on the second page we see all of $d_{\mathrm{contr}}$.
  
  We can now adapt the proof of~\cite[Proposition 3.4]{Willwacher2014} to show that on the part of the complex with bivalent vertices, only the circular graphs contribute to the cohomology (we work dually so we consider a quotient instead of an ideal, but the idea is the same).
  To adapt the proof, one must see the labels of positive degree as formally adding one to the valence of the vertex, thus ``breaking'' a line of bivalent vertices.
  These labels break the symmetry (recall the coinvariants in the definition of the twisting) that allow cohomology classes to be produced.
\end{proof}

\begin{coro}
  \label{coro:zphi-factors}
  The morphism $\Zphi : \fGC_{R} \to \R$ factors through $\fGC'_{R}$ in the homotopy category of \abbr{CDGA}s.
\end{coro}

\begin{proof}
  Let us show that $\Zphi$ is homotopic to zero when restricted to the ideal defining $\fGC'_{R} = \fGC_{R}^{0} / I$ as a quotient of $\fGC_{R}$.
  Up to rescaling $\varepsilon_{A}$ by a real coefficient, we may assume that $\varepsilon_{A} \rho(-)$ and $\int_{M} \sigma(-)$ are homotopic, which induces a homotopy (by derivations) on the sub-\abbr{CDGA} of graphs with no edges.
  Hence $\Zphi$ is homotopic to zero when restricted to the ideal defining $\fGC_{R}^{0}$ from $\fGC_{R}$.
  Moreover the map $\Zphi$ vanishes on graphs with univalent vertices by Corollary~\ref{cor.weak-vanishing}.
  The degree of a circular graph with $k$ vertices is $-k < 0$ (recall that all the labels are $1_{R}$ in a circular graph), but $\Zphi$ vanishes on graphs of nonzero degree.
  Hence $\Zphi$ vanishes on the connected graphs appearing in the definition of $\fLoop_{R}$.
  Therefore, in the homotopy category of \abbr{CDGA}s, $\Zphi$ factors through the homotopy cofiber of the inclusion $\fLoop_{R} \to \fGC_{R}^{0}$, which is quasi-isomorphic to $\fGC_{R}'$ by Proposition~\ref{prop.fgcr-cofib}.
\end{proof}

The statement of the corollary is not concrete, as the ``factorization'' could go through a zigzag of maps.
However, the \abbr{CDGA}s $\fGC_{R}$ and $\fGC'_{R}$ are both cofibrant (see Lemma~\ref{lem:fgc-cofibrant} for $\fGC_{R}$, whose proof can easily be adapted to $\fGC_{R}'$).
Recall from Section~\ref{sec:dg-modules-cdgas} the following definition of homotopy.
Let $\pi : \fGC_{R} \to \fGC'_{R}$ be the quotient map.
Recall that $\APL^{*}(\Delta^{1}) = S(t,dt)$ is a path object for the \abbr{CDGA} $\R$, and $\ev_{0}, \ev_{1} : \APL^{*}(\Delta^{1}) \to \R$ are evaluation at $t = 0$ and $t = 1$.
There exists some morphism $\Zphi' : \fGC'_{R} \to \R$ and some homotopy $h : \fGC_{R} \to \APL^{*}(\Delta^{1})$ such that the following diagram commutes:
\[ \begin{tikzcd}[row sep = small, column sep = large]
    {} & \fGC_{R} \ar[dl, "\Zphi"{swap}] \ar[dr, "\Zphi' \pi"] \ar[d, "h"] \\
    \R & \APL^{*}(\Delta^{1}) \ar[l, "\sim"{near start}, "\ev_{1}" {swap, near start}] \ar[r, "\ev_{0}"{near start}, "\sim"{swap, near start}] & \R
  \end{tikzcd} \]

\begin{defin}
  \label{def:graphs-prime}
  Let $\APL^{*}(\Delta^{1})_{h}$ be the $\fGC_{R}$-module induced by $h$, and let
  \[ \Graphs'_{R}(U) = \APL^{*}(\Delta^{1})_{h} \otimes_{\fGC_{R}} \Tw \Gra_{R}(U). \]
\end{defin}

\begin{defin}
  \label{def:graphs-eps}
  Let $\Ze : \fGC_{R} \to \R$ be the algebra morphism that sends a graph $\gamma$ with a single vertex labeled by $x \in R$ to $\varepsilon_{A}(\rho(x))$, and that vanishes on all the other connected graphs.
  Let $\R_{\varepsilon}$ be the one-dimensional $\fGC_{R}$-module induced by $\Ze$, and let
  \[ \Graphs^{\varepsilon}_{R}(U) = \R_{\varepsilon} \otimes_{\fGC_{R}} \Tw \Gra_{R}(U). \]
\end{defin}

Explicitly, in $\Graphs^{\varepsilon}_{R}$, all internal components with at least two vertices are identified with zero, whereas an internal component with a single vertex labeled by $x \in R$ is identified with the number $\varepsilon_{A}(\rho(x))$.

\begin{lem}
  The morphism $\Zphi' \pi$ is equal to $\Ze$.
\end{lem}

\begin{proof}
  This is a rephrasing of Proposition~\ref{prop.strong-vanishing}.
  Using the same degree counting argument, all the connected graphs with more than one vertex in $\fGC'_{R}$ are of positive degree.
  Since $\R$ is concentrated in degree zero, $\Zphi' \pi$ must vanish on these graphs, just like $\Ze$.
  Moreover the morphism $\pi : \fGC_{R} \to \fGC'_{R} = \fGC^{0}_{R} / I$ factors through $\fGC^{0}_{R}$, where graphs $\gamma$ with a single vertex are already identified with the numbers $\Ze(\gamma)$.
\end{proof}

\begin{prop}
  \label{prop.zig-gra-e-phi}
  For each finite set $U$, we have a zigzag of quasi-isomorphisms of \abbr{CDGA}s:
  \[ \Graphs^{\varepsilon}_{R}(U) \xleftarrow{\sim} \Graphs'_{R}(U) \xrightarrow{\sim} \Graphs^{\varphi}_{R}(U). \]

  If $\chi(M) = 0$, then $\Graphs'_{R}$ and $\Graphs^{\varepsilon}_{R}$ are right Hopf $\Graphs_{n}$-comodules, and the zigzag defines a zigzag of Hopf right comodule morphisms.
\end{prop}

\begin{proof}
  We have a commutative diagram:
  \[ \begin{tikzcd}[column sep = small, row sep = small]
      \Graphs^{\varepsilon}_{R}(U) & \Graphs'_{R}(U) \lar{} \rar{} & \Graphs^{\varphi}_{R}(U) \\
      \Tw\Gra_{R}(U) \otimes_{\fGC_{R}} \R_{\varepsilon} \uar{=}
      & \Tw\Gra_{R}(U) \otimes_{\fGC_{R}} \APL^{*}(\Delta^{1})_{h} \uar{=} \lar[swap]{1 \otimes \ev_{1}} \rar{1 \otimes \ev_{0}}
      & \Tw\Gra_{R}(U) \otimes_{\fGC_{R}} \R_{\varphi} \uar{=}
    \end{tikzcd} \]

  The $\fGC_{R}$-module $\Tw\Gra_{R}(U)$ is cofibrant.
  Indeed, it is quasi-free, because $\Tw \Gra_{R}(U)$ is freely generated as a graded $\fGC_{R}$-module by reduced graphs.
  Moreover, we can adapt the proof of Lemma~\ref{lem:fgc-cofibrant} to filter the space of generators in an appropriate manner and show that $\Tw\Gra_{R}(U)$ is cofibrant.

  Therefore the functor $\Tw \Gra_{R}(U) \otimes_{\fGC_{R}} (-)$ preserves quasi-isomorphisms.
  The two evaluation maps $\ev_{0}, \ev_{1} : \APL^{*}(\Delta^{1}) \to \R$ are quasi-isomorphisms.
  It follows that all the maps in the diagram are quasi-isomorphisms.

  If $\chi(M) = 0$, the proof that $\Graphs'_{R}$ and $\Graphs^{\varepsilon}_{R}$ assemble to $\Graphs_{n}$-comodules is identical to the proof for $\Graphs^{\varphi}_{R}$ (see Proposition~\ref{prop.mph-1}).
  It is also clear that the two zigzags define morphisms of comodules: in $\Graphs_{n}$, as all internal components are identified with zero anyway.
\end{proof}

\begin{prop}
  \label{prop.mph-2}
  The \abbr{CDGA} morphisms $\rho'_{*} : \Gra_{R}(U) \to \GG{A}(U)$ extend to \abbr{CDGA} morphisms $\rho_{*} : \Graphs^{\varepsilon}_{R}(U) \to \GG{A}(U)$ by sending all reduced graphs containing internal vertices to zero.
  If $\chi(M) = 0$ this extension defines a Hopf right comodule morphism.
\end{prop}
 
\begin{proof}
  The submodule of reduced graphs containing internal vertices is a multiplicative ideal and a cooperadic coideal, so all we are left to prove is that $\rho_{*}$ is compatible with differentials.
  Since $\rho'_{*}$ was a chain map, we must only prove that if $\Gamma$ is a reduced graph with internal vertices, then $\rho_{*}(d \Gamma) = 0$.

  If a summand of $d\Gamma$ still contains an internal vertex, then it is mapped to zero by definition of $\rho_{*}$.
  So we need to look for the summands of the differential that can remove all internal vertices at once.

  The differential of $R$ leaves the number of internal vertices constant, therefore if $\Gamma$ already had an internal vertex, so do all the summands of $d_{R}\Gamma$.
  The contracting part $d_{\mathrm{contr}}$ of the differential decreases the number of internal vertices by exactly one, so let us assume that $\Gamma$ has exactly one internal vertex.
  This vertex is at least univalent, as we consider reduced graphs.
  Then there are several cases to consider, depending of the valence of the internal vertex:
  \begin{itemize}
  \item if it is univalent, then the argument of Lemma~\ref{lem.fgczero} shows that contracting the incident edge cancels with the splitting part of the differential;
  \item if it is bivalent, the contracting part has two summands, and both cancel by the symmetry relation $\iota_{u}(a) \omega_{uv} = \iota_{v}(a) \omega_{uv}$ in Equation~\eqref{eq:9};
  \item if it is at least trivalent, then we can use the symmetry relation $\iota_{u}(a) \omega_{uv} = \iota_{v}(a) \omega_{uv}$ to push all the labels on a single vertex, and we see that the sum of graphs that appear is obtained by the Arnold relation (see Figure~\ref{fig.diff-gra-n} for an example in the case of $\Graphs_{n} \to \enV$).
  \end{itemize}

  Finally, the splitting part of the differential leaves the number of internal vertices constant, unless it splits off a whole connected component with only internal vertices, in which case the component is evaluated using the partition function $\Ze$.
  If that connected component consists of a single internal vertex, then we saw in the previous item that splitting the edge connecting this univalent vertex to the rest of the graph cancels with the contraction of that edge.
  Otherwise, if the graph has more than one vertex, then by definition $\Ze$ vanishes on that graph.
\end{proof}

\subsection{The morphisms are quasi-isomorphisms}
\label{sec.morphisms-are-quasi}

In this section we prove that the morphisms constructed in Proposition~\ref{prop.mph-1} and Proposition~\ref{prop.mph-2} are quasi-isomorphisms, completing the proof of Theorem~\ref{thm.A}.

Let us recall our hypotheses and constructions.
Let $M$ be a simply connected closed smooth manifold of dimension at least $4$.
We endow $M$ with a semi-algebraic structure (Section~\ref{sec:semi-algebraic-sets}) and we consider the \abbr{CDGA} $\OmPA^{*}(M)$ of \abbr{PA} forms on $M$, which is a model for the real homotopy type of $M$.
Recall that we fix a zigzag of quasi-isomorphisms of \abbr{CDGA}s $A \xleftarrow{\rho} R \xrightarrow{\sigma} \OmPA^{*}(M)$, where $A$ is a Poincaré duality \abbr{CDGA} (Theorem~\ref{thm.lambrechts-stanley}), and $\sigma$ factors through the quasi-isomorphic sub-\abbr{CDGA} of trivial forms.

Recall that $\varphi \in \OmPA^{n-1}(\FM_{M}(\underline{2}))$ is an (anti-)symmetric trivial form on the compactification of the configuration space of two points in $M$, whose restriction to the sphere bundle $\partial \FM_{M}(\underline{2})$ is a global angular form, and whose differential $d\varphi$ is a representative of the diagonal class of $M$ (Proposition~\ref{def.propagator}).
Recall that we defined the graph complex $\Graphs_{R}^{\varphi}(U)$ using reduced labeled graphs with internal and external vertices (Definition~\ref{def:graphs-phi}) and a partition function built from $\varphi$ (Definition~\ref{def:part-func-phi}).
We also defined the variants $\Graphs_{R}^{\varepsilon}$ and $\Graphs'_{R}$ (Definitions~\ref{def:graphs-prime} and~\ref{def:graphs-eps}).

\begin{thrm}[{Precise version of Theorem~\ref{thm.A}}]
  \label{thm.Abis}
  Let $M$ be a simply connected closed smooth manifold of dimension at least $4$.
  Using the notation recalled above, the following zigzag, where the maps were constructed in Proposition~\ref{prop.mph-1}, Proposition~\ref{prop.zig-gra-e-phi}, and Proposition~\ref{prop.mph-2}, is a zigzag of quasi-isomorphisms of $\Z$-graded \abbr{CDGA}s for all finite sets $U$:
  \[ \GG{A}(U) \xleftarrow{\sim} \Graphs^{\varepsilon}_{R}(U) \xleftarrow{\sim} \Graphs'_{R}(U) \xrightarrow{\sim} \Graphs^{\varphi}_{R}(U) \xrightarrow{\sim} \OmPA^{*}(\FM_{M}(U)).\]

  If $\chi(M) = 0$, then the left-pointing maps form a quasi-isomorphism of Hopf right comodules:
  \[ (\GG{A}, \enV) \xleftarrow\sim (\Graphs^{\varepsilon}_{R}, \Graphs_{n}) \xleftarrow\sim (\Graphs'_{R}, \Graphs_{n}). \]

  If moreover $M$ is framed, then the right-pointing maps also form a quasi-isomorphism of Hopf right comodules:
  \[ (\Graphs'_{R}, \Graphs_{n}) \xrightarrow\sim (\Graphs^{\varphi}_{R}, \Graphs_{n}) \xrightarrow\sim (\OmPA^{*}(\FM_{M}), \OmPA^{*}(\FM_{n})). \]
\end{thrm}

The rest of the section is dedicated to the proof of this theorem.
Let us give a roadmap of this proof.
We first prove that $\Graphs^{\varepsilon}_{R}(U) \to \GG{A}(U)$ is a quasi-isomorphism by an inductive argument on $\#U$ (Proposition~\ref{prop.graph-A-G-A}).
This involves setting up a spectral sequence so that we can reduce the argument to connected graphs.
Then we use explicit homotopies in order to show that both complexes have cohomology of the same dimension, and we show that the morphism is surjective on cohomology by describing a section by explicit arguments.
Then we prove that $\Graphs^{\varphi}_{R}(U) \to \OmPA^{*}(\FM_{M}(U))$ is surjective on cohomology explicitly (Proposition~\ref{prop.surj-cohom}).
Since we know that $\GG{A}(U)$ and $\FM_{M}(U)$ have the same cohomology by the theorem of Lambrechts--Stanley~\cite[Theorem~10.1]{LambrechtsStanley2008a}, this completes the proof that all the maps are quasi-isomorphisms.
Compatibility with the various comodules structures was already shown in Section~\ref{sec.label-graphs-stat}.

\begin{lem}
  \label{lemma.graph-r-graph-A}
  The morphisms $\Graphs^{\varepsilon}_{R}(U) \to \GG{A}(U)$ factor through quasi-isomorphisms $\Graphs^{\varepsilon}_{R}(U) \to \Graphs^{\varepsilon}_{A}(U)$, where $\Graphs^{\varepsilon}_{A}(U)$ is the \abbr{CDGA} obtained by modding graphs with a label in $\ker(\rho : R \to A)$ in $\Graphs^{\varepsilon}_{R}(U)$.
\end{lem}

\begin{proof}
  The morphism $\Graphs^{\varepsilon}_{R} \to \Graphs^{\varepsilon}_{A}$ simply applies the surjective map $\rho : R \to A$ to all the labels.
  Hence $\Graphs^{\varepsilon}_{R} \to \GG{A}$ factors through the quotient.

  We can consider the spectral sequences associated to the filtrations of both $\Graphs^{\varepsilon}_{R}$ and $\Graphs^{\varepsilon}_{A}$ by the number of edges, and we obtain a morphism $\EE^{0} \Graphs^{\varepsilon}_{R} \to \EE^{0} \Graphs^{\varepsilon}_{A}$.
  On both $\EE^{0}$ pages, only the internal differentials coming from $R$ and $A$ remain.
  The chain map $R \to A$ is a quasi-isomorphism; hence we obtain an isomorphism on the $\EE^{1}$ page.
  By standard spectral sequence arguments, it follows that $\Graphs^{\varepsilon}_{R} \to \Graphs^{\varepsilon}_{A}$ is a quasi-isomorphism.
\end{proof}

The \abbr{CDGA} $\Graphs^{\varepsilon}_{A}(U)$ has the same graphical description as the \abbr{CDGA} $\Graphs^{\varepsilon}_{R}(U)$, except that now vertices are labeled by elements of $A$.
An internal component with a single vertex labeled by $a \in A$ is identified with $\varepsilon(a)$, and an internal component with more than one vertex is identified with zero.

\begin{prop}
  \label{prop.graph-A-G-A}
  The morphism $\Graphs^{\varepsilon}_{A} \to \GG{A}$ is a quasi-isomorphism.
\end{prop}

Before starting to prove this proposition, let us outline the different steps.
We filter our complex in such a way that on the $\EE^{0}$ page, only the contracting part of the differential remains (such a technique was already used in the proof of Proposition~\ref{prop.fgcr-cofib}).
Using a splitting result, we can focus on connected graphs.
Finally, we use a ``trick'' (Figure~\ref{fig.move-labels}) for moving labels around in a connected component, reducing ourselves to the case where only one vertex is labeled.
We then get a chain map $A \otimes \Graphs_{n} \to A \otimes \enV(U)$, which is a quasi-isomorphism thanks to the formality theorem.

Let us start with the first part of the outlined program, removing the splitting part of the differential from the picture.
We now define an increasing filtration on $\Graphs^{\varepsilon}_{A}$.
The submodule $F_{s} \Graphs^{\varepsilon}_{A}$ is spanned by reduced graphs such that $\#\text{edges} - \#\text{vertices} \le s$.

\begin{lem}
  \label{lem:abcde}
  The above submodules define a filtration of $\Graphs_{A}^{\varepsilon}$ by subcomplexes, satisfying $F_{-\#U-1} \Graphs^{\varepsilon}_{A}(U) = 0$ for each finite set $U$.
  The $\EE^{0}$ page of the spectral sequence associated to this filtration is isomorphic as a module to $\Graphs^{\varepsilon}_{A}$.
  Under this isomorphism the differential $d^{0}$ is equal to $d_{A} + d'_{\mathrm{contr}}$, where $d_{A}$ is the internal differential coming from $A$ and $d'_{\mathrm{contr}}$ is the part of the differential that contracts all edges but edges connected to a univalent internal vertex.
\end{lem}

\begin{proof}
  Let $\Gamma$ be an internally connected (Definition~\ref{def:internally}) reduced graph.
  If $\Gamma \in \Graphs^{\varepsilon}_{A}(U)$ is the graph with no edges and no internal vertices, then it lives in filtration level $-\#U$. Adding edges can only increase the filtration. Since we consider reduced graphs (i.e.\ no internal components), each time we add an internal vertex (decreasing the filtration) we must add at least one edge (bringing it back up). By induction on the number of internal vertices, each graph is of filtration at least $-\#U$.

  Let us now prove that the differential preserves the filtration and check which parts remain on the associated graded complex.
  The internal differential $d_{A}$ does not change either the number of edges nor the number of vertices and so keeps the filtration constant.
  The contracting part $d_{\mathrm{contr}}$ of the differential decreases both by exactly one, and so keeps the filtration constant too.

  The splitting part $d_{\mathrm{split}}$ of the differential removes one edge.
  If the resulting graph is still connected, then nothing else changes and the filtration is decreased exactly by $1$.
  If the resulting graph is not connected, then we get an internal component $\gamma$ which was connected to the rest of the graph by a single edge, and was then split off and identified with a number in the process.
  If $\gamma$ has a single vertex labeled by $a$ (i.e.\ we split an edge connected to a univalent vertex), then this number is $\varepsilon(a)$, and the filtration is kept constant.
  Otherwise, the summand is zero (and so the filtration is obviously preserved).
  
  In all cases, the differential preserves the filtration, and so we get a filtered chain complex.
  On the associated graded complex, the only remaining parts of the differential are $d_{A}$, $d_{\mathrm{contr}}$, and the part that splits off edges connected to univalent vertices.
  But by the proof of Proposition~\ref{prop.mph-2} this last part cancels out with the part that contracts these edges connected to univalent vertices.
\end{proof}

The symmetric algebra $S(\omega_{uv})_{u \neq v \in U}$ has a weight grading by the word-length on the generators $\omega_{uv}$.
This induces a weight grading on $\enV(U)$, because the ideal defining the relations is compatible with the weight grading.
This grading in turn induces an increasing filtration $F'_{s} \GG{A}$ on $\GG{A}$ (the extra differential strictly decreases the weight).
Define a shifted filtration on $\GG{A}$ by:
\[ F_{s} \GG{A}(U) \coloneqq F'_{s + \# U} \GG{A}(U). \]

\begin{lem}
  \label{lem:fghij}
  The $\EE^{0}$ page of the spectral sequence associated to $F_{*} \GG{A}$ is isomorphic as a module to $\GG{A}$.
  Under this isomorphism the $d^{0}$ differential is just the internal differential of $A$.
  \qed
\end{lem}

\begin{lem}
  \label{lem:klmno}
  The morphism $\Graphs^{\varepsilon}_{A} \to \GG{A}$ preserves the filtration and induces a chain map $\EE^{0} \Graphs^{\varepsilon}_{A}(U) \to \EE^{0} \GG{A}(U)$, for each $U$.
  It maps reduced graphs with internal vertices to zero, an edge $e_{uv}$ between external vertices to $\omega_{uv}$, and a label $a$ of an external vertex $u$ to $\iota_{u}(a)$.
\end{lem}

\begin{proof}
  The morphism $\Graphs^{\varepsilon}_{A}(U) \to \GG{A}(U)$ preserves the filtration by construction.
  If a graph has internal vertices, then its image in $\GG{A}(U)$ is of strictly lower filtration unless the graph is a forest (i.e.\ a product of trees).
  But trees have leaves, therefore by Corollary~\ref{cor.weak-vanishing} and the formula defining $\Graphs^{\varepsilon}_{A} \to \GG{A}$ they are mapped to zero in $\GG{A}(U)$ anyway.
  It is clear that the rest of the morphism preserves filtrations exactly, and so is given on the associated graded complex as stated in the lemma.
\end{proof}

We now use arguments similar to~\cite[Lemma~8.3]{LambrechtsVolic2014}.
For a partition $\pi$ of $U$, define the submodule $\Graphs^{\varepsilon}_{A}\langle \pi \rangle \subset \EE^{0} \Graphs^{\varepsilon}_{A}(U)$ spanned by reduced graphs $\Gamma$ such that the partition of $U$ induced by the connected components of $\Gamma$ is exactly $\pi$.
In particular let $\Graphs^{\varepsilon}_{A}\langle \{U\} \rangle$ be the submodule of connected graphs, where $\{U\}$ is the indiscrete partition of $U$ consisting of a single element.

\begin{lem}
  For each partition $\pi$ of $U$, $\Graphs^{\varepsilon}_{A}\langle\pi\rangle$ is a subcomplex of $\EE^{0} \Graphs^{\varepsilon}_{A}(U)$, and $\EE^{0} \Graphs^{\varepsilon}_{A}(U)$ splits as the sum over all partitions $\pi$:
  \[ \EE^{0} \Graphs^{\varepsilon}_{A}(U) = \bigoplus_{\pi} \bigotimes_{V \in \pi} \Graphs^{\varepsilon}_{A}\langle \{V\} \rangle.
  \]
\end{lem}

\begin{proof}
  Since there is no longer any part of the differential that can split off connected components in $\EE^{0} \Graphs^{\varepsilon}_{A}$, it is clear that $\Graphs^{\varepsilon}_{A}\langle \{U\} \rangle$ is a subcomplex.
  The splitting result is immediate.
\end{proof}

The complex $\EE^{0} \GG{A}(U)$ splits in a similar fashion.
For a monomial in $S(\omega_{uv})_{u \neq v \in U}$, say that $u$ and $v$ are ``connected'' if the term $\omega_{uv}$ appears in the monomial.
Consider the equivalence relation generated by ``$u$ and $v$ are connected''.
The monomial induces in this way a partition $\pi$ of $U$, and this definition factors through the quotient defining $\enV(U)$ (draw a picture of the $3$-term relation).
Finally, for a given monomial in $\GG{A}(U)$, the induced partition of $U$ is still well-defined.

Thus for a given partition $\pi$ of $U$, we can define $\enV \langle \pi \rangle$ and $\GG{A}\langle \pi \rangle$ to be the submodules of $\enV(U)$ and $\EE^{0} \GG{A}(U)$ spanned by monomials inducing the partition $\pi$.
It is a standard fact that $\enV \langle \{U\} \rangle = \Lie_{n}^{\vee}(U)$, see~\cite{Sinh2007}.
The proof of the following lemma is similar to the proof of the previous lemma:

\begin{lem}
  \label{lemma.partition}
  For each partition $\pi$ of $U$, $\GG{A} \langle \pi \rangle$ is a subcomplex of $\EE^{0} \GG{A}(U)$, and $\EE^{0} \GG{A}(U)$ splits as the sum over all partitions $\pi$ of $U$:
  \[ \EE^{0} \GG{A}(U) = \bigoplus_{\pi} \bigotimes_{V \in \pi} \GG{A}\langle \{V\} \rangle.
    \qed \]
\end{lem}

\begin{lem}
  The map $\EE^{0} \Graphs^{\varepsilon}_{A}(U) \to \EE^{0} \GG{A}(U)$ preserves the splitting.
  \qed
\end{lem}

We can now focus on connected graphs to prove Proposition~\ref{prop.graph-A-G-A}.

\begin{lem}
  The complex $\GG{A} \langle \{U\} \rangle$ is isomorphic to $A \otimes \enV \langle \{U\} \rangle$.
\end{lem}

\begin{proof}
  We define explicit isomorphisms in both directions.

  Define $A^{\otimes U} \otimes \enV \langle \{U\} \rangle \to A \otimes \enV \langle \{U\} \rangle$ using the multiplication of $A$.
  This constructions induces a map on the quotient $\EE^{0} \GG{A}(U) \to A \otimes \enV \langle \{U\} \rangle$, which restricts to a map $\GG{A} \langle \{U\} \rangle \to \enV \langle \{U\} \rangle$.
  Since $d_{A}$ is a derivation, this is a chain map.

  Conversely, define $A \otimes \enV \langle \{U\} \rangle \to A^{\otimes U} \otimes \enV \langle \{U\} \rangle$ by $a \otimes x \mapsto \iota_{u}(a) \otimes x$ for some fixed $u \in U$ (it does not matter which one since $x \in \enV \langle \{U\} \rangle$ is ``connected'').
  This construction gives a map $A \otimes \enV \langle \{U\} \rangle \to \GG{A} \langle \{U\} \rangle$, and it is straightforward to check that this map is the inverse isomorphism of the previous map.
\end{proof}

We have a commutative diagram of complexes:
\[ \begin{tikzcd}
    \Graphs^{\varepsilon}_{A} \langle \{U\} \rangle \rar \dar & A \otimes \Graphs'_{n} \langle \{U\} \rangle \dar{\sim} \\
    \GG{A} \langle \{U\} \rangle \rar{\cong} & A \otimes \enV \langle \{U\} \rangle
  \end{tikzcd} \]

Here $\Graphs'_{n}(U)$ is defined similarly to $\Graphs_{n}(U)$ except that multiple edges are allowed.
It is known that the quotient map $\Graphs'_{n}(U) \to \enV(U)$ (which factors through $\Graphs_{n}(U)$) is a quasi-isomorphism~\cite[Proposition 3.9]{Willwacher2014}.
The subcomplex $\Graphs'_{n} \langle \{U\} \rangle$ is spanned by connected graphs.
The upper horizontal map in the diagram multiplies all the labels of a graph.

The right vertical map is the tensor product of $\operatorname{id}_{A}$ and $\Graphs_{n} \langle \{U\} \rangle \xrightarrow{\sim} \enV \langle \{U\} \rangle$ (see~\ref{sec.formality}).
The bottom row is the isomorphism of the previous lemma.

It then remains to prove that $\Graphs^{\varepsilon}_{A} \langle \{U\} \rangle \to A \otimes \Graphs'_{n} \langle \{U\} \rangle$ is a quasi-isomorphism to prove Proposition~\ref{prop.graph-A-G-A}.
If $U = \varnothing$, then $\Graphs'_{A}(\varnothing) = \R = \GG{A}(\varnothing)$ and the morphism is the identity, so there is nothing to do.
From now on we assume that $\# U \geq 1$.

\begin{lem}
  \label{lemma.graphs-a-surj}
  The morphism $\Graphs^{\varepsilon}_{A} \langle \{U\} \rangle \to A \otimes \Graphs'_{n} \langle \{U\} \rangle$ is surjective on cohomology.
\end{lem}

\begin{proof}
  Choose some $u \in U$.
  There is an explicit chain-level section of the morphism, sending $x \otimes \Gamma$ to $\Gamma_{u,x}$, the same graph with the vertex $u$ labeled by $x$ and all the other vertices labeled by $1_{R}$.
  It is a well-defined chain map, which is clearly a section of the morphism in the lemma, hence the morphism of the lemma is surjective on cohomology.
\end{proof}

We now use a proof technique similar to the proof of~\cite[Lemma 8.3]{LambrechtsVolic2014}, working by induction.
The dimension of $H^{*}(\Graphs'_{n}\langle \{U\} \rangle) = \enV \langle \{U\} \rangle = \Lie_{n}^{\vee}(U)$ is well-known:
\begin{equation}
  \label{eq.dim-graphs_n}
  \dim H^{i} ( \Graphs'_{n} \langle \{U\} \rangle ) =
  \begin{cases}
    (\#U - 1) !, & \text{if } i = (n-1)(\#U - 1); \\
    0, & \text{otherwise.}
  \end{cases}
\end{equation}

\begin{lem}
  \label{lemma.same-cohomology}
  For all sets $U$ with $\#U \geq 1$, the dimension of $H^{i}(\Graphs^{\varepsilon}_{A} \langle \{U\} \rangle)$ is the same as the dimension:
  \begin{equation*}
    \dim H^{i}(A \otimes \Graphs'_{n} \langle \{U\} \rangle) = (\#U - 1)! \cdot \dim H^{i - (n-1)(\#U-1)}(A).
  \end{equation*}
\end{lem}

The proof will be by induction on the cardinality of $U$.
Before proving this lemma, we will need two additional sub-lemmas.

\begin{lem}
  \label{lemma.case-u-1}
  The complex $\Graphs^{\varepsilon}_{A} \langle \underline{1} \rangle$ has the same cohomology as $A$.
\end{lem}

\begin{proof}
  Let $\mathcal{I}$ be the subcomplex spanned by graphs with at least one internal vertex.
  We will show that $\mathcal{I}$ is acyclic; as $\Graphs^{\varepsilon}_{A} \langle \underline{1} \rangle / \mathcal{I} \cong A$, this will prove the lemma.

  There is an explicit homotopy $h$ that shows that $\mathcal{I}$ is acyclic.
  Given a graph $\Gamma$ with a single external vertex and at least one internal vertex, define $h(\Gamma)$ to be the same graph with the external vertex replaced by an internal vertex, a new external vertex labeled by $1_{A}$, and an edge connecting the external vertex to the new internal vertex:
  
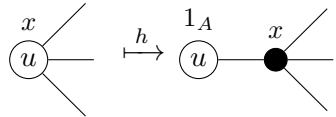
\begin{equation*}
    \begin{tikzpicture}[baseline]
      \node [extv, label = {\small $x$}] (u) {$u$};
      \node (1) [right = .6cm of u] {};
      \node (2) [above right = .8cm of u] {};
      \node (3) [below right = .8cm of u] {};
      \draw (u) -- (1); \draw (u) -- (2); \draw (u) -- (3);
    \end{tikzpicture}
    \xmapsto{\; h \;}
    \begin{tikzpicture}[baseline]
      \node [extv, label = {\small $1_{A}$}] (newu) {$u$};
      \node [intv, label = {\small $x$}] (u) [right = .6cm of newu] {};
      \node (1) [right = .6cm of u] {};
      \node (2) [above right = .8cm of u] {};
      \node (3) [below right = .8cm of u] {};
      \draw (newu) -- (u);
      \draw (u) -- (1); \draw (u) -- (2); \draw (u) -- (3);
    \end{tikzpicture}
  \end{equation*}

  The differential in $\Graphs^{\varepsilon}_{A} \langle \underline{1} \rangle$ only retains the internal differential of $A$ and the contracting part of the differential.
  Contracting the new edge in $h(\Gamma)$ gives $\Gamma$ back, and it is now straightforward to check that $d h(\Gamma) = \Gamma \pm h(d \Gamma)$.
\end{proof}

Now let $U$ be a set with at least two elements, and fix some element $u \in U$.
Let $\Graphs^{u}_{A} \langle \{U\} \rangle \subset \Graphs^{\varepsilon}_{A} \langle \{U\} \rangle$ be the subcomplex spanned by graphs $\Gamma$ such that $u$ has valence $1$, is labeled by $1_{A}$, and is connected to another external vertex.

We now get to the core of the proof of Lemma~\ref{lemma.same-cohomology}.
The idea (adapted from~\cite[Lemma 8.3]{LambrechtsVolic2014}) is to ``push'' the labels of positive degree away from the chosen vertex $u$ through a homotopy.
Roughly speaking, we use Figure~\ref{fig.move-labels} to move labels around up to homotopy.

\begin{figure}[htbp]
  \centering
  \[
    d_{\mathrm{contr}}
    \left(
    \begin{gathered} \begin{tikzpicture}[scale=.7]
      \node [unkv] (1) {};
      \node [intv, label = {$x$}, right = 1cm of 1] (i) {};
      \node [unkv, right = 1cm of i] (2) {};
      \draw (1) -- (i) -- (2);
    \end{tikzpicture} \end{gathered}
    \right)
    =
    \begin{gathered} \begin{tikzpicture}[scale=.7]
      \node [unkv, label = {$x$}] (1) {};
      \node [unkv, right = 1cm of 1] (2) {};
      \draw (1) -- (2);
    \end{tikzpicture} \end{gathered}
    -
    \begin{gathered} \begin{tikzpicture}[scale=.7]
      \node [unkv] (1) {};
      \node [unkv, right = 1cm of 1, label = {$x$}] (2) {};
      \draw (1) -- (2);
    \end{tikzpicture} \end{gathered}
  \]
  \caption{Trick for moving labels around (gray vertices are either internal or external)}
  \label{fig.move-labels}
\end{figure}

\begin{lem}
  \label{lemma.incl-q-iso}
  The inclusion $\Graphs^{u}_{A} \langle \{U\} \rangle \subset \Graphs^{\varepsilon}_{A} \langle \{U\} \rangle$ is a quasi-isomorphism.
\end{lem}

\begin{proof}
  Let $\mathcal{Q}$ be the quotient.
  We will prove that it is acyclic.
  The module $\mathcal{Q}$ further decomposes into a direct sum of modules (but the differential does not preserve the direct sum):
  \begin{itemize}
  \item The module $\mathcal{Q}_{1}$ spanned by graphs where $u$ is of valence $1$, labeled by $1_{A}$, and connected to an internal vertex;
  \item The module $\mathcal{Q}_{2}$ spanned by graphs where $u$ is of valence $\geq 2$ or has a label in $A^{> 0}$.
  \end{itemize}

  We now filter $\mathcal{Q}$ as follows.
  For $s \in \Z$, let $F_{s} \mathcal{Q}_{1}$ be the submodule of $\mathcal{Q}_{1}$ spanned by graphs with at most $s+1$ edges, and let $F_{s} \mathcal{Q}_{2}$ be the submodule spanned by graphs with at most $s$ edges.
  This filtration is preserved by the differential of $\mathcal{Q}$.

  Consider the $\EE^{0}$ page of the spectral sequence associated to this filtration.
  Then the differential $d^{0}$ is a morphism $\EE^{0} \mathcal{Q}_{1} \to \EE^{0} \mathcal{Q}_{2}$ (count the number of edges and use the crucial fact that edges connected to univalent vertices are not contractible when looking at reduced graphs).
  This differential contracts the only edge incident to $u$.
  It is an isomorphism, with an inverse similar to the homotopy defined in Lemma~\ref{lemma.case-u-1}, ``blowing up'' the point $u$ into a new edge connecting $u$ to a new internal vertex that replaces $u$.

  This shows that $(\EE^{0} \mathcal{Q}, d^{0})$ is acyclic, hence $\EE^{1} \mathcal{Q} = 0$.
  It follows that $\mathcal{Q}$ itself is acyclic.
\end{proof}

\begin{proof}[Proof of Lemma~\ref{lemma.same-cohomology}]
  The case $\# U = 0$ is obvious, and the case $\# U = 1$ of the lemma was covered in Lemma~\ref{lemma.case-u-1}.
  We now work by induction and assume the claim proved for $\# U \leq k$, for some $k \geq 1$.

  Let $U$ be of cardinality $k+1$.
  Choose some $u \in U$ and define $\Graphs^{u}_{A} \langle \{U\} \rangle$ as before.
  By Lemma~\ref{lemma.incl-q-iso} we only need to show that this complex has the right cohomology.
  It splits as:
  \begin{equation}
    \label{eq:13}
    \Graphs^{u}_{A} \langle \{U\} \rangle \cong \bigoplus_{\mathclap{v \in U \setminus \{u\}}} e_{uv} \cdot \Graphs^{\varepsilon}_{A} \langle \{U \setminus \{u\}\} \rangle,
  \end{equation}
  and therefore using the induction hypothesis:
  \begin{align*}
    \dim H^{i}(\Graphs^{u}_{A} \langle \{U\} \rangle)
    & = k \cdot \dim H^{i - (n-1)} ( \Graphs^{\varepsilon}_{A} \langle \{U \setminus \{u\}\} \rangle ) \\
    & = k! \cdot \dim H^{i - k (n-1)}(A).
      \qedhere
  \end{align*}
\end{proof}

\begin{proof}[Proof of Proposition~\ref{prop.graph-A-G-A}]
  By Lemma~\ref{lemma.graphs-a-surj}, the morphism induced by $\Graphs^{\varepsilon}_{A} \to \GG{A}$ on the $\EE^{0}$ page is surjective on cohomology.
  By Lemma~\ref{lemma.same-cohomology} and Equation~\eqref{eq.dim-graphs_n}, both $\EE^{0}$ pages have the same cohomology, and so the induced morphism is a quasi-isomorphism.
  Standard spectral arguments imply the proposition.
\end{proof}

\begin{prop}
  \label{prop.surj-cohom}
  The morphism $\omega : \Graphs'_{R}(U) \to \OmPA^{*}(\FM_{M}(U))$ is a quasi-iso\-mor\-phism.
\end{prop}

\begin{proof}
  By Equation~\eqref{eq.ea-fmm-same-cohom}, Proposition~\ref{prop.zig-gra-e-phi}, Lemma~\ref{lemma.graph-r-graph-A}, and Proposition~\ref{prop.graph-A-G-A}, both \abbr{CDGA}s have the same cohomology of finite type, so it will suffice to show that the map is surjective on cohomology to prove that it is a quasi-isomorphism.

  We work by induction.
  The case $U = \varnothing$ is immediate, as $\Graphs'_{R}(\varnothing) \xrightarrow{\sim} \Graphs^{\varphi}_{R}(\varnothing) = \OmPA^{*}(\FM_{M}(\varnothing)) = \R$ and the last map is the identity.

  Suppose that $U = \{ u \}$ is a singleton.
  Since $\rho$ is a quasi-isomorphism, for every cocycle $\alpha \in \OmPA^{*}(\FM_{M}(U)) = \OmPA^{*}(M)$ there is some cocycle $x \in R$ such that $\rho(x)$ is cohomologous to $\alpha$.
  Then the graph $\gamma_{x}$ with a single (external) vertex labeled by $x$ is a cocycle in $\Graphs'_{R}(U)$, and $\omega(\gamma_{x}) = \rho(x)$ is cohomologous to $\alpha$.
  This proves that $\Graphs'_{R}(\{u\}) \to \OmPA^{*}(M)$ is surjective on cohomology, and hence is a quasi-isomorphism.

  Now assume that $U = \{ u \} \sqcup V$, where $\# V \geq 1$, and assume that the claim is proven for sets of vertices of size at most $\# V = \#U - 1$.
  By Equation~\eqref{eq.ea-fmm-same-cohom}, we may represent any cohomology class of $\FM_{M}(U)$ by an element $z \in \GG{A}(U)$ satisfying $dz = 0$.
  Using the relations defining $\GG{A}(U)$, we may write $z$ as
  \[ z = z' + \sum_{v \in V} \omega_{uv} z_{v}, \]
  where $z' \in A \otimes \GG{A}(V)$ and $z_{v} \in \GG{A}(V)$.
  The relation $dz = 0$ is equivalent to
  \begin{align}
    dz' + \sum_{v \in V} (p_{u} \times p_{v})^{*}(\Delta_{A}) \cdot z_{v} = 0, \label{eq.some-label} \\
    \text{and } dz_{v} = 0 \text{ for all } v.
  \end{align}

  By the induction hypothesis, for all $v \in V$ there exists a cocycle $\gamma_{v} \in \Graphs'_{R}(V)$ such that $\omega(\gamma_{v})$ represents the cohomology class of the cocycle $z_{v}$ in $H^{*}(\FM_{M}(V))$, and such that $\sigma_{*}(\gamma_{v})$ is equal to $z_{v}$ up to a coboundary.

  By Equation~\eqref{eq.some-label}, the cocycle
  \[ \tilde\gamma = \sum_{v \in V} (p_{u} \times p_{v})^{*}(\Delta_{R}) \cdot \gamma_{v} \in R \otimes \Graphs'_{R}(V) \]
  is mapped to a coboundary in $A \otimes \GG{A}(V)$.
  The map $\sigma_{*} : R \otimes \Graphs'_{R}(V) \to A \otimes \GG{A}(V)$ is a quasi-isomorphism, hence $\tilde\gamma = d \tilde\gamma_{1}$ is a coboundary too.

  It follows that $z' - \sigma_{*}(\tilde\gamma_{1}) \in A \otimes \GG{A}(V)$ is a cocycle.
  Thus by the induction hypothesis there exists some $\tilde\gamma_{2} \in R \otimes \Graphs'_{R}(V)$ whose cohomology class represents the same cohomology class as $z' - \sigma_{*}(\tilde\gamma_{1})$ in $H^{*}(A \otimes \GG{A}(V)) = H^{*}(M \times \FM_{M}(V))$.

  We now let $\gamma' = -\tilde\gamma_{1} + \tilde\gamma_{2}$, hence $d\gamma' = -\tilde\gamma + 0 = - \tilde\gamma$ and $\sigma_{*}(\gamma')$ is equal to $z'$ up to a coboundary.
  By abuse of notation we still let $\gamma'$ be the image of $\gamma'$ under the obvious map $R \otimes \Graphs'_{R}(V) \to \Graphs'_{R}(U)$, $x \otimes \Gamma \mapsto \iota_{u}(x) \cdot \Gamma$.
  Then
  \[ \gamma = \gamma' + \sum_{v \in V} e_{uv} \cdot \gamma_{v} \in \Graphs'_{R}(U) \]
  is a cocycle, and $\omega(\gamma)$ represents the cohomology class of $z$ in $\OmPA^{*}(\FM_{M}(U))$.
  We have shown that the morphism $\Graphs'_{R}(U) \to \OmPA^{*}(\FM_{M}(U))$ is surjective on cohomology, and hence it is a quasi-isomorphism.
\end{proof}

\begin{proof}[Proof of Theorem~\ref{thm.Abis}]
  The zigzag of the theorem becomes, after factorizing the first map through $\Graphs^{\varepsilon}_{A}$:
  \[ \GG{A}(U) \gets \Graphs^{\varepsilon}_{A}(U) \gets \Graphs^{\varepsilon}_{R}(U) \gets \Graphs'_{R}(U) \to \Graphs^{\varphi}_{R}(U) \to \OmPA^{*}(\FM_{M}(U)) \]

  All these maps are quasi-isomorphisms by Lemma~\ref{lemma.graph-r-graph-A}, Proposition~\ref{prop.zig-gra-e-phi}, Proposition~\ref{prop.graph-A-G-A}, and Proposition\ref{prop.surj-cohom}.
  Their compatibility with the comodule structures (under the relevant hypotheses) are due to Proposition~\ref{prop.mph-1}, Proposition~\ref{prop.zig-gra-e-phi}, and Proposition~\ref{prop.mph-2}.
\end{proof}

The last thing we need to check is the following proposition, which shows that that we can choose any Poincaré duality model.

\begin{prop}
  \label{cnf.prop.ga-gb}
  If $A$ and $A'$ are two quasi-isomorphic simply connected Poincaré duality \abbr{CDGA}s, then there is a weak equivalence of symmetric collections $\GG{A} \simeq \GG{A'}$.
  If moreover $\chi(A) = 0$ then this weak equivalence is a weak equivalence of right Hopf $\enV$-comodules.
\end{prop}

\begin{proof}
  The \abbr{CDGA}s $A$ and $A'$ are quasi-isomorphic, hence there exists some cofibrant $S$ and quasi-isomorphisms $f : S \xrightarrow{\sim} A$ and $f' : S \xrightarrow{\sim} A'$.
  This yields two chain maps $\varepsilon = \varepsilon_{A} \circ f,\, \varepsilon' = \varepsilon_{A'} \circ f' : S \to \R[-n]$.
  Mimicking the proof of Proposition~\ref{lemma.zigzag-a-r-omega}, we can also find (anti-)symmetric cocycles $\Delta, \Delta' \in S \otimes S$ and such that $(f \otimes f)\Delta = \Delta_{A}$ and $(f' \otimes f')\Delta' = \Delta_{A'}$.

  We can then build symmetric collections $\Graphs_{S}^{\varepsilon, \Delta}$ and a quasi-isomorphism $f_{*} : \Graphs_{S}^{\varepsilon, \Delta} \to \GG{A}$ similarly to Section~\ref{sec.label-graphs-stat}.
  The differential of an edge $e_{uv}$ in $\Graphs_{S}^{\varepsilon,\Delta}$ is $\iota_{uv}(\Delta)$, and an isolated internal vertex labeled by $x \in S$ is identified with $\varepsilon(x)$.
  In parallel, we can build $f_{*}' : \Graphs_{S}^{\varepsilon',\Delta'} \xrightarrow{\sim} \GG{A'}$.

  If moreover $\chi(A) = 0$, then we can choose $\Delta, \, \Delta'$ such that both graph complexes become right Hopf $\Graphs_{n}$-comodules, and $f_{*}$, $f'_{*}$ are compatible with the comodule structure.
  It thus suffices to find a quasi-isomorphism $\Graphs_{S}^{\varepsilon,\Delta} \simeq \Graphs_{S}^{\varepsilon',\Delta'}$ to prove the proposition.

  We first have an isomorphism $\Graphs_{S}^{\varepsilon',\Delta'} \cong \Graphs_{S}^{\varepsilon',\Delta}$ (with the obvious notation).
  Indeed, the two cocycles $\Delta$ and $\Delta'$ are cohomologous, say $\Delta' - \Delta = d\alpha$ for some $\alpha \in S \otimes S$ of degree $n-1$.
  If we replace $\alpha$ by $(\alpha + (-1)^{n} \alpha^{21})/2$, then we can assume that $\alpha^{21} = (-1)^{n} \alpha$.
  Moreover if $\chi(A) = 0$, then we can replace $\alpha$ by $\alpha - (\mu_{S}(\alpha) \otimes 1 + (-1)^{n} 1 \otimes \mu_{S}(\alpha))/2$ to get $\mu_{S}(\alpha) = 0$.
  We then obtain an isomorphism by mapping an edge $e_{uv}$ to $e_{uv} \pm \iota_{uv}(\alpha)$ (the sign depending on the direction of the isomorphism).
  This map is compatible with differentials, with products, and with the comodule structures if $\chi(A) = 0$.

  The dg-module $S$ is cofibrant and $\R[-n]$ is fibrant (like all dg-modules).
  We can assume that $\varepsilon$ and $\varepsilon'$ induce the same map on cohomology (it suffices to rescale one map, say $\varepsilon'$, and there is an automorphism of $\Graphs_{S}^{\varepsilon',\Delta}$ which takes care of this rescaling).
  Thus the two maps $\varepsilon, \varepsilon' : S \to \R[-n]$ are homotopic, i.e.\ there exists some $h : S[1] \to \R[-n]$ such that $\varepsilon(x) - \varepsilon'(x) = h(dx)$ for all $x \in S$.
  This homotopy induces a homotopy between the two morphisms $\Ze, Z_{\varepsilon'} : \fGC_{S} \to \R$.
  Because $\Tw \Gra_{S}^{\Delta}(U)$ and $\Tw \Gra_{S}^{\Delta'}(U)$ are cofibrant as modules over $\fGC_{S}$, we obtain quasi-isomorphisms $\Graphs_{S}^{\varepsilon,\Delta} \simeq \Graphs_{S}^{\varepsilon',\Delta}$ (compare with Proposition~\ref{prop.zig-gra-e-phi}).
\end{proof}

\begin{coro}
  \label{cor:main-cor}
  Let $M$ be a smooth simply connected closed manifold and $A$ be any Poincaré duality model of $M$.
  Then $\GG{A}(\underline{k})$ is a real model for $\Conf_{k}(M)$.
\end{coro}

\begin{proof}
  The corollary follows from Theorem~\ref{thm.Abis} in the case where $\dim M \ge 4$ (together with the previous proposition to ensure that we can choose any Poincaré duality model $A$ in our constructions).
  Note that the graph complexes are, in general, nonzero even in negative degrees, but by Proposition~\ref{prop:grading} this does not change the result.
  In dimension at most $3$, the only examples of simply connected closed manifolds are $S^{2}$ and $S^{3}$.
  We address these examples in Section~\ref{sec:models-conf-2}.
\end{proof}

\begin{coro}
  \label{cor.only-depends}
  The real homotopy types of the configuration spaces of a smooth simply connected closed manifold only depends on the real homotopy type of the manifold.
\end{coro}

\begin{proof}
  When $\dim M \geq 3$, the Fadell--Neuwirth fibrations~\cite{FadellNeuwirth1962} $\Conf_{k-1}(M \setminus *) \hookrightarrow \Conf_{k}(M) \to M$ show by induction that if $M$ is simply connected, then so is $\Conf_{k}(M)$ for all $k \geq 1$.
  Hence the real model $\GG{A}(\underline{k})$ completely encodes the real homotopy type of $\Conf_{k}(M)$.
\end{proof}

\subsection{Models for configurations on the 2- and 3-spheres}
\label{sec:models-conf-2}

The degree-counting argument of Proposition~\ref{prop.strong-vanishing} does not work in dimension less than $4$, so we have to use other means to prove that the Lambrechts--Stanley \abbr{CDGA}s are models for the configuration spaces.

There are no simply connected closed manifolds of dimension $1$.
In dimension $2$, the only simply connected closed manifold is the $2$-sphere, $S^{2}$.
This manifold is a complex projective variety: $S^{2} = \mathbb{CP}^{1}$.
Hence the result of Kriz~\cite{Kriz1994} shows that $\GG{H^{*}(S^{2})}(\underline{k})$ (denoted $E(k)$ there) is a rational model for $\Conf_{k}(S^{2})$.
The $2$-sphere $S^{2}$ is studied in greater detail in Section~\ref{sec.oriented-surfaces}, where we study the action of the framed little $2$-disks operad on a framed version of $\FM_{S^{2}}$.

In dimension $3$, the only simply connected smooth closed manifold is the $3$-sphere by Perelman's proof of the Poincaré conjecture~\cite{Perelman2002,Perelman2003}.
we also the following partial result, communicated to us by Thomas Willwacher:

\begin{prop}
  \label{prop.result-s3}
  The \abbr{CDGA} $\GG{A}(\underline{k})$, where $A = H^{*}(S^{3}; \Q)$, is a rational model of $\Conf_{k}(S^{3})$ for all $k \ge 0$.
\end{prop}

\begin{proof}
  The claim is clear for $k = 0$.
  Since $S^{3}$ is a Lie  group, the Fadell--Neuwirth fibration is trivial~\cite[Theorem 4]{FadellNeuwirth1962}:

  \[ \Conf_{k}(\R^{3}) \hookrightarrow \Conf_{k+1}(S^{3}) \to S^{3} \]
  The space $\Conf_{k+1}(S^{3})$ is thus identified with $S^{3} \times \Conf_{k}(\R^{3})$, which is rationally formal with cohomology $H^{*}(S^{3}) \otimes \enV[3](\underline{k})$.
  It thus suffices to build a quasi-isomorphism between $\GG{A}(\underline{k+1})$ and $H^{*}(S_{3}) \otimes \enV(k)$.

  To simplify notation, we consider $\GG{A}(\underline{k}_{+})$ (where $\underline{k}_{+} = \{ 0, \dots, k \}$), which is obviously isomorphic to $\GG{A}(\underline{k+1})$.
  Let us denote by $\upsilon \in H^{3}(S^{3}) = A^{3}$ the volume form of $S^{3}$, and recall that the diagonal class $\Delta_{A}$ is given by $1 \otimes \upsilon - \upsilon \otimes 1$.
  We have an explicit map $f : H^{*}(S^{3}) \to \enV[3](\underline{k})$ given on generators by $f(\nu \otimes 1) = \iota_{0}(\nu)$ and $f(1 \otimes \omega_{ij}) = \omega_{ij} + \omega_{0i} - \omega_{0j}$.

  The Arnold relations show that this is a well-defined algebra morphism.
  Let us prove that $d \circ f = 0$ on the generator $\omega_{ij}$ (the vanishing on $\upsilon \otimes 1$ is clear).
  We may assume that $k = 2$ and $(i,j) = (1,2)$, and then apply $\iota_{ij}$ to get the general case.
  Then we have:
  \[ (d \circ f)(\omega_{12}) = (1 \otimes 1 \otimes \upsilon - 1 \otimes \upsilon \otimes 1) + (1 \otimes \upsilon \otimes 1 - \upsilon \otimes 1 \otimes 1) - (1 \otimes 1 \otimes \upsilon - \upsilon \otimes 1 \otimes 1) = 0 \]

  We know that both \abbr{CDGA}s have the same cohomology, so to check that $f$ is a quasi-isomorphism it suffices to check that it is surjective in cohomology.
  The cohomology $H^{*}(\GG{A}(\underline{k}_{+})) \cong H^{*}(S^{3}) \otimes \enV[3](\underline{k})$ is generated in degrees $2$ (by the $\omega_{ij}$'s) and $3$ (by the $\iota_{i}(\upsilon)$'s), so it suffices to check surjectivity in these degrees.

  In degree $3$, the cocycle $\upsilon \otimes 1$ is sent to a generator of $H^{3}(\GG{A}(\underline{k}_{+})) \cong H^{3}(S^{3}) = \Q$.
  Indeed, assume $\iota_{0}(\upsilon) = d\omega$, where $\omega$ is a linear combination of the $\omega_{ij}$ for degree reasons.
  In $d\omega$, the sum of the coefficients of each $\iota_{i}(\upsilon)$ is zero, because they all come in pairs ($d\omega_{ij} = \iota_{j}(\upsilon) - \iota_{i}(\upsilon)$).
  We want the coefficient of $\iota_{0}(\upsilon)$ to be $1$, so at least one of the other coefficient must be nonzero to compensate, hence $d\omega \neq \iota_{0}(\upsilon)$.

  It remains to prove that $H^{2}(f)$ is surjective.
  We consider the quotient map $p : \GG{A}(\underline{k}_{+}) \to \enV[3](\underline{k})$ that maps $\iota_{i}(\upsilon)$ and $\omega_{0i}$ to zero for all $1 \leq i \leq k$.
  We also consider the quotient map $q : H^{*}(S^{3}) \otimes \enV[3](\underline{k}) \to \enV[3](\underline{k})$ sending $\upsilon \otimes 1$ to zero.
  We get a morphism of short exact sequences:
  \[ \begin{tikzcd}[column sep = small, row sep = small]
      0 \rar & \ker q \rar \dar & H^{*}(S^{3}) \otimes \enV[3](\underline{k}) \rar{q} \dar{f} & \enV[3](\underline{k}) \dar{=} \rar & 0 \\
      0 \rar & \ker p \rar & \GG{A}(\underline{k}) \rar{p} & \enV[3](\underline{k}) \rar & 0
    \end{tikzcd} \]

  We consider part of the long exact sequence in cohomology induced by these short exact sequences of complexes:
  \[ \begin{tikzcd}[column sep = small, row sep = small]
      \enV[3](\underline{k})^{1} \rar \dar{=} & H^{2}(\ker q) \rar \dar{(1)} & H^{2}(H^{*}(S^{3}) \otimes \enV[3](\underline{k})) = \enV[3](\underline{k})^{2} \rar \dar{H^{2}(f)} & \enV[3](\underline{k})^{2} \dar{=} \\
      \enV[3](\underline{k})^{1} \rar & H^{2}(\ker p) \rar & H^{2}(\GG{A}(\underline{k}_{+})) \rar & \enV[3](\underline{k})^{2}
    \end{tikzcd} \]

  For degree reasons, $H^{2}(\ker q) = 0$ and so the map (1) is injective.
  By the four lemma, it follows that $H^{2}(f)$ is injective.
  Since both domain and codomain have the same finite dimension, it follows that $H^{2}(f)$ is an isomorphism.
\end{proof}

\section{Factorization homology of universal enveloping \texorpdfstring{$E_{n}$}{E\_n}-algebras}
\label{sec.fact-homol-e_n}

\subsection{Factorization homology and formality}
\label{sec.fact-homol}

The manifold $\R^{n}$ is framed.
Let $U$ be a finite set and consider the space of framed embeddings (i.e.\ such that the differential at each point preserves the given trivializations of the tangent bundles) of $U$ copies of $\R^{n}$ in itself, with the compact open topology:
\begin{equation}
  \Disk_{n}^{\mathrm{fr}}(U) :=  \Embfr(\R^{n} \times U, \R^{n}) \subset \operatorname{Map}(\R^{n} \times U, \R^{n}).
\end{equation}

Using composition of embeddings, these spaces assemble to form a topological operad $\Disk_{n}^{\mathrm{fr}}$.
This operad is weakly equivalent to the operad of little $n$-disks~\cite[Remark 2.10]{AyalaFrancis2015}, and the application that takes $f \in \Disk_{n}^{\mathrm{fr}}(U)$ to $\{ f(0 \times u) \}_{u \in U} \in \Conf_{U}(\R^{n})$ is a homotopy equivalence.

Similarly if $M$ is a framed manifold, then the spaces $\Embfr(\R^{n} \times -, M)$ assemble to form a topological right $\Disk_{n}^{\mathrm{fr}}$-module, again given by composition of embeddings.
We call it $\Disk_{M}^{\mathrm{fr}}$.
If $B$ is a $\Disk_{n}^{\mathrm{fr}}$-algebra, factorization homology is given by a derived composition product~\cite[Definition 3.2]{AyalaFrancis2015}:
\begin{equation}
  \int_{M} B := \Disk_{M}^{\mathrm{fr}} \circ_{\Disk_{n}^{\mathrm{fr}}}^{\mathbb{L}} B \coloneqq \operatorname{hocoeq} \bigl( \Disk_{M}^{\mathrm{fr}} \circ \Disk_{n}^{\mathrm{fr}} \circ B \rightrightarrows \Disk_{M}^{\mathrm{fr}} \circ B \bigr).
\end{equation}

Using~\cite[Section 2]{Turchin2013}, the pair $(\FM_{M}, \FM_{n})$ is weakly equivalent to the pair $(\Disk_{M}^{\mathrm{fr}}, \Disk_{n}^{\mathrm{fr}})$.
So if $B$ is an $\FM_{n}$-algebra, then its factorization homology is:
\begin{equation}
  \label{eq:15}
  \int_{M} B \simeq \FM_{M} \circ_{\FM_{n}}^{\mathbb{L}} B \coloneqq \operatorname{hocoeq} \bigl( \FM_{M} \circ \FM_{n} \circ B \rightrightarrows \FM_{M} \circ B \bigr).
\end{equation}

We now work in the category of chain complexes over $\R$.
We use the formality theorem (Section~\ref{sec.formality}) and the fact that weak equivalences of operads induce Quillen equivalence between categories of right modules (resp.\ categories of algebras) by~\cite[Theorems 16.A, 16.B]{Fresse2009}.
Thus, to any homotopy class $[B]$ of $E_{n}$-algebras in the category of chain complexes, there corresponds a homotopy class $[\tilde{B}]$ of $\mathtt{e}_{n}$-algebras (which is generally not easy to describe).

Using Theorem~\ref{thm.Abis}, a game of adjunctions~\cite[Theorems 15.1.A and 15.2.A]{Fresse2009} shows that:
\begin{equation}
  \label{eq:16}
  \int_{M} B \simeq \GG{A}^{\vee} \circ_{\mathtt{e}_{n}}^{\mathbb{L}} \tilde{B},
\end{equation}
where $A$ is the Poincaré duality model of $M$ mentioned in the theorem, and $\GG{A}^{\vee}$ is the right $\mathtt{e}_{n}$-module dual to $\GG{A}$ viewed as a chain complex.

\subsection{Higher enveloping algebras}
\label{sec.high-envel-algebr}

Knudsen~\cite[Theorem A]{Knudsen2016} considers a higher enveloping algebra functor $U_n$ from homotopy Lie algebras to nonunital $E_n$-algebras.
This functor generalizes the standard enveloping algebra functor from the category of Lie algebras to the category of associative algebras.

Let $n$ be at least $2$.
We can again use Kontsevich's theorem on the formality of the little disks operads to identity the category of non-unital $\mathtt{E}_{n}$-algebras with the category of $\mathtt{e}_{n}$-algebras in homotopy.
We also use that a homotopy Lie algebra is equivalent, in homotopy, to an ordinary Lie algebra.
Then we get that Knudsen's higher enveloping algebra functor is equivalent to the left adjoint of the obvious forgetful functor $\mathtt{e}_{n}\text{-Alg} \to \Lie\text{-Alg}$, which maps an $n$-Poisson algebra $B$ to its underlying shifted Lie algebra $B[1-n]$.
This model $\tilde{U}_{n} : \Lie\text{-Alg} \to \mathtt{e}_{n}\text{-Alg}$ maps a Lie algebra $\mathfrak{g}$ to the $n$-Poisson algebra given by $\tilde{U}_{n}(\mathfrak{g}) = S(\mathfrak{g}[n-1])$, with the shifted Lie bracket defined using the Leibniz rule.

Knudsen~\cite[Theorem 3.16]{Knudsen2014} also gives a way of computing factorization homology of higher enveloping algebras.
If $\mathfrak{g}$ is a Lie algebra, then so is $A \otimes \mathfrak{g}$ for any \abbr{CDGA} $A$.
Then the factorization homology of $U_{n}(\mathfrak{g})$ on $M$ is given by:
\begin{equation}
  \label{eq.thm-knudsen}
  \int_{M} U_{n}(\mathfrak{g}) \simeq \CCE(\APL^{-*}(M) \otimes \mathfrak{g})
\end{equation}
where $\CCE$ is the Chevalley--Eilenberg complex and $\APL^{-*}(M)$ is the \abbr{CDGA} of rational piecewise polynomial differential forms, with the usual grading reversed.

\begin{prop}
  \label{prop.cmp-knudsen}
  Let $A$ be a Poincaré duality \abbr{CDGA}.
  Then we have a quasi-isomorphism of chain complexes:
  \[ \GG{A}^{\vee} \circ^{\mathbb{L}}_{\mathtt{e}_{n}} S(\mathfrak{g}[n-1]) \xrightarrow{\sim} \CCE(A^{-*} \otimes \mathfrak{g}).
  \]
\end{prop}

If $A$ is a Poincaré duality model of $M$, we have $A \simeq \OmPA^{*}(M) \simeq \APL^{*}(M) \otimes_{\Q} \R$~\cite[Theorem 6.1]{HardtLambrechtsTurchinVolic2011}.
It follows that the Chevalley--Eilenberg complex of the previous proposition is weakly equivalent to the Chevalley--Eilenberg complex of Equation~\eqref{eq.thm-knudsen}.
By Equation~\eqref{eq:15}, the derived circle product over $\mathtt{e}_{n}$ computes the factorization homology of $U_{n}(\mathfrak{g})$ on $M$, and so we recover Knudsen's theorem (over the reals) for closed framed simply connected manifolds.

Let $\mathtt{I}$ be the unit of the composition product, defined by $\mathtt{I}(\underline{1}) = \R$ and $\mathtt{I}(U) = 0$ for $\# U \neq 1$.
Let $\Lambda$ be the suspension of operads, satisfying
\[ \Lambda \PP \circ (X [-1]) = (\PP \circ X) [-1] = \mathtt{I}[-1] \circ (\PP \circ X). \]
As as symmetric collection, $\Lambda \PP$ is simply given by $\Lambda \PP = \mathtt{I}[-1] \circ \PP \circ \mathtt{I}[1]$.
Recall that we let $\Lie_{n} = \Lambda^{1-n} \Lie$.
The symmetric collection
\begin{equation}
  \label{eq:20}
  \mathtt{L}_{n} := \Lie \circ \mathtt{I}[1-n] = \mathtt{I}[1-n] \circ \Lie_{n}
\end{equation}
is a $(\Lie, \Lie_{n})$-bimodule, i.e.\ a $\Lie$-algebra in the category of $\Lie_{n}$-right modules.
We have $\mathtt{L}_{n}(U) = (\Lie_{n}(U))[1-n]$.
This bimodule satisfies, for any Lie algebra $\mathfrak{g}$,
\begin{equation}
  \label{eq:19}
  \mathtt{L}_{n} \circ_{\Lie_{n}} \mathfrak{g}[n-1] \cong \mathfrak{g} \text{ as Lie algebras.}
\end{equation}

We can view the \abbr{CDGA} $A^{-*}$ as a symmetric collection concentrated in arity $0$, and as such it is a commutative algebra in the category of symmetric collections.
Thus the tensor product
\[ A^{-*} \otimes \mathtt{L}_{n} = \{ A^{-*} \otimes \mathtt{L}_{n}(k) \}_{k \geq 0} \]
becomes a $\Lie$-algebra in right $\Lie_{n}$-modules, where the right $\Lie_{n}$-module structure comes from $\mathtt{L}_{n}$ and the Lie algebra structure combines the Lie algebra structure of $\mathtt{L}_{n}$ and the \abbr{CDGA} structure of $A^{-*}$.
Its Chevalley--Eilenberg complex $\CCE(A^{-*} \otimes \mathtt{L}_{n})$ is well-defined, and by functoriality of $\CCE$, it is a right $\Lie_{n}$-module.

The proof of the following lemma is essentially found (in a different language) in~\cite[Section 2]{FelixThomas2004}.

\begin{lem}
  \label{lemma.cmp-felix-thomas}
  The right $\Lie_{n}$-modules $\GG{A}^{\vee}$ and $\CCE(A^{-*} \otimes \mathtt{L}_{n})$ are isomorphic.
\end{lem}

\begin{proof}
  We will actually define a non-degenerate pairing
  \[ \langle -, - \rangle : \GG{A}(U) \otimes \CCE(A^{-*} \otimes \mathtt{L}_{n})(U) \to \R, \]
  for each finite set $U$, compatible with differentials and the right $\Lie_{n}$-(co)module structures. As both complexes are finite-dimensional in each degree, this is sufficient to prove that they are isomorphic.

  Recall that the Chevalley--Eilenberg complex $\CCE(\mathfrak{g})$ is given by the cofree cocommutative conilpotent coalgebra $S^{c}(\mathfrak{g}[-1])$, together with a differential induced by the Koszul duality morphism $\Lambda^{-1} \Com^{\vee} \to \Lie$.
  It follows that as a module, $\CCE(A^{-*} \otimes \mathtt{L}_{n})(U)$ is given by:
  \begin{align}
    \CCE(A^{-*} \otimes \mathtt{L}_{n})(U)
    & = \bigoplus_{r \geq 0} \left( \bigoplus_{\pi \in \operatorname{Part}_{r}(U)} A^{-*} \otimes \mathtt{L}_{n}(U_{1})[-1] \otimes \dots \otimes A^{-*} \otimes \mathtt{L}_{n}(U_{r})[-1] \right)^{\Sigma_{r}}
      \nonumber
    \\
    & = \bigoplus_{r \geq 0} \left( \bigoplus_{\pi \in \operatorname{Part}_{r}(U)} (A^{n-*})^{\otimes r} \otimes \Lie_{n}(U_{1}) \otimes \dots \otimes \Lie_{n}(U_{r}) \right)^{\Sigma_{r}}
      \label{eq.split-cce}
  \end{align}
  where the sums run over all partitions $\pi = \{ U_{1} \sqcup \dots \sqcup U_{r} \}$ of $U$ and $A^{n-*} = A^{-*}[-n]$ (which is a \abbr{CDGA}, Poincaré dual to $A$).

  Fix some $r \geq 0$ and some partition $\pi = \{ U_{1} \sqcup \dots \sqcup U_{r} \}$.
  We define a first pairing:
  \begin{equation}
    \label{eq:18}
    \bigr( A^{\otimes U} \otimes \enV(U) \bigl) \otimes \bigl( (A^{n-*})^{\otimes r} \otimes \Lie_{n}(U_{1}) \otimes \dots \otimes \Lie_{n}(U_{r}) \bigr) \to \R
  \end{equation}
  as follows:
  \begin{itemize}
  \item On the $A$ factors, the pairing uses the Poincaré duality pairing $\varepsilon_{A}$.
    It is given by the following formula (where $a_{U_{i}} = \prod_{u \in U_{i}} a_{u}$):
    \[ (a_{u})_{u \in U} \otimes (a'_{1} \otimes \dots \otimes a'_{r}) \mapsto \pm \varepsilon_{A} (a_{U_{1}} \cdot a'_{1}) \dots \varepsilon_{A} (a_{U_{r}} \cdot a'_{r}), \]
  \item On the factor $\enV(U) \otimes \bigotimes_{i=1}^{r} \Lie_{n}(U_{i})$, it uses the duality pairing on $\enV(U) \otimes \mathtt{e}_{n}(U)$ (recalling that $\mathtt{e}_{n} = \Com \circ \Lie_{n}$ so that we can view $\bigotimes_{i=1}^{r} \Lie_{n}(U_{i})$ as a submodule of $\mathtt{e}_{n}(U)$).
  \end{itemize}

  The pairing in Equation~\eqref{eq:18} is the product of the two pairings we just defined.
  It is extended linearly on all of $(A^{\otimes U} \otimes \enV(U)) \otimes \CCE(A^{-*} \otimes \mathtt{L}_{n})(U)$, and it factors through the quotient defining $\GG{A}(U)$ from $A^{\otimes U} \otimes \enV(U)$.

  To check the non-degeneracy of this pairing, we use the vector subspaces $\GG{A}\langle \pi \rangle$ of Lemma~\ref{lemma.partition}, which are well-defined even though they are not preserved by the differential if we do not consider the graded space $\EE^{0} \GG{A}$.
  Fix some partition $\pi = \{ U_{1}, \dots, U_{r} \}$ of $U$, then we have an isomorphism of vector spaces:
  \[ \GG{A} \langle \pi \rangle \cong A^{\otimes r} \otimes \Lie_{n}^{\vee}(U_{1}) \otimes \dots \otimes \Lie_{n}^{\vee}(U_{r}).
  \]

  It is clear that $\GG{A}\langle \pi \rangle$ is paired with the factor corresponding to $\pi$ in Equation~\eqref{eq.split-cce}, using the Poincaré duality pairing of $A$ and the pairing between $\Lie_{n}$ and its dual; and if two elements correspond to different partitions, then their pairing is equal to zero.
  Since both $\varepsilon_{A}$ and the pairing between $\Lie_{n}$ and its dual are non-degenerate, the total pairing is non-degenerate.

  The pairing is compatible with the $\Lie_{n}$-(co)module structures, i.e.\ the following diagram commutes (a relatively easy but notationally tedious check):
  \[ \begin{tikzcd}[ampersand replacement=\&]
      \begin{matrix}
        \GG{A}(U) & \otimes & \CCE(A^{-*} \otimes \mathtt{L}_{n})(U/W) \\
        & \otimes & \Lie_{n}(W)
      \end{matrix}
      \rar{1 \otimes \circ_{W}} \dar{\circ_{W}^{\vee} \otimes 1}
      \& \GG{A}(U) \otimes \CCE(A^{-*} \otimes \mathtt{L}_{n})(U) \dar{\langle -, -\rangle} \\
      \begin{matrix}
        \GG{A}(U/W) & \otimes & \CCE(A^{-*} \otimes \mathtt{L}_{n})(U/W) \\
        \Lie_{n}^{\vee}(W) & \otimes & \Lie_{n}(W)
      \end{matrix}
      \rar{%
        \begin{smallmatrix}
          \langle -, -\rangle \\ \langle -, -\rangle_{\Lie_{n}}
        \end{smallmatrix}%
      }
      \& \R
    \end{tikzcd} \]

  Finally, we easily check, using the identity $\varepsilon_{A}(aa') = \sum_{(\Delta_{A})} \pm \varepsilon_{A}(a \Delta_{A}') \varepsilon_{A}(a' \Delta_{A}'')$ (which in turns follows from the definition of $\Delta_{A}$) that the pairing commutes with differentials (i.e.\ $\langle d(-), - \rangle = \pm \langle -, d(-) \rangle$).
\end{proof}

\begin{proof}[Proof of Proposition~\ref{prop.cmp-knudsen}]
  The operad $\mathtt{e}_{n}$ is given by the composition product $\Com \circ \Lie_{n}$ equipped with a distributive law that encodes the Leibniz rule.
  We get the following isomorphism (natural in $\mathfrak{g}$):
  \begin{align*}
    \GG{A}^{\vee} \circ_{\mathtt{e}_{n}} S(\mathfrak{g}[n-1])
    & = \GG{A}^{\vee} \circ_{\mathtt{e}_{n}} (\Com \circ \mathfrak{g}[n-1]) \\
    & \cong \GG{A}^{\vee} \circ_{\mathtt{e}_{n}} (\mathtt{e}_{n} \circ_{\Lie_{n}} \mathfrak{g}[n-1]) \\
    & \cong \GG{A}^{\vee} \circ_{\Lie_{n}} \mathfrak{g}[n-1].
  \end{align*}

  According to Lemma~\ref{lemma.cmp-felix-thomas}, the right $\Lie_{n}$-module $\GG{A}^{\vee}$ is isomorphic to $\CCE(A^{-*} \otimes \mathtt{L}_{n})$.
  The functoriality of $A^{-*} \otimes -$ and $\CCE(-)$, as well as Equation~\eqref{eq:19}, imply that we have the following isomorphism (natural in $\mathfrak{g}$):
  \begin{align*}
    \GG{A}^{\vee} \circ_{\Lie_{n}} \mathfrak{g}[n-1]
    & \cong \CCE(A^{-*} \otimes \mathtt{L}_{n}) \circ_{\Lie_{n}} \mathfrak{g}[n-1] \\
    & \cong \CCE \bigl( A^{-*} \otimes ((\mathtt{L}_{n}) \circ_{\Lie_{n}} \mathfrak{g}[n-1]) \bigr) \\
    & \cong \CCE(A^{-*} \otimes \mathfrak{g}).
  \end{align*}

  The derived circle product is computed by taking a cofibrant resolution of $S(\mathfrak{g}[n-1])$.
  Let $Q_{\mathfrak{g}} \xrightarrow{\sim} \mathfrak{g}$ be a cofibrant resolution of the Lie algebra $\mathfrak{g}$.
  Then $S(Q_{\mathfrak{g}}[n-1])$ is a cofibrant $\mathtt{e}_{n}$-algebra, and by Künneth's formula $S(Q_{\mathfrak{g}}[n-1]) \to S(\mathfrak{g}[n-1])$ is a quasi-isomorphism.
  It follows that:
  \[ \GG{A}^{\vee} \circ_{\mathtt{e}_{n}}^{\mathbb{L}} S(\mathfrak{g}[n-1]) = \GG{A}^{\vee} \circ_{\mathtt{e}_{n}} S(Q_{\mathfrak{g}}[n-1]).
  \]

  We therefore have a commutative diagram:
  \[ \begin{tikzcd}[row sep = small, column sep = small]
      \GG{A}^{\vee} \circ_{\mathtt{e}_{n}}^{\mathbb{L}} S(\mathfrak{g}[n-1])
      \dar{\cong} \rar
      & \GG{A}^{\vee} \circ_{\mathtt{e}_{n}} S(\mathfrak{g}[n-1])
      \dar{\cong} \\
      \CCE(A^{-*} \otimes Q_{\mathfrak{g}})
      \rar
      & \CCE(A^{-*} \otimes \mathfrak{g})
    \end{tikzcd} \]

  The functor $\CCE$ preserves quasi-isomorphisms of Lie algebras, hence the bottom map is a quasi-isomorphism.
  The proposition follows.
\end{proof}

\section{Outlook: The case of the 2-sphere and oriented manifolds}
\label{sec.oriented-surfaces}

Up to now, we were considering framed manifolds $M$ in order to define the action of the (unframed) Fulton--MacPherson $\FM_{n}$ on $\FM_{M}$.
When $M$ is not framed, it is not possible to coherently define insertion of infinitesimal configurations from $\FM_{n}$ into the tangent space of $M$, because we lack a coherent identification of the tangent space at every point with $\R^{n}$.
Instead, for an oriented (but not necessarily framed) manifold $M$, there exists an action of the \emph{framed} Fulton--MacPherson operad obtained by considering infinitesimal configurations together with rotations of $\mathrm{SO}(n)$ (see below for precise definitions).

In dimension $2$, the formality of $\FM_{2}$ was extended to a proof of the formality of the framed version of $\FM_{2}$ in~\cite{GiansiracusaSalvatore2010} (see also~\cite{Severa2010} for an alternative proof and~\cite{KhoroshkinWillwacher2017} for a generalization for even $n$).
We now provide a generalization of the previous work for the $2$-sphere, and we formulate a conjecture for higher dimensional closed manifolds that are not necessarily framed.

\subsection{Framed little disks and framed configurations}
\label{sec.framed-little-disks}

Following Salvatore--Wahl~\cite[Definition 2.1]{SalvatoreWahl2003}, we describe the framed little disks operad as a semi-direct product.
If $G$ is a topological group and $\PP$ is an operad in $G$-spaces, the semi-direct product $\PP \rtimes G$ is the topological operad defined by $(\PP \rtimes G)(n) = \PP(n) \times G^{n}$ and explicit formulas for the composition.
If $H$ is a commutative Hopf algebra and $\CC$ is a Hopf cooperad in $H$-comodules, then the semi-direct product $\CC \rtimes H$ is defined by formally dual formulas.

The operad $\FM_{n}$ is an operad in $\mathrm{SO}(n)$-spaces, the action rotating configurations.
Thus we can form an operad $\mathtt{f}\FM_{n} = \FM_{n} \rtimes \mathrm{SO}(n)$, the framed Fulton--MacPherson operad, weakly equivalent to the standard framed little disks operad.

Given an oriented $n$-manifold $M$, there is a corresponding right module over $\mathtt{f}\FM_{n}$, which we call $\mathtt{f}\FM_{M}$~\cite[Section 2]{Turchin2013}.
The space $\mathtt{f}\FM_{M}(U)$ is a principal $\mathrm{SO}(n)^{\times U}$-bundle over $\FM_{M}(U)$.
Since $\mathrm{SO}(n)$ is an algebraic group, $\mathtt{f}\FM_{n}$ and $\mathtt{f}\FM_{M}(U)$ are respectively an operad and a module in semi-algebraic spaces.

\subsection{Cohomology of \texorpdfstring{$\mathtt{f}\FM_n$}{fFM\_n} and potential model}
\label{sec.cohom-ffm_n-right}

The cohomology of $\mathrm{SO}(n)$ is classically given by Pontryagin and Euler classes:
\begin{align*}
  H^{*}(\mathrm{SO}(2n); \Q)
  & = S(\beta_{1}, \dots, \beta_{n-1}, \alpha_{2n-1})
  & (\deg \alpha_{2n-1} = 2n-1)
  \\
  H^{*}(\mathrm{SO}(2n+1))
  & = S(\beta_{1}, \dots, \beta_{n})
  & (\deg \beta_{i} = 4i-1)
\end{align*}

By the Künneth formula, $\mathtt{f}\enV(U) = \enV(U) \otimes H^{*}(\mathrm{SO}(n))^{\otimes U}$.
We now provide explicit formulas for the cocomposition~\cite{SalvatoreWahl2003}.
If $x \in H^{*}(\mathrm{SO}(n))$ and $u \in U$, then denote as before $\iota_{u}(x) \in H^{*}(\mathrm{SO}(n))^{\otimes U}$.
Let $W \subset U$.
If $x$ is either $\beta_{i}$ or $\alpha_{2n-1}$ in the even case, then we have:
\begin{equation}
  \label{eq.cocomp-alpha}
  \circ_{W}^{\vee}(\iota_{u}(x)) =
  \begin{cases}
    \iota_{*}(x) \otimes 1 + 1 \otimes \iota_{u}(x), & \text{if } u \in W; \\
    \iota_{u}(x) \otimes 1, & \text{otherwise.}
  \end{cases}
\end{equation}

The formula for $\circ_{W}^{\vee}(\omega_{uv})$ depends on the parity of $n$.
If $n$ is odd, then $\circ_{W}^{\vee}(\omega_{uv})$ is still given by Equation~\eqref{eq.coop-gra-n}.
Otherwise, in $\mathtt{f}\enV[2n]$ we have:
\begin{equation}
  \label{eq.cocomp-w}
  \circ_{W}^{\vee}(\omega_{uv}) = \begin{cases}
    \iota_{*}(\alpha_{2n-1}) \otimes 1 + 1 \otimes \omega_{uv}, & \text{if } u,v \in W; \\
    \omega_{[u][v]} \otimes 1, & \text{otherwise}. \\
  \end{cases}
\end{equation}

From now on, we focus on oriented surfaces.
The only simply connected compact surface is $M = S^{2}$.
We can choose $A = H^{*}(S^{2}) = S(\upsilon)/(\upsilon^{2})$ as its Poincaré duality model.
The Euler class of $A$ is $e_{A} = \chi(S^{2}) \vol_{A} = 2 \upsilon$, and the diagonal class is $\Delta_{A} = \upsilon \otimes 1 + 1 \otimes \upsilon$.
Recall that $\mu_{A}(\Delta_{A}) = e_{A}$.

\begin{defin}
  \label{def.framed-kls}
  The \textbf{framed \abbr{LS CDGA}} $\mathtt{f}\GG{A}(U)$ is given by:
  \[ \mathtt{f}\GG{A}(U) = (A^{\otimes U} \otimes \mathtt{f}\enV[2](U) / (\iota_{u} (a) \cdot \omega_{uv} = \iota_{v}(a) \cdot \omega_{uv}), d), \]
  where the differential is given by $d\omega_{uv} = \iota_{uv}(\Delta_{A})$ and $d \iota_{u}(\alpha) = \iota_{u}(e_{A})$.
\end{defin}

\begin{prop}
  \label{prop.fea-comod}
  The collection $\{ \mathtt{f}\GG{A}(U) \}_{U}$ is a Hopf right $\mathtt{f}\enV[2]$-comodule, with cocomposition given by the same formula as Equation~\eqref{eq.comodule}.
\end{prop}

\begin{proof}
  The proofs that the cocomposition is compatible with the cooperad structure of $\mathtt{f}\enV[2]$, and that this is compatible with the quotient, is the same as in the proof of Proposition~\ref{prop.comodule-eA}.
  It remains to check compatibility with differentials.

  We check this compatibility on generators.
  The internal differential of $A = H^{*}(S^{2})$ is zero, so it is easy to check that $\circ_{W}^{\vee}(d(\iota_{u}(a))) = d(\circ_{W}^{\vee}(\iota_{u}(a))) = 0$.
  Similarly, by Equation~\eqref{eq.cocomp-alpha}, checking the equality on $\alpha$ is immediate.
  As before there are several cases to check for $\omega_{uv}$.
  If $u,v \in W$, then by Equation~\eqref{eq.cocomp-w},
  \begin{align*}
    d(\circ_{W}^{\vee}(\omega_{uv}))
    & = d(\iota_{*}(\alpha) \otimes 1 + 1 \otimes \omega_{uv}) = \iota_{*}(e_{A}) \otimes 1 \\
    & = \iota_{*}(\mu_{A}(\Delta_{A})) \otimes 1 = \circ_{W}^{\vee}(d \omega_{uv}),
  \end{align*}
  and otherwise the proof is identical to the proof of Proposition~\ref{prop.comodule-eA}.
\end{proof}

\subsection{Connecting \texorpdfstring{$\mathtt{f}\GG{A}$}{G\_A} to \texorpdfstring{$\OmPA^*(\mathtt{f}\FM_{S^{2}})$}{Omega(fFM\_S2)}}
\label{sec.connecting-e_a-with}

The framed little $2$-disks operad is formal~\cite{GiansiracusaSalvatore2010,Severa2010}.
We focus on the proof of Giansiracusa--Salvatore~\cite{GiansiracusaSalvatore2010}, which goes along the same line as the proof of Kontsevich of the formality of $\FM_{n}$.
To simplify notations, let $H = H^{*}(S^{1})$, which is a Hopf algebra.
The operad $\Graphs_{2}$ is an operad in $H$-comodules, so there is a semi-direct product $\Graphs_{2} \rtimes H$.
Giansiracusa and Salvatore construct a zigzag:
\begin{equation}
  \label{eq:14}
  \mathtt{f}\enV[2] \xleftarrow{\sim} \Graphs_{2} \rtimes H \xrightarrow{\sim} \OmPA(\mathtt{f}\FM_{2}).
\end{equation}

The first map is the tensor product of $\Graphs_{2} \xrightarrow{\sim} \enV[2]$ and the identity of $H$.
The second map is given by the Kontsevich integral on $\Graphs_{2}$ and by sending the generator $\alpha \in H$ to the volume form of $\OmPA^{*}(S^{1})$ (pulled back by the relevant projection).
They check that both maps are maps of Hopf (almost) cooperads, and they use the Künneth formula to conclude that these maps are quasi-isomorphisms.

\begin{thrm}
  The Hopf right comodule $(\mathtt{f}\GG{A}, \mathtt{f}\enV[2])$, where $A = H^{*}(S^{2}; \R)$, is quasi-isomorphic to the Hopf right comodule $(\OmPA^{*}(\mathtt{f}\FM_{S^2}), \OmPA^{*}(\mathtt{f}\FM_{2}))$.
\end{thrm}
\begin{proof}
  It is now straightforward to adapt the proof of Theorem~\ref{thm.A} to this setting, reusing the proof of Giansiracusa--Salvatore~\cite{GiansiracusaSalvatore2010}.
  We build the zigzag:
  \[ \mathtt{f}\GG{A} \gets \Graphs^{\varepsilon}_{A} \rtimes H \to \OmPA^{*}(\mathtt{f}\FM_{S^2}). \]

  We simply choose $R = A = H^{*}(S^{2})$, mapping $\upsilon \in H^{2}(S^{2})$ to the volume form of $S^{2}$.
  Note that the propagator can be made completely explicit on $S^{2}$, and it can be checked that $Z_{\varphi}$ vanishes on all connected graphs with more than one vertex~\cite[Proposition 80]{CamposWillwacher2016}.
  The middle term is a Hopf right $(\Graphs_{2} \rtimes H)$-comodule built out of $\Graphs^{\varepsilon}_{A}$ and $H$, using formulas similar to the formulas defining $\Graphs_{2} \rtimes H$ out of $\Graphs_{2}$ and $H$.
  The first map is given by the tensor product of $\Graphs_{R} \to \GG{A}$ and the identity of $H$.

  The second map is given by the morphism of Proposition~\ref{prop.mph-1} on the $\Graphs^{\varepsilon}_{A}$ factor, composed with the pullback along the projection $\mathtt{f}\FM_{S^2} \to \FM_{S^2}$.
  The generator $\alpha \in H$ is sent to a pullback of a global angular form $\psi$ of the principal $\mathrm{SO}(2)$-bundle $\mathtt{f}\FM_{S^{2}}(\underline{1}) \to \FM_{S^{2}}(\underline{1}) = S^2$ induced by the orientation of $S^2$.
  This form satisfies $d\psi = \chi(S^2) \vol_{S^2}$.

  The proof of Giansiracusa--Salvatore~\cite{GiansiracusaSalvatore2010} then adapts itself to prove that these two maps are maps of Hopf right comodules.
  The Künneth formula implies that the first map is a quasi-isomorphism, and the second map induces an isomorphism on the $\EE^{2}$-page of the Serre spectral sequence associated to the bundle $\mathtt{f}\FM_{S^2} \to \FM_{S^2}$ and hence is itself a quasi-isomorphism.
\end{proof}

\begin{coro}
  The \abbr{CDGA} $\mathtt{f}\GG{H^{*}(S^2)}(\underline{k})$ of Definition~\ref{def.framed-kls} is a real model for $\Conf^{\mathrm{or}}_{k}(S^2)$, the $\mathrm{SO}(2)^{\times k}$-principal bundle over $\Conf_{k}(S^2)$ induced by the orientation of $S^2$.
\end{coro}

If $M$ is an oriented $n$-manifold with $n > 2$, Definition~\ref{def.framed-kls} readily adapts to define $\mathtt{f}\GG{H^{*}(M)}$, by setting $d\alpha$ to be the Euler class of $M$ (when $n$ is even), and $d\beta_{i}$ to be the $i$th Pontryagin class of $M$.
The proof of Proposition~\ref{prop.fea-comod} adapts easily to this new setting, and $\mathtt{f}\GG{H^{*}(M)}$ becomes a Hopf right $\mathtt{f}\enV$-comodule.

\begin{conj}
  If $M$ is a formal, simply connected, oriented closed $n$-manifold and if the framed little $n$-disks operad $\mathtt{fe}_{n}$ is formal, then the pair $(\mathtt{f}\GG{H^{*}(M)}, \mathtt{f}\enV)$ is quasi-isomorphic to the pair $(\OmPA^{*}(\mathtt{f}\FM_{M}), \OmPA^{*}(\mathtt{f}\FM_{n}))$.
\end{conj}

To directly adapt our proof for the conjecture, the difficulty would be the same as encountered by Giansiracusa--Salvatore~\cite{GiansiracusaSalvatore2010}, namely finding forms in $\OmPA^{*}(\mathtt{f}\FM_{n})$ corresponding to the generators of $H^{*}(\mathrm{SO}(n))$ and compatible with the Kontsevich integral.
It was recently proved that the framed Fulton--MacPherson is formal for even $n$ and not formal for odd $n \ge 3$ \cite{Moriya2016,KhoroshkinWillwacher2017}.
However, the proof that $\mathtt{f}\FM_{n}$ is formal for even $n \ge 4$, due to Khoroshkin and Willwacher~\cite{KhoroshkinWillwacher2017}, is much more involved than the proof of the formality of $\mathtt{fFM}_{2}$.
In particular, the zigzag of maps is not completely explicit and relies on obstruction-theoretical arguments.
It would be interesting to try and adapt the conjecture in this setting.

If $M$ itself is not formal then it is also not clear how to define Pontryagin classes in some Poincaré duality model of $M$ (the Euler class is canonically given by $\chi(A) \vol_{A}$).
Nevertheless, for any oriented manifold $M$ we get invariants of $\mathtt{fe}_{n}$-algebras by considering the functor $\mathtt{f}\GG{H^{*}(M)}^{\vee} \circ^{\mathbb{L}}_{\mathtt{fe}_{n}} (-)$.
Despite not necessarily computing factorization homology, they could prove interesting.

\paragraph{Acknowledgments}

I would like to thank several people:
my (now former) advisor Benoit Fresse for giving me the opportunity to study this topic and for numerous helpful discussions regarding the content of this paper;
Thomas Willwacher and Ricardo Campos for helpful discussions about their own model for configuration spaces ktheir explanation of propagators and partition functors, and for several helpful remarks;
Ben Knudsen for explaining the relationship between the \abbr{LS CDGA}s and the Chevalley--Eilenberg complex;
Pascal Lambrechts for several helpful discussions;
Ivo Dell'Ambrogio, Julien Ducoulombier, Matteo Felder, and Antoine Touzé for their comments;
and finally the anonymous referee, for a thorough and detailed report with numerous suggestions that greatly improved this paper.
The author was supported by ERC StG 678156--GRAPHCPX during part of the completion of this manuscript.

\printbibliography

\appendix
\small
\section{Glossary of notation}

\noindent\textbf{DG-modules and \abbr{CDGA}s}
\begin{itemize}[nosep]
\item[] $V[k] = \bigoplus_{n \in \Z} V^{n+k}$: desuspension of a dg-module (Section~\ref{sec:dg-modules-cdgas})
\item[] $(v \otimes w)^{21} \coloneqq \pm w \otimes v$ (Section~\ref{sec:dg-modules-cdgas})
\item[] $X = \sum_{(X)} X' \otimes X'' \in V \otimes W$: Sweedler's notation (Section~\ref{sec:dg-modules-cdgas})
\end{itemize}

\noindent\textbf{Cooperads and comodules}
\begin{itemize}[nosep]
\item[] $\underline{k} = \{ 1, \dots, k \}$ (Section~\ref{sec:coop-their-right})
\item[] $\circ^{\vee}_{W} : \CC(U) \to \CC(U/W) \otimes \CC(W)$: cooperadic cocomposition (Section~\ref{sec:coop-their-right})
\item[] $\circ^{\vee}_{W} : \NN(U) \to \NN(U/W) \otimes \CC(W)$: right comodule structure map (Section~\ref{sec:coop-their-right})
\end{itemize}

\noindent\textbf{Semi-algebraic sets and \abbr{PA} forms}
\begin{itemize}[nosep]
\item[] $\OmPA^{*}(-)$: \abbr{CDGA} of piecewise semi-algebraic (PA) forms (Section~\ref{sec:semi-algebraic-sets})
\item[] $p_{*}(-) = \int_{p : E \to B} (-)$: integral along the fibers of the \textsc{pa} bundle $p$ (Section~\ref{sec:semi-algebraic-sets})
\end{itemize}

\noindent\textbf{Little disks and related objects}
\begin{itemize}[nosep]
\item[] $\FM_{n}(k)$: Fulton--MacPherson compactification of $\Conf_{k}(\R^{n})$ (Section~\ref{sec.little-disks})
\item[] $\mathtt{e}_{n} \coloneqq H_{*}(\FM_{n})$, $\enV \coloneqq H^{*}(\FM_{n})$ homology and cohomology of $\FM_{n}$ (Section~\ref{sec.little-disks})
\item[] $\vol_{n-1} \in \OmPA^{n-1}(\FM_{n}(\underline{2}))$ volume form (Section~\ref{sec.little-disks})
\item[] $\FM_{M}(k)$: Fulton--MacPherson compactification of $\Conf_{k}(M)$ (Section~\ref{sec.little-disks})
\item[] $p : \partial \FM_{M}(\underline{2}) \to M$ sphere bundle of rank $n-1$ (Section~\ref{sec.little-disks})
\end{itemize}

\noindent\textbf{Poincaré duality \abbr{CDGA}s}
\begin{itemize}[nosep]
\item[] $(A,\varepsilon_{A})$: Poincaré duality \abbr{CDGA} with its orientation (Section~\ref{sec.poinc-dual-models})
\item[] $\vol_{A} \in A^{n}$: volume form (Section~\ref{sec.poinc-dual-models})
\item[] $\Delta_{A} \in (A \otimes A)^{n}$: diagonal cocycle (Section~\ref{sec.poinc-dual-models})
\item[] $\GG{A}(\underline{k})$: Lambrechts--Stanley \abbr{CDGA}s (Section~\ref{sec.kriz-lambr-stanl})
\end{itemize}

\noindent\textbf{Graph complexes for $\R^{n}$}
\begin{itemize}[nosep]
\item[] $\Gra_{n}$: graphs with only external vertices (Section~\ref{sec.formality})
\item[] $\Tw \Gra_{n}$: graphs with external and internal vertices (Section~\ref{sec.formality})
\item[] $\Graphs_{n}$: reduced graphs with external and internal vertices (Section~\ref{sec.formality})
\item[] $\Gra_{n}^{\circlearrowleft}$, $\Graphs_{n}^{\circlearrowleft}$: variants with loops and multiple edges (Section~\ref{sec.graphs-with-loops})
\item[] $\mu = e_{12}^{\vee}$: Maurer--Cartan element used to twist the graphs cooperad (Section~\ref{sec.formality})
\item[] $\omega : \Graphs_{n} \to \OmPA^{*}(\FM_{n})$: Kontsevich's integrals (Section~\ref{sec.formality})
\end{itemize}

\noindent\textbf{Graph complexes for a closed manifold $M$}
\begin{itemize}[nosep]
\item[] $\Gra_{R}$: labeled graphs with only external vertices (Section~\ref{sec.external-vertices})
\item[] $\Gra_{R}^{\circlearrowleft}$: variant with loops and multiple edges (Section~\ref{sec.external-vertices})
\item[] $\Tw\Gra_{R}$: labeled graphs with internal and external vertices (Section~\ref{sec.twisting-twgra_r})
\item[] $\varphi \in \OmPA^{n-1}(\FM_{M}(\underline{2}))$: propagator (Section~\ref{sec.propagator})
\item[] $\fGC_{R}$: full labeled graph complex (Definition~\ref{def:fgc})
\item[] $\Zphi : \fGC_{R} \to \R$: partition function (Section~\ref{sec.label-graphs-stat})
\item[] $\Graphs_{R}^{\varphi}$: reduced labeled graphs with internal and external vertices (Section~\ref{sec.reduct-graphs_r})
\item[] $\omega : \Graphs_{R}^{\varphi}$: integrals (Section~\ref{sec.reduct-graphs_r})
\item[] $\Ze : \fGC_{R} \to \R$: almost trivial partition function (Definition~\ref{def:graphs-eps})
\item[] $\Graphs_{R}^{\varepsilon}$: reduced labeled graphs with internal and external vertices (Definition~\ref{def:graphs-eps})
\end{itemize}

\noindent\textbf{Factorization homology}
\begin{itemize}[nosep]
\item[] $\Disk_{n}^{\mathrm{fr}}$: operad of framed embeddings (Section~\ref{sec.fact-homol-e_n})
\item[] $\Disk_{M}^{\mathrm{fr}}$: module of framed embeddings for a framed $M$ (Section~\ref{sec.fact-homol-e_n})
\item[] $\int_{M} A \coloneqq \Disk_{M}^{\mathrm{fr}} \circ_{\Disk_{n}^{\mathrm{fr}}}^{\mathbb{L}} A$: factorization homology (Section~\ref{sec.fact-homol-e_n})
\item[] $\CCE$: Chevalley--Eilenberg complex (Section~\ref{sec.fact-homol-e_n})
\end{itemize}

\noindent\textbf{Framed case}
\begin{itemize}[nosep]
\item[] $\mathtt{fFM}_{n} = \FM_{n} \rtimes \mathrm{SO}(n)$ framed Fulton--MacPherson operad (Section~\ref{sec.oriented-surfaces})
\item[] $\mathtt{fFM}_{M}$: framed Fulton--MacPherson compactification (Section~\ref{sec.oriented-surfaces})
\item[] $\mathtt{fG}_{A}(\underline{k})$: framed Lambrechts--Stanley \abbr{CDGA}s (Section~\ref{sec.oriented-surfaces})
\end{itemize}
\end{document}